\documentclass[12pt]{article}
\usepackage{amssymb,amsmath,amscd, amsthm,stmaryrd,verbatim, mathrsfs}
\usepackage{booktabs}
\usepackage[dvipsnames]{xcolor}
\definecolor{mydarksky}{RGB}{26,71,120}
\usepackage[colorlinks,
linkcolor=mydarksky,
anchorcolor=cyan,
citecolor=WildStrawberry
]{hyperref}
\usepackage{enumerate}
\usepackage{hyperref}
\usepackage{appendix}
\numberwithin{equation}{section}
\newtheorem{theorem}{Theorem}[section]
\newtheorem{lemma}[theorem]{Lemma}
\newtheorem{corollary}[theorem]{Corollary}
\newtheorem{proposition}[theorem]{Proposition}
\theoremstyle{definition}

\theoremstyle{remark}
\newtheorem{remark}[theorem]{Remark}

\title{}
\AtEndDocument{\bigskip{\footnotesize%
  (Z.~Zhou) \textsc{Shenzhen International Center for Mathematics, Southern University of Science and Technology, Shenzhen 518055, P. R. China} \par
  \textit{E-mail address}: \texttt{\href{mailto:zhouzhenyu14@mails.ucas.ac.cn}{ zhouzhenyu14@mails.ucas.ac.cn}} \par
}}
\begin{document}

\baselineskip=15pt

\begin{center}{\large \bf Inhomogeneous First-order Differential Operator Representations of $E_6$ and Graded  Betti Numbers} \footnote{\href{https://mathscinet.ams.org/mathscinet/msc/msc2020.html}{2020 Mathematical Subject Classification}: Primary 17B10; Secondary 17B25, 13D02, 17B55, 05A19.}
\end{center}

\begin{center}{\large Zhenyu Zhou}\end{center}

\begin{abstract}
We realize the restricted dual Verma modules of arbitrary highest weight for the complex simple Lie algebra of type $E_6$ by inhomogeneous first-order differential operators in $36$ variables. The construction is obtained by applying Xu’s functor (2016) to the full conformal oscillator representation of $\mathfrak{o}(10)$. We construct polynomial bases for all six fundamental $E_6$-modules and obtain a finite spanning set for every finite-dimensional irreducible $E_6$-module $V(\lambda)$. The degree in $16$ selected variables gives $V(\lambda)$ a grading over the symmetric algebra of the opposite abelian nilradical. We compute its Hilbert series and graded Betti numbers, and verify the Buchsbaum-Eisenbud-Horrocks conjecture for this family. The minimal graded free resolution of $V(\lambda)$ is pure precisely when $\lambda$ is a nonnegative multiple of the Weyl vector of $E_6$. Among nonzero dominant integral $\lambda$'s, the minimum total Betti number is $435204$, attained exactly at two minuscule weights $\omega_1$ and $\omega_6$.
\vspace{0.3cm}

\noindent{\it Keywords}:\hspace{0.3cm} Lie algebra of type $E_6$; differential operators; finite-dimensional irreducible modules; singular vectors;  graded Betti numbers.
\end{abstract}

\tableofcontents

\section{Introduction}

Finite-dimensional irreducible representations of finite-dimensional simple Lie algebra over $\mathbb{C}$ were classified abstractly by Cartan and Weyl in the early twentieth century.  Nevertheless, explicit formulas for root vectors in these algebras are generally difficult to obtain. In 1950, Gelfand and Tsetlin \cite{GT1, GT2} introduced bases labeled by specific tableaux for finite-dimensional irreducible representations of the general linear and orthogonal Lie algebras using successive multiplicity-free restrictions. They also gave explicit formulas for the action of a set of Lie algebra generators on these bases, now known as the Gelfand–Tsetlin formulas. Zhelobenko \cite{DPZ1,DPZ2} derived the well-known branching rules for all classical Lie algebras by realizing representations in spaces of polynomials satisfying certain systems of differential equations. He also constructed lowering operators that yield matrix element formulas for the general linear Lie algebra. Using the theory of Mickelsson algebras and the representation theory of Yangians, Molev \cite{M1,M2,M3,M4} constructed a Gelfand–Tsetlin weight basis for finite-dimensional irreducible representations of classical Lie algebras and gave explicit formulas for the matrix elements of the generators. Together with Yakimova \cite{MY}, he also derived branching rules for finite-dimensional irreducible representations of reductive Lie algebras via monomial bases of multiplicity spaces associated to restrictions to reductive subalgebras. 

The Lie algebra and group $E_6$ are prominent objects in modern mathematics and appear increasingly often in various areas of science. The $78$-dimensional simple Lie algebra of type $E_6$ can be realized via derivations and multiplication operators with zero trace on the $27$-dimensional exceptional simple Jordan algebra (see e.g., \cite{A,T}). Dickson \cite{D} (1901) first observed the existence of an $E_6$-invariant symmetric trilinear form on its $27$-dimensional fundamental irreducible module. Aschbacher \cite{As1,As2,As3,As4} used this symmetric trilinear form to study the subgroup structure of the group $E_6$ and the geometry induced on the $27$-dimensional module. Binegar and Zierau \cite{BZ} found a singular representation of $E_6$. In \cite{X1}, Xu proved that the space of homogeneous polynomial solutions of degree $m$ to the dual cubic Dickson invariant differential operator decomposes as a direct sum of $\llbracket m/2 \rrbracket+1$ explicitly determined irreducible $E_6$-submodules and that the full polynomial algebra is a free module over the polynomial algebra in the Dickson invariant generated by these solutions. Anderson and Blažek \cite{AB1,AB2,AB3} found certain Clebsch-Gordan coefficients relevant to $E_6$ unification model building.  Fernández-Núñez, Garcia-Fuertes and Perelomov \cite{FGP} used the quantum Calogero-Sutherland model corresponding to the $E_6$ root system to compute Clebsch-Gordan series for this algebra. There are also numerous other interesting applications of $E_6$ (see e.g., \cite{BK, GSA, G}).

Despite the richness of the representations of classical Lie algebras, results for $E_6$ are comparatively scarce due to its intricate structure. The motivation of this paper is to realize natural explicit representations of $E_6$, in particular to construct finite-dimensional irreducible representations for arbitrary dominant integral highest weights. Throughout this paper, we work over the base field $\mathbb{C}$. 

Let $\mathscr{Y}=\mathbb{C}[y_1,\cdots,y_{16}]$ be the polynomial algebra in $16$ variables, let $M$ be a module of the orthogonal Lie algebra $\mathfrak{o}(10)$, and let $\mathscr{G}^{E_6}$ be the simple Lie algebra of type $E_6$. Denote
\begin{eqnarray}
	\widehat{M}=\mathscr{Y}\otimes_\mathbb{C} M.
\end{eqnarray}
Xu \cite{Xu} discovered a new inhomogeneous first-order differential operator (oscillator) representation of $\mathscr{G}^{E_6}$ on $\mathscr{Y}$, which yields a fractional representation of the corresponding Lie group on the $16$-dimensional space. Moreover, a specific subalgebra of $\mathscr{G}^{E_6}$, denoted by $\mathscr{G}^{D_5}$ which is isomorphic to $\mathfrak{o}(10)$, acts on $\mathscr{Y}$ by homogeneous first-order differential operators. Using this representation and Shen's mixed product for Witt algebra in \cite{Shen}, Xu obtained a $\mathscr{G}^{E_6}$-module structure on $\widehat{M}$. Furthermore,
\begin{eqnarray}
	\overline{M}=U(\mathscr{G}^{E_6})(1\otimes M)
\end{eqnarray}
is a $\mathscr{G}^{E_6}$-submodule. For $n\in\mathbb{N}$, let $\mathscr{Y}[n]$ be the subspace of homogeneous polynomials in $\mathscr{Y}$ of degree $n$, and set
\begin{eqnarray}
	\widehat{M}[n]=\mathscr{Y}[n]\otimes_{\mathbb{C}}M,\quad \overline{M}[n]=\overline{M}\cap\widehat{M}[n]. \label{1.3}
\end{eqnarray}
We say that the elements of $\widehat{M}[n]$ has {\it $y$-degree} $n$, and call $\widehat{M}[n]$ as the {\it $y$-homogeneous} subspace of $\widehat{M}$ with $y$-degree $n$.
Then $\widehat{M}[n]$ is the tensor product of the $\mathfrak{o}(10)$-modules $\mathscr{Y}[n]$ and $M$, and $\overline{M}[n]$ is an $\mathfrak{o}(10)$-submodule of $\widehat{M}[n]$. If $M$ is an irreducible $\mathfrak{o}(10)$-module, then $\overline{M}$ is naturally an irreducible $\mathscr{G}^{E_6}$-module. This framework provides a direct polynomial extension from irreducible $\mathfrak{o}(10)$-modules to irreducible $\mathscr{G}^{E_6}$-modules. Xu also found a sufficient condition for $\widehat{M}$ to be irreducible when $M$ is a finite-dimensional irreducible $\mathfrak{o}(10)$-module. Moreover, $\overline{M}$ is a highest-weight irreducible $\mathscr{G}^{E_6}$-module whenever $M$ is a highest-weight irreducible $\mathfrak{o}(10)$-module. 

On the other hand, the author constructed in previous work \cite{Zhou} the so-called full conformal oscillator representation of the even orthogonal Lie algebras $\mathfrak{o}(2n+2)$ on the polynomial algebra $\mathscr{X}$ in $n(n+1)$ variables, write as $\{x_{i,j}\mid 1\leqslant |j|\leqslant i\leqslant n\}$, and explicitly described the corresponding finite-dimensional irreducible module for dominant integral highest weights. In this paper, we take $M=\mathscr{X}$ and construct representations of $\mathscr{G}^{E_6}$ on the polynomial algebra $\mathscr{C}$ in $36$ variables, where the action of every root vector is given explicitly as an inhomogeneous first-order differential operator (Theorem \ref{thm3.1} and Theorem \ref{thm3.2}). 

This differential operators realization also admits a geometric interpretation on the big cell of the flag variety. Let $G$ be the simply connected complex algebraic group with Lie algebra $\mathscr{G}^{E_6}$, and let $B_-$ and $N_+$ be the negative Borel subgroup and the positive unipotent subgroup determined by the chosen triangular decomposition. The flag variety $X=B_-\setminus G$ contains the open cell $U=B_-\setminus B_-N_+\cong N_+\cong \mathbb{A}^{36}$ which is a $36$-dimensional affine space. Theorem \ref{thm3.1} yields a homomorphism
\begin{align}
    \pi_{\lambda}:U(\mathscr{G}^{E_6})\longrightarrow \mathscr{D}(\mathbb{A}^{36})
\end{align}
of associative algebra, where $\mathscr{D}(\mathbb{A}^{36})$ is the Weyl algebra of polynomial coefficient differential operators. The constructed module $\mathscr{C}$ is therefore obtained by restricting the natural module $\mathscr{O}(\mathbb{A}^{36})$ along $\pi_{\lambda}$. This viewpoint is related to the realization of representations of Lie algebras by twisted differential operators on flag varieties, as in Beilinson-Bernstein localization \cite{BB}. On the big cell, the infinitesimal action takes the form of vector fields together with multiplication operators which are dependent on the weight.

A $singular$ $vector$ in a $\mathscr{G}^{E_6}$-module $V$ is a nonzero weight vector annihilated by all positive root vectors. We find that the constant $1\in\mathscr{C}$ is the unique singular vector up to scalar, and that the cyclic $\mathscr{G}^{E_6}$-submodule $V(\lambda)=U(\mathscr{G}^{E_6})(1)$ is an irreducible highest-weight $\mathscr{G}^{E_6}$-module whose highest weight $\lambda$ can be arbitrary. We prove that the restricted dual module $\mathscr{C}^{\vee}$ is isomorphic to the Verma module with the same highest weight. This result establishes a precise relationship between the polynomial realization and the algebraic theory of highest-weight modules. Moreover, we find all the polynomials in $\mathscr{C}$ which are annihilated by $E_{\alpha_i}$ for $i\in\overline{1,5}$, thus all the singular vectors in $\mathscr{C}$ when it is restricted to the Lie subalgebra $\mathscr{G}^{D_5}$ of corank $1$, showing that they span a polynomial algebra in $16$ explicit variables. Such singular vectors also have a direct geometric interpretation. Let $I$ be the left ideal generated by $\pi_{\lambda}(E_{\alpha_1}),\cdots,\pi_{\lambda}(E_{\alpha_5})$ of $\mathscr D(\mathbb{A}^{36})$. Evaluation at $1+I$ identifies $\mathrm{Hom}_{\mathscr D(\mathbb{A}^{36})}(\mathscr D(\mathbb{A}^{36})/I,\mathscr{O}(\mathbb{A}^{36}))$ with the kernel of the actions of these positive root vectors of $\mathscr{G}^{D_5}$.

Denote by $\omega_i$ the $i$-th fundamental weight of $\mathscr{G}^{E_6}$ with respect to its selected simple positive roots ($i=1,\cdots,6$). We prove that for dominant integral weights $\lambda$ and $\lambda'$ the pairwise products of elements of spanning sets of $V(\lambda)$ and $V(\lambda')$ span $V(\lambda+\lambda')$. Consequently, for any dominant integral weight $\lambda$, the module $V(\lambda)$ is spanned by the products of bases of the fundamental modules $V(\omega_i)$ (Theorem \ref{thmspan}). In particular, we also provide explicit polynomial bases for all six fundamental $\mathscr{G}^{E_6}$-modules. This gives a purely algebraic and constructive way to generate high-dimensional representations from the fundamental ones.

The polynomial realization also leads to modules of finite length over a commutative algebra. Let $\mathscr{G}^{E_6}_-$ be the $16$-dimensional abelian nilradical opposite to the maximal parabolic subalgebra whose Levi subalgebra has semisimple part of type $D_5$. For each dominant integral weight $\lambda$, restricting $V(\lambda)$ to $\mathscr{G}^{E_6}_{-}$ makes it a module over the polynomial algebra $\mathscr{S}=\mathbb{C}[z_1,\cdots,z_{16}]\cong U(\mathscr{G}^{E_6}_-)$. The degree in (\ref{1.3}) gives a grading structure on this $\mathscr{S}$-module. This construction relates the explicit representation to questions about the graded free resolution and lower bounds on Betti numbers $b_i(V(\lambda))$. The Buchsbaum-Eisenbud-Horrocks conjecture predicts that every nonzero finitely generated graded module $M$ of codimension $c$ over a polynomial ring $\mathscr{P}$ satisfies
\begin{eqnarray}
	b_i(M)\geqslant\binom{c}{i}.
\end{eqnarray}
This conjecture originates in the work of  Buchsbaum and Eisenbud \cite{BE} and in Horrocks' question on the rank of syzygies recorded in \cite{Ha}, see also \cite[Conjecture 4.1]{BG}. These inequalities remain open in general, although their summation, $\sum_i b_i(M)\geqslant 2^c$, is known by the works of Walker \cite{W} in characteristic $\neq2$, and of VandeBogert and Walker \cite{VW} in characteristic $=2$. A related question of Charalambous, Evans and Miller \cite{CEM} asks for a larger bound for the total Betti number when the module is not a complete intersection quotient (also see \cite[Question B]{BG}). Recent progress includes the verification of this bound for Gorenstein quotients of codimension $5$ or $6$ with pure resolution by Boocher, Huang and Wolf \cite{BHW}. Using Lie algebra homology, Sam, VandeBogert and Weyman \cite{SVW} established a result for the modules of finite length associated with the abelian parabolic nilradicals  of classical types. 

For the exceptional type $E_6$, we use the homology form of Kostant's theorem \cite{K} with the grading structure above to determine the graded Betti numbers of $V(\lambda)$. We prove that the minimal graded free resolution of $V(\lambda)$ is pure if and only if $\lambda=k(\omega_1+\cdots+\omega_6)$ for $k\in\mathbb{N}$. We also show that among the modules with $\lambda\neq0$, the minimum total Betti number is $435204$, reached exactly at $\lambda=\omega_1$ or $\omega_6$. In particular, the Buchsbaum-Eisenbud-Horrocks conjecture and the question of Charalambous-Evans-Miller are both verified for the $\mathscr{S}$-modules we considered. The construction provides an explicit $E_6$-type family for which purity and the minimum total Betti number can be completely determined.

The paper is organized as follows. In Section \ref{section2}, we recall some essential ingredients, emphasizing on the full conformal oscillator representation of $\mathfrak{o}(10)$ and Xu's functor from $\mathfrak{o}(10)$-modules to $\mathscr{G}^{E_6}$-modules. Section \ref{section3} contains the main general constructions. We first give explicit representation formulas for all root vectors of $\mathscr{G}^{E_6}$. Then we show that the restricted dual of the polynomial algebra $\mathscr{C}$ is isomorphic to the Verma module, and the cyclic submodule $U(\mathscr{G}^{E_6})(1)$ is an irreducible highest‑weight module. We also determine all $\mathscr{G}^{D_5}$-singular vectors in $\mathscr{C}$. In section \ref{section4}, we provide explicit bases for all the six fundamental modules and a finite spanning set for arbitrary finite-dimensional irreducible $\mathscr{G}^{E_6}$-modules $V(\lambda)$. Section \ref{section5} is devoted to study $V(\lambda)$ more deeply when it is finite dimensional, we determine the dimension of the $y$-homogeneous subspace, and compute the graded Betti numbers as an $U(\mathscr{G}^{E_6}_-)$-module.

\textbf{Conventions.} Throughout, $\mathbb{C},\mathbb{R},\mathbb{Z}$, $\mathbb{N}$ denote the complex number field, the real number field, the ring of integers, the set of nonnegative integers, respectively. For $i,j\in\mathbb{Z}$, we write $\overline{i,j}=\{k\in\mathbb{Z}\mid i\leqslant k\leqslant j\}$, this set is nonempty if and only if $i\leqslant j$. The bracket $\llbracket x\rrbracket$ for $x\in\mathbb{R}$ denotes the greatest integer not exceeding $x$. Empty sums $\sum$ and products $\prod$ are taken to be $0$ and $1$, respectively.
Let $\delta_{i,j}$ be the Kronecker $\delta$. For a Lie algebra $\mathfrak{g}$, $U(\mathfrak{g})$ denotes its universal enveloping algebra. 

\section{Preliminaries}\label{section2}
In this section, we recall two essential ingredients. Subsection \ref{section2.1} reviews the basic concepts of polynomials and differential operators. Subsection \ref{section2.2} presents the full conformal oscillator representation of $\mathfrak{o}(10)$ on the polynomial algebra $\mathscr{X}$ in $20$ variables, together with the explicit formulas for the action of the orthogonal Lie algebra. Subsection \ref{section2.3} reviews Xu’s functor from the category of $\mathfrak{o}(10)$-modules to the category of $\mathscr{G}^{E_6}$-modules, including the constructed module
$\widehat{M}=\mathscr{Y}\otimes M$ and the explicit differential operator representation of the $\mathscr{G}^{E_6}$ generators on $\widehat{M}$.
\subsection{Polynomial algebras and differential operators}\label{section2.1}
Let $n$ be a positive integer, and let $\mathscr{P}=\mathbb{C}[z_1,\cdots,z_n]$ be the polynomial algebra in $n$ variables. Denote $\partial_{z_i}\in\operatorname{End}_{\mathbb{C}}(\mathscr{P})$ to be the unique derivation satisfying
\begin{eqnarray}
	\partial_{z_i}(z_j)=\delta_{i,j}.
\end{eqnarray}
Identify each $f\in\mathscr{P}$ with the multiplication operator $g\mapsto fg$ on $\mathscr{P}$. The polynomial Witt algebra associated with $\mathscr{P}$ is
\begin{eqnarray}
	\mathscr{W}_{\mathscr{P}}=\operatorname{Der}_{\mathbb{C}}(\mathscr{P})=\bigoplus_{i=1}^n\mathscr{P}\partial_{z_i},
\end{eqnarray}
equipped with the commutator bracket. Explicitly,
\begin{eqnarray}
	[f_1\partial_{z_i},f_2\partial_{z_j}]=f_1\partial_{z_i}(f_2)\partial_{z_j}-f_2\partial_{z_j}(f_1)\partial_{z_i}
\end{eqnarray}
for $f_1,f_2\in\mathscr{P}$ and $i,j\in\overline{1,n}$. Thus $\mathscr{W}_{\mathscr{P}}$ is the Lie algebra of vector fields on the affine space $\mathbb{A}^n$.

The algebra of differential operators on $\mathscr{P}$ is the associative algebra
\begin{equation}
	\mathscr{D}_{\mathscr{P}}=\bigoplus_{a\in\mathbb{N}^n}\mathscr{P}\partial^{a}, \quad \partial^{a}=\partial_{z_1}^{a_1}\cdots\partial_{z_n}^{a_n},
\end{equation}
with the multiplication given by composition. Equivalently, $\mathscr{D}_{\mathscr{P}}$ is the $n$-th Weyl algebra generated by $z_1,\cdots,z_n,\partial_{z_1},\cdots,\partial_{z_n}$. For $m\in\mathbb{N}$, let
\begin{eqnarray*}
	\mathscr{D}_{\mathscr{P}}^{\leqslant m}=\bigoplus_{|a|\leqslant m}\mathscr{P}\partial^{a}\quad (|a|=a_1+\cdots+a_n)
\end{eqnarray*}
be the subspace of differential operators of order at most $m$. In particular
\begin{eqnarray}
	\mathscr{D}_{\mathscr{P}}^{\leqslant 1}=\mathscr{P}\oplus\mathscr{W}_{\mathscr{P}}
\end{eqnarray}
is a Lie subalgebra of $\mathscr{D}_{\mathscr{P}}$, with the Lie bracket given by
\begin{eqnarray}
	[f_1+d_1,f_2+d_2]=[d_1,d_2]+d_1(f_2)-d_2(f_1),\quad f_1,f_2\in\mathscr{P},\;d_1,d_2\in\mathscr{W}_{\mathscr{P}}.
\end{eqnarray}
A realization of a Lie algebra $\mathfrak{g}$ by inhomogeneous first-order differential operators is a Lie algebra homomorphism
\begin{eqnarray*}
	\mathfrak{g}\longrightarrow \mathscr{D}_{\mathscr{P}}^{\leqslant1},
\end{eqnarray*}
which extends uniquely to an associative algebra homomorphism
\begin{eqnarray}
	U(\mathfrak{g})\longrightarrow \mathscr{D}_{\mathscr{P}}.
\end{eqnarray}

\subsection{Differential operator realization for $D_5$}\label{section2.2}
We reorder the rows and columns of $10\times 10$ matrices in $M_{10\times 10}(\mathbb{C})$ by the indices
\begin{eqnarray}
	1,\:2,\:3,\:4,\:5,\:-1,\:-2,\:-3,\:-4,\:-5.
\end{eqnarray}
Denote by $E_{i,j}$ the $10\times 10$ matrix with $1$ as its $(i,j)$-entry and $0$ elsewhere. Set
\begin{eqnarray}
	A_{i,j}=E_{i,j}-E_{-j,-i}\;\;\;\text{for}\;i,j\in\overline{-5,-1}\cup\overline{1,5}.
\end{eqnarray}
Note that $A_{i,j}=-A_{-j,-i}$ and $A_{i,-i}=0_{10\times 10}$. The orthogonal Lie algebra is given by
\begin{eqnarray}
	&\mathfrak{o}(10)=\sum\limits_{i,j=1}^5\mathbb{C}A_{i,j}+\sum\limits_{1\leqslant i<j\leqslant 5}(\mathbb{C}A_{i,-j}+\mathbb{C}A_{-i,j}), \label{even}
\end{eqnarray}
which is isomorphic to the simple Lie algebra of type $D_5$.
We take $\mathfrak{h}_0=\oplus_{r=1}^5\mathbb{C}A_{r,r}$ as the Cartan subalgebra and define $\varepsilon_{r}(A_{s,s})=\delta_{r,s}$. The positive root system is
\begin{eqnarray}
	\Phi^{+}_{\mathfrak{o}(10)}=\{\varepsilon_i-\varepsilon_j,\varepsilon_i+\varepsilon_j\mid 1\leqslant j<i\leqslant 5\},
\end{eqnarray}
the corresponding positive root vectors are $A_{i,j}$ with $1\leqslant |j| <i\leqslant 5$. For arbitrary $\mu_1,\cdots,\mu_5\in\mathbb{C}$, write $\mu=\sum\limits_{r=1}^5\mu_r\varepsilon_r\in\mathfrak{h}_0^*$ .

Denote $\mathscr{X}=\mathbb{C}[x_{i,j}\mid 1\leqslant|j|\leqslant i\leqslant 4]$, and put
\begin{eqnarray}
	p_{i,j}=x_{i,-j}+\sum\limits_{r=1}^{j-1}x_{j-1,r}x_{i,-r},\qquad \partial_{x_{i,j}}=\frac{\partial}{\partial x_{i,j}}.
\end{eqnarray} 
Define $5\times 5$ matrices in $M_{5\times 5}(\mathscr{X})$ with entries in $\mathscr{X}$,
\begin{eqnarray}
	L=\begin{pmatrix}
		1&0&0&0&0\\
		x_{1,1}&1&0&0&0\\
		x_{2,1}&x_{2,2}&1&0&0\\
		x_{3,1}&x_{3,2}&x_{3,3}&1&0\\
		x_{4,1}&x_{4,2}&x_{4,3}&x_{4,4}&1
	\end{pmatrix}
\text{~and~}
	K=\begin{pmatrix}
		0&-p_{1,1}&-p_{2,1}&-p_{3,1}&-p_{4,1}\\
		p_{1,1}&0&-p_{2,2}&-p_{3,2}&-p_{4,2}\\
		p_{2,1}&p_{2,2}&0&-p_{3,3}&-p_{4,3}\\
		p_{3,1}&p_{3,2}&p_{3,3}&0&-p_{4,4}\\
		p_{4,1}&p_{4,2}&p_{4,3}&p_{4,4}&0
	\end{pmatrix}. \label{2.8}
\end{eqnarray}
Note that $L^{-1}=I_5+\sum_{r=1}^4(I_5-L)^r$ also has polynomial entries, we define
\begin{eqnarray}
	C=K(L^{-1})^t. \label{2.9}
\end{eqnarray}
In particular, we can direct compute
\begin{align}
	L_{i+1,j}=x_{i,j}\quad\text{and}\quad C_{i+1,j}=x_{i,-j}\quad\text{for }1\leqslant j\leqslant i\leqslant 4,\label{lc}
\end{align}
and from $CL^t+LC^t=0$ we solve
\begin{eqnarray}
	C_{i,i}=-\sum_{j=1}^{i-1}x_{i-1,j}x_{i-1,-j}\quad\text{for }i\in\overline{2,5}. \label{cc}
\end{eqnarray}
For $u\in\mathfrak{o}(10)$, we write
\begin{align}
	u=\begin{pmatrix}
		P&Q\\R&-P^{t}
	\end{pmatrix}\;\text{ with }P,R,Q\in M_{5\times 5}(\mathbb{C})\text{ and }Q^{t}=-Q,~R^{t}=-R.
\end{align}
Moreover. Let $U$ be the upper triangular part of $(LP+CR)L^{-1}$ including the diagonal. Then $(LP+CR)L^{-1}-U$ is a strictly lower triangular matrix, and so is $LP+CR-UL$.
Define a map $\varrho_{\mu}:\mathfrak{o}(10)\longrightarrow\operatorname{End}_{\mathbb{C}}(\mathscr{X})$ by
\begin{eqnarray}\label{2.11}
	\varrho_{\mu}(u)=\sum_{1\leqslant s\leqslant t\leqslant 4}\!\!&\!\!\!\!&\!\!((LP+CR-UL)_{t+1,s}\partial_{x_{t,s}}\nonumber\\
	\!\!&\!\!\!\!&\!\!+(LQ-CP^t-UC)_{t+1,s}\partial_{x_{t,-s}})+\sum_{r=1}^5\mu_rU_{r,r}.\label{rho}
\end{eqnarray}
The last summation of polynomials acts on $\mathscr{X}$ by multiplication. Denote by $d_u$ the part of differential operator of $\varrho_{\mu}(u)$. Equivalently,
\begin{eqnarray}
	d_u=\varrho_{\mu}(u)-\sum_{r=1}^5\mu_rU_{r,r}=\varrho_0(u).
\end{eqnarray}
Let $d_u$ acts entry-wise on the matrices. By above arguments we have
\begin{eqnarray}
	d_u(L)=LP+CR-UL. \label{L}
\end{eqnarray}
Let $E=d_u(C)-LQ+CP^t+UC$, then $E$ is upper triangular. Moreover, acting $d_u$ on $CL^t+LC^t=0$ we obtain
\begin{align}
	0=&(E+LQ-CP^t-UC)L^t+C(LP+CR-UL)^t\nonumber\\
	 &+(LP+CR-UL)C^t+L(E+LQ-CP^t-UC)^t\nonumber\\
	 =&EL^t+LE^t.
\end{align}
Thus $EL^t$ is skew-symmetric while it is also upper triangular, hence is zero. Since $L$ is invertible, $E=0$. We therefore have
\begin{eqnarray}
	d_u(C)=LQ-CP^t-UC. \label{C}
\end{eqnarray}
Hence (\ref{L}) and (\ref{C}) gives
\begin{eqnarray}
	d_u\begin{pmatrix}
		L&C
	\end{pmatrix}=\begin{pmatrix}
	L&C
	\end{pmatrix}u-U\begin{pmatrix}
	L&C
	\end{pmatrix}. \label{LC}
\end{eqnarray}
\begin{proposition}
	For any $\mu\in\mathfrak{h}_0^*$, the map $\varrho_\mu$ is a Lie algebra homomorphism.
\end{proposition}
\begin{proof}
	The linearity of $\varrho_{\mu}$ is easy to see since $L$ and $C$ are fixed and the projection onto the upper triangular part is linear.
	It leaves to show that $\varrho_{\mu}$ preserves the Lie bracket. For $u,v\in\mathfrak{o}(10)$, write $U_u$, $U_v$ respectively for the upper triangular matrix $U$ defined above. Eq. (\ref{LC}) gives
	\begin{eqnarray}
		[d_u,d_v]\begin{pmatrix}
			L&C
		\end{pmatrix}=\begin{pmatrix}
		L&C
		\end{pmatrix}[u,v]-(d_u(U_v)-d_v(U_u)+[U_u,U_v])\begin{pmatrix}
		L&C
		\end{pmatrix}, \label{dudv}
	\end{eqnarray}
	from which we get
	\begin{eqnarray}
		([d_u,d_v]-d_{[u,v]})\begin{pmatrix}
			L&C
		\end{pmatrix}=(U_{[u,v]}-d_u(U_v)+d_v(U_u)-[U_u,U_v])\begin{pmatrix}
		L&C
		\end{pmatrix}
	\end{eqnarray}
	Hence
	\begin{eqnarray}
		([d_u,d_v]-d_{[u,v]})(L)L^{-1}=U_{[u,v]}-d_u(U_v)+d_v(U_u)-[U_u,U_v]. \label{dd}
	\end{eqnarray}
	Since $[d_u,d_v]-d_{[u,v]}$ is a derivation, $([d_u,d_v]-d_{[u,v]})(L)$ is strictly lower triangular by (\ref{2.8}). Thus in (\ref{dd}) the left hand side is strictly lower triangular, while the right hand side is upper triangular. Both sides must be zero. Consequently,
	\begin{eqnarray}
		U_{[u,v]}=d_u(U_v)-d_v(U_u)+[U_u,U_v]. \label{uuv}
	\end{eqnarray}
	By (\ref{dudv}) and (\ref{uuv}) we have
	\begin{eqnarray}
		[d_u,d_v]\begin{pmatrix}
			L&C
		\end{pmatrix}=d_{[u,v]}\begin{pmatrix}
		L&C
		\end{pmatrix}.
	\end{eqnarray}
	Eq. (\ref{lc}) shows that the entries $L_{i+1,j}$ and $C_{i+1,j}$ for $1\leqslant j\leqslant i\leqslant 4$ generate $\mathscr{X}$, hence
	\begin{eqnarray}
	    [d_u,d_v]=d_{[u,v]}.
	\end{eqnarray}
	Moreover, the matrix $[U_u,U_v]$ is strictly upper triangular since each factor of the bracket is upper triangular, so (\ref{uuv}) also implies that the $(r,r)$-entry on both sides satisfy
	\begin{eqnarray}
		(U_{[u,v]})_{r,r}=(d_u(U_v))_{r,r}-(d_v(U_u))_{r,r}.
	\end{eqnarray}
	Now by (\ref{rho}) we have
	\begin{align}
		[\varrho_{\mu}(u),\varrho_{\mu}(v)]=&\:[d_u+\sum_{r=1}^{5}\mu_r(U_u)_{r,r},d_v+\sum_{r=1}^{5}\mu_r(U_v)_{r,r}]\nonumber\\
		=&\:[d_u,d_v]+\sum_{r=1}^5\mu_r\big((d_u(U_v))_{r,r}-(d_v(U_u))_{r,r}\big)\nonumber\\
		=&\:d_{[u,v]}+\sum_{r=1}^5\mu_r(U_{[u,v]})_{r,r}\nonumber\\
		=&\:\varrho_{\mu}([u,v]).
	\end{align}
	This proves the proposition.	
\end{proof}
The positive root vectors and Cartan elements take the particularly simple forms,
\begin{eqnarray}
	\varrho_{\mu}(A_{i,j})=\sum_{k=i}^4(x_{k,i}\partial_{x_{k,j}}-x_{k,-j}\partial_{x_{k,-i}})+\partial_{x_{i-1,j}}\quad (1\leqslant |j|<i\leqslant 5), \label{2.24}
\end{eqnarray}
and
\begin{eqnarray}
	\varrho_{\mu}(A_{i,i})=\sum\limits_{k=i}^4(x_{k,i}\partial_{x_{k,i}}-x_{k,-i}\partial_{x_{k,-i}})-\sum\limits_{|s|=1}^{i-1}x_{i-1,s}\partial_{x_{i-1,s}}+\mu_i\quad(i\in\overline{1,5}). \label{2.25}
\end{eqnarray}
Eq. (\ref{2.11}) specifies the action of each $A_{i,j}$ by polynomial differential operators of order at most one, thus the image $\varrho_{\mu}(\mathfrak{o}(10))$ is indeed a Lie subalgebra of $\mathscr{D}_{\mathscr{X}}^{\leqslant 1}=\mathscr{X}\oplus{\mathscr{W}}_{\mathscr{X}}$. 

Using the methods developed by Xu and Zhao in \cite{XuZhao}, we realized the so-called full conformal oscillator representations of $\mathfrak{o}(10)$ on $\mathscr{X}$ \cite[Theorem 2.3]{Zhou}, where the actions of root vectors and the Cartan subalgebra agree with that given by $\varrho_{\mu}$ for fixed $\mu$ here. The interested reader is referred to the reference for details. 

\subsection{The functor from $D_5$-$\textbf{Mod}$ to $E_6$-$\textbf{Mod}$}\label{section2.3}
Recall the Dynkin diagram of type $E_6$:
\begin{center}
	\begin{picture}(20,40)
		\put(-130,0){$E_6$:}\put(-84,0){\circle{6}}\put(-84,
		-14){1}\put(-81,0){\line(1,0){36}}\put(-42,0){\circle{6}}\put(-42,
		-14){3}\put(-39,0){\line(1,0){36}}\put(0,0){\circle{6}}\put(0,
		-14){4}\put(0,3){\line(0,1){32}}\put(0,38){\circle{6}}\put(7,36){2}\put(3,0){\line(1,0){36}}
		\put(42,0){\circle{6}}\put(42,-14){5}\put(45,0){\line(1,0){36}}\put(84,0){\circle{6}}\put(84,
		-14){6}
	\end{picture}	
\end{center}

~

Let $\{\alpha_i\mid i\in\overline{1,6}\}$ be the positive simple roots corresponding to the vertices, let $\Phi_{E_6}$ be the root system of $E_6$, and let $Q_{E_6}=\sum_{i=1}^6\mathbb{Z}\alpha_i$ be the root lattice, and let $Q_{E_6}^+=\sum_{i=1}^6\mathbb{N}\alpha_i$ be the positive root lattice of type $E_6$. For $\alpha=\sum_{i=1}^6n_i\alpha_i\in Q_{E_6}$, define the height of $\alpha$ to be
\begin{eqnarray}
	&\operatorname{ht}(\alpha)=\sum\limits_{i=1}^6n_i.
\end{eqnarray}
Denote by $(\cdot\;,\cdot)$ the symmetric $\mathbb{Z}$-bilinear form on $Q_{E_6}$ such that $(\alpha,\alpha)=2$ for all $\alpha\in\Phi_{E_6}$. Define a map $F:Q_{E_6}\times Q_{E_6}\longrightarrow\{-1,1\}$ by
\begin{eqnarray}
	&F(\sum\limits_{i=1}^6\lambda_i\alpha_i,\sum\limits_{j=1}^6l_j\alpha_j)=(-1)^{\sum\limits_{i=1}^6\lambda_il_i+\lambda_1l_3+\lambda_4l_2+\lambda_3l_4+\lambda_5l_4+\lambda_6l_5},\;\;\;\lambda_i,l_j\in\mathbb{Z}. \label{2.37}
\end{eqnarray}
Extending $(\cdot\;,\cdot)$ on $\mathscr{H}=\oplus_{i=1}^6\mathbb{C}\alpha_i$, the simple Lie algebra of type $E_6$ is 
\begin{eqnarray}
	&\mathscr{G}^{E_6}=\mathscr{H}\oplus\bigoplus\limits_{\alpha\in\Phi_{E_6}}\mathbb{C}E_{\alpha} \label{2.26}
\end{eqnarray}
with the Lie bracket given by
\begin{eqnarray}
	[\mathscr{H},\mathscr{H}]=0,\quad[h,E_{\alpha}]=-[E_{\alpha},h]=(h,\alpha)E_{\alpha},\quad[E_{\alpha},E_{-\alpha}]=-\alpha \label{Lie}
\end{eqnarray}
for $h\in \mathscr{H}$, and
\begin{eqnarray}
	[E_{\alpha},E_{\beta}]=\Big\{\begin{array}{ll}
		0,&\text{if }\alpha+\beta\notin\Phi_{E_{6}},\\
		F(\alpha,\beta)E_{\alpha+\beta},&\text{if }\alpha+\beta\in\Phi_{E_6},
	\end{array} \label{2.28}
\end{eqnarray}
for $\alpha+\beta\neq0$ (see e.g., \cite{Xu}). The subspace $\mathscr{H}$ is a Cartan subalgebra of $\mathscr{G}^{E_6}$.
The corresponding Cartan matrix is 
\begin{align}
	A^{E_6}=\begin{pmatrix}
		2&0&-1&0&0&0\\
		0&2&0&-1&0&0\\
		-1&0&2&-1&0&0\\
		0&-1&-1&2&-1&0\\
		0&0&0&-1&2&-1\\
		0&0&0&0&-1&2
	\end{pmatrix}. \label{CM}
\end{align}
Set
\begin{eqnarray}
	&Q^{D_5}=\sum\limits_{i=1}^5\mathbb{Z}\alpha_i,\;\;\;\Phi_{D_5}=\Phi_{E_6}\bigcap Q^{D_5}.
\end{eqnarray}
Then
\begin{eqnarray}
	&\mathscr{G}^{D_5}=\sum\limits_{i=1}^5\mathbb{C}\alpha_i\oplus\sum\limits_{\alpha\in\Phi_{D_5}}\mathbb{C}E_{\alpha}
\end{eqnarray}
is a Lie subalgebra of $\mathscr{G}^{E_6}$ isomorphic to the simple Lie algebra of type $D_5$. Moreover, $\mathscr{H}^{D_5}=\oplus_{i=1}^5\mathbb{C}\alpha_i$ is the Cartan subalgebra of $\mathscr{G}^{D_5}$.

Denote by $\Phi_{E_6}^+$ the set of positive roots of $E_6$ and by $\Phi_{D_5}^+$ the set of positive roots of $D_5$. The elements in $\Phi^+_{E_6}\setminus\Phi^+_{D_5}$ are
\begin{eqnarray}
	&\beta_1=\alpha_6,\;\beta_2=\alpha_5+\alpha_6,\;\beta_3=\alpha_4+\alpha_5+\alpha_6,\;\beta_4=\alpha_3+\alpha_4+\alpha_5+\alpha_6,\label{2.29}\\
	&\beta_5=\alpha_2+\alpha_4+\alpha_5+\alpha_6,\;\beta_6=\alpha_2+\alpha_3+\alpha_4+\alpha_5+\alpha_6,\\
	&\beta_7=\alpha_1+\alpha_3+\alpha_4+\alpha_5+\alpha_6,\;\beta_8=\alpha_1+\alpha_2+\alpha_3+\alpha_4+\alpha_5+\alpha_6,\\
	&\beta_9=\alpha_2+\alpha_3+2\alpha_4+\alpha_5+\alpha_6,\;\beta_{10}=\alpha_2+\alpha_3+2\alpha_4+2\alpha_5+\alpha_6,\\
	&\beta_{11}=\alpha_1+\alpha_2+\alpha_3+2\alpha_4+\alpha_5+\alpha_6,\;\beta_{12}=\alpha_1+\alpha_2+2\alpha_3+2\alpha_4+\alpha_5+\alpha_6,\\
	&\beta_{13}=\alpha_1+\alpha_2+\alpha_3+2\alpha_4+2\alpha_5+\alpha_6,\;\beta_{14}=\alpha_1+\alpha_2+2\alpha_3+2\alpha_4+2\alpha_5+\alpha_6,\\
	&\beta_{15}=\alpha_1+\alpha_2+2\alpha_3+3\alpha_4+2\alpha_5+\alpha_6,\\
	&\beta_{16}=\alpha_1+2\alpha_2+2\alpha_3+3\alpha_4+2\alpha_5+\alpha_6,
\end{eqnarray}
and the elements in $\Phi^+_{D_5}$ are
\begin{eqnarray}
	&\alpha_1,\;\;\alpha_2,\;\;\alpha_3,\;\;\alpha_4,\;\;\alpha_5,\;\;\alpha_1+\alpha_3,\;\;\alpha_2+\alpha_4,\;\;\alpha_3+\alpha_4,\;\;\alpha_4+\alpha_5,\label{2.38}\\
	&\alpha_1+\alpha_3+\alpha_4,\;\;\alpha_2+\alpha_3+\alpha_4,\;\;\alpha_2+\alpha_4+\alpha_5,\;\;\alpha_3+\alpha_4+\alpha_5,\\
	&\alpha_1+\alpha_2+\alpha_3+\alpha_4,\;\;\alpha_1+\alpha_3+\alpha_4+\alpha_5,\;\;\alpha_2+\alpha_3+\alpha_4+\alpha_5,\\
	&\alpha_1+\alpha_2+\alpha_3+\alpha_4+\alpha_5,\;\;\alpha_2+\alpha_3+2\alpha_4+\alpha_5,\\
	&\alpha_1+\alpha_2+\alpha_3+2\alpha_4+\alpha_5,\;\;\alpha_1+\alpha_2+2\alpha_3+2\alpha_4+\alpha_5.\label{2.42}
\end{eqnarray}
In particular,
\begin{eqnarray}
	&\sum\limits_{\alpha\in\Phi_{E_6}^+}\mathbb{C}E_{\alpha}\;\;\;\text{and}\;\;\;\sum\limits_{\alpha\in\Phi_{E_6}^+}\mathbb{C}E_{-\alpha}
\end{eqnarray}
are the nilpotent subalgebra of positive root vectors and the nilpotent subalgebra of negative root vectors, respectively. The fundamental weights $\omega_i\in \mathscr{H}^*$ satisfy
\begin{eqnarray}
	\omega_i(\alpha_j)=\delta_{i,j}.
\end{eqnarray}
The isomorphism classes of finite-dimensional irreducible $\mathscr{G}^{E_6}$-modules are in one-to-one correspondence with the set of weight $\lambda=\sum_{i=1}^6\lambda_i\omega_i$ such that $\lambda_1,\lambda_2,\cdots,\lambda_6\in\mathbb{N}$. 
Write
\begin{eqnarray}
	\xi_i=E_{\beta_i} \;\;\;\text{and}\;\;\;\eta_i=E_{-\beta_i}\quad (i\in\overline{1,16}).  \label{2.43}
\end{eqnarray}
Set
\begin{eqnarray}
	&\mathscr{G}^{E_{6}}_+=\sum\limits_{i=1}^{16}\mathbb{C}\xi_i,\;\;\mathscr{G}^{E_{6}}_-=\sum\limits_{i=1}^{16}\mathbb{C}\eta_i,\;\;\mathscr{G}^{E_6}_0=\mathscr{G}^{D_5}+\mathbb{C}\alpha_6. \label{2.44}
\end{eqnarray}
Then $\mathscr{G}^{E_6}_{\pm}$ are abelian Lie subalgebra of $\mathscr{G}^{E_{6}}$ and $\mathscr{G}^{E_6}_0$ is a reductive Lie subalgebra of $\mathscr{G}^{E_{6}}$. Moreover,
\begin{eqnarray}
	\mathscr{G}^{E_6}=\mathscr{G}^{E_6}_-+\mathscr{G}^{E_6}_0+\mathscr{G}^{E_6}_+,\;\;[\mathscr{G}^{E_6}_0,\mathscr{G}^{E_6}_{\pm}]=\mathscr{G}^{E_6}_{\pm},\;\;[\mathscr{G}^{E_6}_+,\mathscr{G}^{E_6}_-]\subset\mathscr{G}^{E_6}_0. \label{2.41}
\end{eqnarray}
Furthermore, there exists a Lie algebra isomorphism $\nu:\mathfrak{o}(10)\longrightarrow \mathscr{G}^{D_5}$ determined by the generators:
\begin{eqnarray}
	&\nu(A_{1,1})=(\alpha_2-\alpha_5)/2,\;\nu(A_{2,2})=(\alpha_2+\alpha_5)/2,\;\nu(A_{3,3})=\alpha_4+(\alpha_2+\alpha_5)/2,\label{2.46}\\
	&\nu(A_{4,4})=\alpha_3+\alpha_4+(\alpha_2+\alpha_5)/2,\;\nu(A_{5,5})=\alpha_1+\alpha_3+\alpha_4+(\alpha_2+\alpha_5)/2,\label{2.47}\\
	&\nu(A_{5,4})=E_{\alpha_1},\;\;\nu(A_{4,5})=-E_{-\alpha_1},\;\;\nu(A_{4,3})=E_{\alpha_3},\;\;\nu(A_{3,4})=-E_{-\alpha_3},\label{2.48}\\
	&\nu(A_{3,2})=E_{\alpha_4},\;\;\nu(A_{2,3})=-E_{-\alpha_4},\;\;\nu(A_{2,1})=E_{\alpha_5},\label{2.49}\\
	&\nu(A_{1,2})=-E_{-\alpha_5},\;\;\nu(A_{2,-1})=E_{\alpha_2},\;\;\nu(A_{-1,2})=-E_{-\alpha_{2}}.\label{2.50}
\end{eqnarray}
From now on, we identify $\mathfrak{o}(10)$ with $\mathscr{G}^{D_5}$ via $\nu$.

Denote $\mathscr{Y}=\mathbb{C}[y_1,\cdots,y_{16}]$ and let $\mathscr{W}_{\mathscr{Y}}=\oplus_{k=1}^{16}\mathscr{Y}\partial_{y_k}$ be the Witt algebra associated to $\mathscr{Y}$. By (\ref{2.41}), we may write
\begin{eqnarray}
	&[u,\eta_i]=\sum\limits_{j=1}^{16}c_{j,i}(u)\eta_j,\quad (u\in\mathscr{G}^{E_6}_0,\;i\in\overline{1,16}),  \label{247}
\end{eqnarray}
with the coefficients $c_{j,i}(u)\in \mathbb{C}$. Then we obtain a homomorphism $\sigma:\mathscr{G}^{E_6}_0\longrightarrow\mathscr{W}_{\mathscr{Y}}$ of Lie algebras given by
\begin{eqnarray}
	&\sigma(u)=\sum\limits_{i,j=1}^{16}c_{j,i}(u)y_j\partial_{y_i}. \label{248}
\end{eqnarray}
Appendix \ref{app} exhibits these explicit coefficients $c_{j,i}(u)$ for the basis elements of $\mathscr{G}^{E_6}_0$.

Denote
\begin{eqnarray}
	&\zeta_1=y_1y_{10}+y_2y_9-y_3y_6+y_4y_5,\label{zeta1}\\
	&\zeta_2=y_1y_{13}+y_2y_{11}-y_3y_8+y_5y_7,\label{zeta2}\\
	&\zeta_3=y_1y_{14}+y_2y_{12}-y_4y_8+y_6y_7,\label{zeta3}\\
	&\zeta_4=y_1y_{15}-y_3y_{12}+y_4y_{11}-y_7y_9,\label{zeta4}\\
	&\zeta_5=-y_2y_{15}-y_3y_{14}+y_4y_{13}-y_7y_{10},\label{zeta5}\\
	&\zeta_6=y_7y_{16}+y_8y_{15}+y_{11}y_{14}-y_{12}y_{13},\label{zeta6}\\
	&\zeta_7=-y_{4}y_{16}-y_{6}y_{15}-y_9y_{14}+y_{10}y_{12},\label{zeta7}\\
	&\zeta_8=y_3y_{16}+y_5y_{15}+y_9y_{13}-y_{10}y_{11},\label{zeta8}\\
	&\zeta_9=y_2y_{16}-y_5y_{14}+y_6y_{13}-y_8y_{10},\label{zeta9}\\
	&\zeta_{10}=y_1y_{16}+y_5y_{12}-y_6y_{11}+y_8y_9,\label{zeta10}
\end{eqnarray}
and let
\begin{eqnarray}
	&\mathcal{D}_{y}=\sum\limits_{i=1}^{16}y_i\partial_{y_i}.
\end{eqnarray}
We extend the map $\sigma$ from $\mathscr{G}^{E_6}_0$ to $\mathscr{G}^{E_6}$ by letting
\begin{eqnarray}
	&\sigma(\xi_i)=\partial_{y_i},\quad \sigma(\eta_i)=y_i\mathcal{D}_{y}-\sum\limits_{\ell=1}^{10}\sum\limits_{k=1}^{16}\zeta_{\ell}\partial_{y_k}(\partial_{y_i}(\zeta_{\ell}))\partial_{y_k},\quad(i\in\overline{1,16}). \label{2.65}
\end{eqnarray}
Direct computations with (\ref{zeta1})-(\ref{zeta10}) show that each $\partial_{y_k}(\partial_{y_i}(\zeta_{\ell}))$ is a constant in $\{-1,0,1\}$, the explicit valuations are also listed in Appendix \ref{app}.
Moreover, we have
\begin{proposition}[\text{\cite[Theorem 4.1]{Xu}}]
	The map $\sigma:\mathscr{G}^{E_6}\longrightarrow\mathscr{W}_{\mathscr{Y}}$ is a monomorphism of Lie algebras. Therefore, the polynomial algebra $\mathscr{Y}$ becomes a $\mathscr{G}^{E_6}$-module via $\sigma$.
\end{proposition}

Write
\begin{eqnarray}
	\widehat{\alpha}=2\alpha_1+3\alpha_2+4\alpha_3+6\alpha_4+5\alpha_5+4\alpha_6.\label{2.104}
\end{eqnarray}
Then
\begin{eqnarray}
	[\widehat{\alpha},\mathscr{G}^{D_5}]=0 \text{ and } \sigma(\widehat{\alpha})=-3\mathcal{D}_{y}.
\end{eqnarray}
Let $M$ be an $\mathscr{G}^{D_5}$-module and fix $\kappa\in\mathbb{C}$. We extend $M$ to be a module over the Lie subalgebra $\mathscr{G}^{E_6}_0=\mathscr{G}^{D_5}\oplus\mathbb{C}\widehat{\alpha}$ by letting
\begin{eqnarray}
	\widehat{\alpha}\cdot m=2\kappa m \label{2.68}
\end{eqnarray}
for any $m\in M$.

Using Shen's mixed product in \cite{Shen}, Xu obtained the following $\mathscr{G}^{E_6}$-module structure on $\widehat{M}=\mathscr{Y}\otimes_{\mathbb{C}}M$:
\begin{eqnarray}
	&\widehat{\alpha}|_{\widehat{M}}=-3\mathcal{D}_{y}\otimes \text{Id}_M+2\kappa\text{Id}_{\mathscr{Y}}\otimes\text{Id}_M,\label{2.107}\\
	&A_{i,j}|_{\widehat{M}}=\sigma(A_{i,j})\otimes\text{Id}_M+\text{Id}_{\mathscr{Y}}\otimes A_{i,j}|_M,\label{2.109} \\
	&\xi_k|_{\widehat{M}}=\partial_{y_k}\otimes \text{Id}_M,\label{2.108}\\
	&\eta_k|_{\widehat{M}}= \sigma(\eta_k)\otimes\operatorname{Id}_{M}+\sum\limits_{r=1}^{16}y_r\otimes [\xi_r,\eta_k]|_M \label{2.67}
\end{eqnarray}
for $k\in\overline{1,16}$ and $i,j\in\overline{-5,-1}\cup\overline{1,5}$.
By (\ref{2.107})-(\ref{2.108}), we have
\begin{eqnarray}
	\mathscr{G}^{E_6}_+(1\otimes M)=\{0\},\;\;\mathscr{G}^{E_6}_0(1\otimes M)\subseteq1\otimes M.
\end{eqnarray}
By the PBW Theorem, we have
\begin{eqnarray}
	U(\mathscr{G}^{E_6})(1\otimes M)=U(\mathscr{G}^{E_6}_-)(1\otimes M).
\end{eqnarray}
\begin{proposition}[\text{\cite[Theorem 5.1 and Proposition 5.2]{Xu}}]  The map $M\longrightarrow\widehat{M}$ gives a functor from the category of  $\mathscr{G}^{D_5}$-modules to the category of $\mathscr{G}^{E_6}$-modules. Moreover, if $M$ is an irreducible $\mathscr{G}^{D_5}$-module, then $U(\mathscr{G}^{E_6}_-)(1\otimes M)$ is an irreducible $\mathscr{G}^{E_6}$-module.
\end{proposition}

\section{General constructions}\label{section3}
In this section, we realize $\mathscr{G}^{E_6}$ as inhomogeneous first-order differential operators in $36$ variables. We first describe the action of $\mathscr{G}^{E_6}$ on $\mathscr{C}$, then show that for any weight $\lambda\in \mathscr{H}^*$, the $\mathscr{G}^{E_6}$-module $\mathscr{C}$ is isomorphic to the restricted dual Verma module, and the cyclic module $U(\mathscr{G}^{E_6})(1)$ is an irreducible highest-weight module. 

\subsection{Differential operator realization for $E_6$}\label{section3.1}
Take $M=\mathscr{X}$ to be the $\mathfrak{o}(10)$-module constructed in Section \ref{section2.2}. Applying the functor from $D_5$-\textbf{Mod} to $E_6$-\textbf{Mod} in Section \ref{section2.3}, Eqs. (\ref{2.107})-(\ref{2.67}) yield a representation of $\mathscr{G}^{E_{6}}$  on
\begin{eqnarray}
	\mathscr{Y}\otimes_{\mathbb{C}}\mathscr{X}\cong \mathbb{C}[y_k,x_{i,j}\mid 1\leqslant k\leqslant 16, 1\leqslant |j|\leqslant i\leqslant 4]=:\mathscr{C}.
\end{eqnarray}
Equivalently, we obtain a Lie algebra homomorphism
\begin{eqnarray}
	\pi_{\mu,\kappa}:\mathscr{G}^{E_6}\longrightarrow \operatorname{End}_{\mathbb{C}}(\mathscr{C}).
\end{eqnarray}
We may treat $\mathscr{X}$ and $\mathscr{Y}$ as subalgebras of $\mathscr{C}$, and treat $\mathscr{W}_{\mathscr{X}}$ and $\mathscr{W}_{\mathscr{Y}}$ as Lie subalgebras of $\mathscr{W}_{\mathscr{C}}\subset\mathscr{D}_{\mathscr{C}}$.

Recall the symmetric invariant form $(\cdot,\cdot)$ on $\mathscr{G}^{E_6}$ such that $(\alpha,\alpha)=2$ for all $\alpha\in\Phi_{E_6}$. Moreover,
\begin{eqnarray}
	(\xi_i,\eta_j)=-\delta_{i,j},\quad (\widehat{\alpha},\widehat{\alpha})=12,\quad (\widehat{\alpha},\mathscr{G}^{D_5})=0.\label{3.2}
\end{eqnarray}
By (\ref{247})-(\ref{248}) and (\ref{3.2}) we have
\begin{eqnarray}
	&([\xi_j,\eta_i],\widehat{\alpha})=(\xi_j,[\eta_i,\widehat{\alpha}])=(\xi_j,3\eta_i)=-3\delta_{i,j}, \label{3.3}
\end{eqnarray}
and
\begin{eqnarray}
	&([\xi_j,\eta_i],u)=(\xi_j,[\eta_i,u])=(\xi_j,-\sum\limits_{\ell=1}^{16}c_{\ell,i}(u)\eta_\ell)=c_{j,i}(u). \label{3.4}
\end{eqnarray}
Note that the matrices $A_{r,r}$ with $r\in\overline{1,5}$ and $A_{p,q},A_{q,p}$ with $1\leqslant |q|<p\leqslant 5$ form a basis of $\mathfrak{o}(10)$, and $A_{r,r}$ with $r\in\overline{1,5}$ and $A_{q,p},A_{p,q}$ with $1\leqslant |q|<p\leqslant 5$ form a dual basis with respect to $(\cdot,\cdot)$. Since $[\xi_j,\eta_i]\in\mathscr{G}^{E_6}_0$, we deduce from (\ref{3.2})-(\ref{3.4}) that
\begin{eqnarray}
	\sum\limits_{j=1}^{16}y_j[\xi_j,\eta_i]&=&-\frac{1}{4}y_i\widehat{\alpha}+\sum\limits_{r=1}^5\sigma(A_{r,r})(y_i)A_{r,r}\nonumber\\
	&&+\sum\limits_{1\leqslant |q|<p\leqslant 5} (\sigma(A_{p,q})(y_i)A_{q,p}+\sigma(A_{q,p})(y_i)A_{p,q}) \label{3.6}
\end{eqnarray}
for $i\in\overline{1,16}$. Let 
\begin{eqnarray}
	&Z_i=\sum\limits_{r=1}^5\sigma(A_{r,r})(y_i)A_{r,r}+\sum\limits_{1\leqslant |q|<p\leqslant 5} (\sigma(A_{p,q})(y_i)A_{q,p}+\sigma(A_{q,p})(y_i)A_{p,q}), \label{zi}
\end{eqnarray}
then $Z_i$ is an element of $\mathscr{Y}\otimes_{\mathbb{C}}\mathfrak{o}(10)$, so we may write
\begin{eqnarray}
	Z_i=\begin{pmatrix}
		P_i&Q_i\\R_i&-P_i^{t} \label{zi1}
	\end{pmatrix},
\end{eqnarray}
where $P_i,Q_i,R_i\in M_{5\times 5}(\mathscr{Y})$ and $Q_i^t=-Q_i$, $R_i^t=-R_i$. In particular, we have
\begin{eqnarray}
	(P_i)_{p,q}=\sigma(A_{q,p})(y_i),\;\;(Q_i)_{p,q}=\sigma(A_{-q,p})(y_i),\;\;(R_i)_{p,q}=\sigma(A_{q,-p})(y_i)
\end{eqnarray}
for $1\leqslant p,q\leqslant 5$. Recall the matrices $L$, $K$ and $C$ defined in (\ref{2.8})-(\ref{2.9}). Let $U_i$ be the upper triangular part of $(LP_i+CR_i)L^{-1}$ including the diagonal. 
\begin{theorem}\label{thm3.1}
	For $\mu=(\mu_1,\cdots,\mu_5)\in\mathbb{C}^5$ and $\kappa\in\mathbb{C}$, the homomorphism $\pi_{\mu,\kappa}$ is determined by 
	\begin{eqnarray}
		&\pi_{\mu,\kappa}(\widehat{\alpha})=-3\mathcal{D}_y+2\kappa,\label{r1}\\
		&\pi_{\mu,\kappa}(u)=\sigma(u)+\varrho_{\mu}(u),\label{r2}\\
		&\pi_{\mu,\kappa}(\xi_i)=\partial_{y_i},\label{r3}
	\end{eqnarray}
	and
	\begin{eqnarray}
		\!\!\!\!\pi_{\mu,\kappa}(\eta_i)\!\!&\!\!=\!\!&\!\!
		\sum\limits_{1\leqslant s\leqslant t\leqslant 4}\big((LP_i+CR_i-U_iL)_{t+1,s}\partial_{x_{t,s}}+(LQ_i-CP_i^t-U_iC)_{t+1,s}\partial_{x_{t,-s}}\big)  \nonumber\\
		\!\!&\!\!\!\!&\!\!+y_i\mathcal{D}_y-\sum_{\ell=1}^{10}\sum\limits_{k=1}^{16}\zeta_{\ell}\partial_{y_k}(\partial_{y_i}(\zeta_{\ell}))\partial_{y_k}\!-\!\frac{\kappa}{2}y_i+\!\sum_{r=1}^5\!\mu_r(U_i)_{r,r} 	 \label{r4}	
	\end{eqnarray}
	for $u\in\mathscr{G}^{D_5}$, $i\in\overline{1,16}$. Consequently, its image $\pi_{\mu,\kappa}(\mathscr{G}^{E_6})$ is a Lie subalgebra of $\mathscr{D}_{\mathscr{C}}^{\leqslant 1}$.
\end{theorem}
\begin{proof}
	Eqs. (\ref{r1})-(\ref{r3}) are derived from (\ref{2.107})-(\ref{2.108}) directly. By (\ref{2.67}) we have
	\begin{eqnarray}
		\pi_{\mu,\kappa}(\eta_i)=\sigma(\eta_i)+\sum\limits_{j=1}^{16}y_j\varrho_{\mu}([\xi_j,\eta_i]). \label{3.14}
	\end{eqnarray}
	By (\ref{2.68}) and (\ref{3.6}) we have
	\begin{eqnarray}
		\!\!&\!\!\!\!&\!\!\sum\limits_{j=1}^{16}y_j\varrho_{\mu}([\xi_j,\eta_i])+\frac{\kappa}{2}y_i\nonumber\\
		\!\!&\!\!=\!\!&\!\!\sum\limits_{r=1}^5\!\sigma(A_{r,r})(y_i)\varrho_{\mu}(A_{r,r})+\!\sum\limits_{1\leqslant |q|<p\leqslant 5}\! (\sigma(A_{p,q})(y_i)\varrho_{\mu}(A_{q,p})+\sigma(A_{q,p})(y_i)\varrho_{\mu}(A_{p,q})).\label{3.15}
	\end{eqnarray}
	We may extend the map $\varrho_{\mu}$ to $\mathscr{Y}\otimes_{\mathbb{C}}\mathfrak{o}(10)$ by $\varrho_{\mu}(fu)=f\varrho_{\mu}(u)$ with $f\in\mathscr{Y}$ and $u\in\mathfrak{o}(10)$. Then (\ref{zi}) and (\ref{3.15}) give
	\begin{eqnarray}
		\sum\limits_{j=1}^{16}y_j\varrho_{\mu}([\xi_j,\eta_i])+\frac{\kappa}{2}y_i=\varrho_{\mu}(Z_i).
	\end{eqnarray}
	The linearity of $\varrho_{\mu}$ and (\ref{rho}) imply
	\begin{eqnarray}
		\varrho_{\mu}(Z_i)&=&\sum\limits_{1\leqslant s\leqslant t\leqslant 4}\big((LP_i+CR_i-U_iL)_{t+1,s}\partial_{x_{t,s}}\nonumber\\
		&&+(LQ_i-CP_i^t-U_iC)_{t+1,s}\partial_{x_{t,-s}}\big)+\!\sum_{r=1}^5\!\mu_r(U_i)_{r,r}. \label{3.17}
	\end{eqnarray}
	Now (\ref{2.65}), (\ref{3.14}) and (\ref{3.17}) together give (\ref{r4}).
	
	Since $\mathscr{G}^{E_6}$ is spanned by $\mathscr{G}^{D_5}$, $\widehat{\alpha}$, $\xi_i$ and $\eta_i$ ($i\in\overline{1,16}$), we immediately have
	\begin{eqnarray}
		\pi_{\mu,\kappa}(\mathscr{G}^{E_6})\subset \mathscr{D}_{\mathscr{C}}^{\leqslant 1}.
	\end{eqnarray}
	The proof is completed.
\end{proof}

Next we write down more explicit representation formulas for the Chevalley generators of $\mathscr{G}^{E_6}$. Denote by $\omega_1,\cdots,\omega_6\in \mathscr{H}^*$ the fundamental weights of $\mathscr{G}^{E_6}$ characterized by $\omega_i(\alpha_j)=\delta_{i,j}$. It is convenient to express the parameters of the construction in coefficients of fundamental weights.
For arbitrary $\lambda_1,\cdots,\lambda_6\in\mathbb{C}$, put
\begin{eqnarray}
	\lambda=\sum_{i=1}^6\lambda_i\omega_i.
\end{eqnarray}
Take the parameters $\mu$ and $\kappa$ in Theorem \ref{thm3.1} as
\begin{eqnarray}
	&\mu_1=\:\frac{1}{2}(\lambda_2-\lambda_5),\quad\mu_2=\frac{1}{2}(\lambda_2+\lambda_5),\quad \mu_3=\:\frac{1}{2}(\lambda_2+\lambda_5)+\lambda_4, \label{3.122}\\
	&\mu_4=\frac{1}{2}(\lambda_2+\lambda_5)+\lambda_4+\lambda_3,\quad
	\mu_5=\:\frac{1}{2}(\lambda_2+\lambda_5)+\lambda_4+\lambda_3+\lambda_1,\\
	&\kappa=\:\lambda_1+\frac{3}{2}\lambda_2+2\lambda_3+3\lambda_4+\frac{5}{2}\lambda_5+2\lambda_6,\label{3.124}
\end{eqnarray}
and write $\pi_{\lambda}$ in place of $\pi_{\mu,\kappa}$.
Inversely, we have
\begin{eqnarray}
	&\lambda_1=\mu_5-\mu_4,\;\;\lambda_2=\mu_2+\mu_1,\;\;\lambda_3=\mu_4-\mu_3,\;\;\lambda_4=\mu_3-\mu_2,\\
	&\lambda_5=\mu_2-\mu_1,\;\;\;\lambda_6=\frac{1}{2}(\mu_1-\mu_2-\mu_3-\mu_4-\mu_5+\kappa). \label{lam6}
\end{eqnarray}
\begin{theorem}\label{thm3.2}
	For parameters $\lambda_1,\cdots,\lambda_6\in\mathbb{C}$, the following formulas hold.
	\begin{enumerate}[(i)]
		\item The simple positive root vectors.
		\begin{align}
			\pi_{\lambda}(E_{\alpha_1})=&\:y_4\partial_{y_7}+y_6\partial_{y_8}+y_9\partial_{y_{11}}+y_{10}\partial_{y_{13}}+\partial_{x_{4,4}},\label{3.91}\\
			\pi_{\lambda}(E_{\alpha_2})=&\:-y_3\partial_{y_5}-y_4\partial_{y_6}-y_7\partial_{y_8}+y_{15}\partial_{y_{16}}+x_{4,2}\partial_{x_{4,-1}}-x_{4,1}\partial_{x_{4,-2}}\nonumber\\
			&\:+x_{3,2}\partial_{x_{3,-1}}-x_{3,1}\partial_{x_{3,-2}}+x_{2,2}\partial_{x_{2,-1}}-x_{2,1}\partial_{x_{2,-2}}+\partial_{x_{1,-1}},\label{ea2}\\
			\pi_{\lambda}(E_{\alpha_3})=&\:y_3\partial_{y_4}+y_5\partial_{y_6}+y_{11}\partial_{y_{12}}+y_{13}\partial_{y_{14}}+x_{4,4}\partial_{x_{4,3}}-x_{4,-3}\partial_{x_{4,-4}}+\partial_{x_{3,3}},\label{ea3}\\
			\pi_{\lambda}(E_{\alpha_4})=&\:-y_2\partial_{y_3}-y_6\partial_{y_9}-y_8\partial_{y_{11}}+y_{14}\partial_{y_{15}}+x_{4,3}\partial_{x_{4,2}}-x_{4,-2}\partial_{x_{4,-3}}\nonumber\\
			&\:+x_{3,3}\partial_{x_{3,2}}-x_{3,-2}\partial_{x_{3,-3}}+\partial_{x_{2,2}},\label{ea4}\\
			\pi_{\lambda}(E_{\alpha_5})=&\:-y_1\partial_{y_2}+y_9\partial_{y_{10}}+y_{11}\partial_{y_{13}}+y_{12}\partial_{y_{14}}+x_{4,2}\partial_{x_{4,1}}-x_{4,-1}\partial_{x_{4,-2}}\nonumber\\
			&\:+x_{3,2}\partial_{x_{3,1}}-x_{3,-1}\partial_{x_{3,-2}}+x_{2,2}\partial_{x_{2,1}}-x_{2,-1}\partial_{x_{2,-2}}+\partial_{x_{1,1}},\label{3.95}\\
			\pi_{\lambda}(E_{\alpha_6})=&\:\partial_{y_1}.\label{3.96}
		\end{align}
		\item The simple coroots.
		\begin{align}
			\pi_{\lambda}(\alpha_1)=&\:y_4\partial_{y_4}+y_6\partial_{y_6}-y_7\partial_{y_7}-y_8\partial_{y_8}+y_9\partial_{y_9}+y_{10}\partial_{y_{10}}-y_{11}\partial_{y_{11}}-y_{13}\partial_{y_{13}}\nonumber \\
			&\:-\mathop{\textstyle\sum}\limits_{|r|=1}^3(x_{4,r}\partial_{x_{4,r}}-x_{3,r}\partial_{x_{3,r}})-2x_{4,4}\partial_{x_{4,4}}+\lambda_1,\label{3.97}\\
			\pi_{\lambda}(\alpha_2)=&\:y_3\partial_{y_3}+y_4\partial_{y_4}-y_5\partial_{y_5}-y_6\partial_{y_6}+y_7\partial_{y_7}-y_8\partial_{y_8}+y_{15}\partial_{y_{15}}-y_{16}\partial_{y_{16}}\nonumber\\
			&\:+\mathop{\textstyle\sum}\limits_{i=2}^4(x_{i,2}\partial_{x_{i,2}}-x_{i,-2}\partial_{x_{i,-2}}+x_{i,1}\partial_{x_{i,1}}-x_{i,-1}\partial_{x_{i,-1}})\nonumber\\
			&\:-2x_{1,-1}\partial_{x_{1,-1}}+\lambda_2,
		\end{align}
		\begin{align}
			\pi_{\lambda}(\alpha_3)=&\:y_3\partial_{y_3}-y_4\partial_{y_4}+y_{5}\partial_{y_5}-y_6\partial_{y_6}+y_{11}\partial_{y_{11}}-y_{12}\partial_{y_{12}}+y_{13}\partial_{y_{13}}-y_{14}\partial_{y_{14}}\nonumber\\
			&\:+x_{4,4}\partial_{x_{4,4}}-x_{4,-4}\partial_{x_{4,-4}}-x_{4,3}\partial_{x_{4,3}}+x_{4,-3}\partial_{x_{4,-3}}\nonumber\\
			&\:-\mathop{\textstyle\sum}\limits_{|r|=1}^2(x_{3,r}\partial_{x_{3,r}}-x_{2,r}\partial_{x_{2,r}})-2x_{3,3}\partial_{x_{3,3}}+\lambda_3,  \\
			\pi_{\lambda}(\alpha_4)=&\:y_2\partial_{y_2}-y_{3}\partial_{y_3}+y_6\partial_{y_6}+y_{8}\partial_{y_8}-y_9\partial_{y_9}-y_{11}\partial_{y_{11}}+y_{14}\partial_{y_{14}}-y_{15}\partial_{y_{15}}\nonumber\\
			&\:+\mathop{\textstyle\sum}\limits_{i=3}^4(x_{i,3}\partial_{x_{i,3}}-x_{i,-3}\partial_{x_{i,-3}}-x_{i,2}\partial_{x_{i,2}}+x_{i,-2}\partial_{x_{i,-2}})\nonumber\\
			&\:-\mathop{\textstyle\sum}\limits_{r\in\{1,-1\}}(x_{2,r}\partial_{x_{2,r}}-x_{1,r}\partial_{x_{1,r}})-2x_{2,2}\partial_{x_{2,2}}+\lambda_4,  \\
			\pi_{\lambda}(\alpha_5)=&\:y_1\partial_{y_{1}}-y_{2}\partial_{y_2}+y_9\partial_{y_9}-y_{10}\partial_{y_{10}}+y_{11}\partial_{y_{11}}+y_{12}\partial_{y_{12}}-y_{13}\partial_{y_{13}}-y_{14}\partial_{y_{14}}\nonumber\\
			&\:+\mathop{\textstyle\sum}\limits_{i=2}^4(x_{i,2}\partial_{x_{i,2}}-x_{i,-2}\partial_{x_{i,-2}}-x_{i,1}\partial_{x_{i,1}}+x_{i,-1}\partial_{x_{i,-1}})-2x_{1,1}\partial_{x_{1,1}}+\lambda_5,\label{3.101}\\
			\pi_{\lambda}(\alpha_6)=&\:-2y_1\partial_{y_{1}}-\mathop{\textstyle\sum}\limits_{i=2}^{9}y_i\partial_{y_i}-y_{11}\partial_{y_{11}}-y_{12}\partial_{y_{12}}+x_{3,1}\partial_{x_{3,1}}+x_{3,-2}\partial_{x_{3,-2}}\nonumber\\
			&\:+x_{4,1}\partial_{x_{4,1}}+x_{4,-2}\partial_{x_{4,-2}}+x_{4,-3}\partial_{x_{4,-3}}+x_{4,-4}\partial_{x_{4,-4}}\nonumber\\
			&\:+x_{3,-3}\partial_{x_{3,-3}}+x_{2,1}\partial_{x_{2,1}}+x_{2,-2}\partial_{x_{2,-2}}+x_{1,1}\partial_{x_{1,1}}+\lambda_6.\label{3.102}
		\end{align}
			\item The simple negative root vectors.
		\begin{align}
			\pi_{\lambda}(E_{-\alpha_1})=&\:-y_7\partial_{y_4}-y_8\partial_{y_6}-y_{11}\partial_{y_9}-y_{13}\partial_{y_{10}}+x_{4,4}\mathop{\textstyle\sum}\limits_{|r|=1}^4x_{4,r}\partial_{x_{4,r}}\nonumber\\
			&\:-\mathop{\textstyle\sum}\limits_{s=1}^4x_{4,s}x_{4,-s}\partial_{x_{4,-4}}+\mathop{\textstyle\sum}\limits_{|t|=1}^3(x_{4,t}-x_{4,4}x_{3,t})\partial_{x_{3,t}}-\lambda_1x_{4,4},\label{3.85}\\
			\pi_{\lambda}(E_{-\alpha_2})=&\:y_5\partial_{y_3}+y_6\partial_{y_4}+y_8\partial_{y_7}-y_{16}\partial_{y_{15}}-\mathop{\textstyle\sum}\limits_{i=2}^4(x_{i,-1}\partial_{x_{i,2}}-x_{i,-2}\partial_{x_{i,1}})\nonumber\\
			&\:+x_{1,-1}^2\partial_{x_{1,-1}}-\lambda_2x_{1,-1},\\
			\pi_{\lambda}(E_{-\alpha_3})=&\:-y_4\partial_{y_3}-y_{6}\partial_{y_5}-y_{12}\partial_{y_{11}}-y_{14}\partial_{y_{13}}-x_{4,3}\partial_{x_{4,4}}+x_{4,-4}\partial_{x_{4,-3}}\nonumber\\
			&\:+x_{3,3}\mathop{\textstyle\sum}\limits_{|r|=1}^3x_{3,r}\partial_{x_{3,r}}-\mathop{\textstyle\sum}\limits_{s=1}^3x_{3,s}x_{3,-s}\partial_{x_{3,-3}}+\mathop{\textstyle\sum}\limits_{|t|=1}^2(x_{3,t}-x_{3,3}x_{2,t})\partial_{x_{2,t}}\nonumber\\
			&\:-\lambda_3x_{3,3},\\
			\pi_{\lambda}(E_{-\alpha_4})=&\:y_3\partial_{y_2}+y_{9}\partial_{y_6}+y_{11}\partial_{y_8}-y_{15}\partial_{y_{14}}-x_{4,2}\partial_{x_{4,3}}+x_{4,-3}\partial_{x_{4,-2}}\nonumber\\
			&\:-x_{3,2}\partial_{x_{3,3}}+x_{3,-3}\partial_{x_{3,-2}}+x_{2,2}\mathop{\textstyle\sum}\limits_{|r|=1}^2x_{2,r}\partial_{x_{2,r}}-\mathop{\textstyle\sum}\limits_{s=1}^2x_{2,s}x_{2,-s}\partial_{x_{2,-2}}\nonumber\\
			&\:+(x_{2,1}-x_{2,2}x_{1,1})\partial_{x_{1,1}}+(x_{2,-1}-x_{2,2}x_{1,-1})\partial_{x_{1,-1}}-\lambda_4x_{2,2},\\
			\pi_{\lambda}(E_{-\alpha_5})=&\:y_2\partial_{y_1}-y_{10}\partial_{y_9}-y_{13}\partial_{y_{11}}-y_{14}\partial_{y_{12}}-\mathop{\textstyle\sum}\limits_{i=2}^4(x_{i,1}\partial_{x_{i,2}}-x_{i,-2}\partial_{x_{i,-1}})\nonumber\\
			&\:+x_{1,1}^2\partial_{x_{1,1}}-\lambda_5x_{1,1},\label{3.89}
		\end{align}
		\begin{align}
			\pi_{\lambda}(E_{-\alpha_6})=&\:y_1\mathop{\textstyle\sum}\limits_{i=1}^{16}y_i\partial_{y_{i}}-\zeta_1\partial_{y_{10}}-\zeta_2\partial_{y_{13}}-\zeta_3\partial_{y_{14}}-\zeta_4\partial_{y_{15}}-\zeta_{10}\partial_{y_{16}}\nonumber\\
			&\:-(y_1x_{4,1}+y_2x_{4,2}+y_3x_{4,3}-y_4x_{4,4}+y_7)\partial_{x_{4,1}}\nonumber\\
			&\:+(y_2x_{4,-1}-y_1x_{4,-2}-y_5x_{4,3}+y_6x_{4,4}-y_8)\partial_{x_{4,-2}}\nonumber\\
			&\:+(y_3x_{4,-1}+y_5x_{4,2}-y_1x_{4,-3}+y_9x_{4,4}-y_{11})\partial_{x_{4,-3}}\nonumber\\
			&\:-(y_4x_{4,-1}+y_6x_{4,2}+y_9x_{4,3}+y_1x_{4,-4}-y_{12})\partial_{x_{4,-4}}\nonumber\\
			&\:-(y_1x_{3,1}+y_2x_{3,2}+y_3x_{3,3}-y_4)\partial_{x_{3,1}}\nonumber\\
			&\:+(y_2x_{3,-1}-y_1x_{3,-2}-y_5x_{3,3}+y_6)\partial_{x_{3,-2}}\nonumber\\
			&\:+(y_3x_{3,-1}+y_5x_{3,2}-y_1x_{3,-3}+y_9)\partial_{x_{3,-3}}\nonumber\\
			&\:-(y_1x_{2,1}+y_2x_{2,2}+y_3)\partial_{x_{2,1}}+(y_2x_{2,-1}-y_1x_{2,-2}-y_5)\partial_{x_{2,-2}}\nonumber\\
			&\:-(y_1x_{1,1}+y_2)\partial_{x_{1,1}}-\lambda_6y_1.\label{3.90}
		\end{align}		
	\end{enumerate}
\end{theorem}
	\begin{proof}
		Here we only compute the formulas for the triple $\{E_{\alpha_6},\alpha_6,E_{-\alpha_6}\}$, and the remaining formulas can be obtained by the same method.
		
		By (\ref{2.43}) we have $E_{\alpha_6}=\xi_1$ and $E_{-\alpha_6}=\eta_1$. Firstly, Eq. (\ref{r3}) in Theorem \ref{thm3.1} gives $\pi_{\lambda}(E_{\alpha_6})=\partial_{y_1}$, which is exactly (\ref{3.96}). 
		By (\ref{zeta1})-(\ref{zeta10}) and (\ref{2.65}), we have
		\begin{eqnarray}
			\sigma(E_{-\alpha_6})=y_1\mathop{\textstyle\sum}\limits_{i=1}^{16}y_i\partial_{y_{i}}-\zeta_1\partial_{y_{10}}-\zeta_2\partial_{y_{13}}-\zeta_3\partial_{y_{14}}-\zeta_4\partial_{y_{15}}-\zeta_{10}\partial_{y_{16}}.\label{e61}
		\end{eqnarray}
		Eqs. (\ref{248}) and (\ref{zi}) give the explicit expression of $Z_1=\begin{pmatrix}
			P_1&Q_1\\R_1&-P_1^t
		\end{pmatrix}$ as
		\begin{eqnarray}
			P_1=\begin{pmatrix}
				-\frac{y_1}{2}&0&0&0&0\\
				-y_2&\frac{y_1}{2}&0&0&0\\
				-y_3&0&\frac{y_1}{2}&0&0\\
				y_4&0&0&\frac{y_1}{2}&0\\
				-y_7&0&0&0&\frac{y_1}{2}
			\end{pmatrix},\;
			Q_1=\begin{pmatrix}
				0&0&0&0&0\\
				0&0&y_5&-y_6&y_8\\
				0&-y_5&0&-y_9&y_{11}\\
				0&y_6&y_9&0&-y_{12}\\
				0&-y_8&-y_{11}&y_{12}&0
			\end{pmatrix},\; R_1=0.\label{PQR1}
		\end{eqnarray}
		Moreover, the strict upper triangular part of $LP_1L^{-1}$ is zero since $L$, $P_1$ and $L^{-1}$ are all lower triangular, hence
		\begin{eqnarray}
			U_1=\operatorname{diag}P_1=\operatorname{diag}(-\frac{y_1}{2},\frac{y_1}{2},\frac{y_1}{2},\frac{y_1}{2},\frac{y_1}{2}) \label{U1}
		\end{eqnarray}
		By (\ref{lam6}), the multiplication term of $\pi_{\lambda}(E_{-\alpha_6})$ is
		\begin{eqnarray}
			-\frac{\kappa}{2}y_1+\sum_{r=1}^{5}\mu_r(U_1)_{r,r}=-\lambda_6y_1. \label{e62}
		\end{eqnarray} 
		Applying the expressions (\ref{2.8}), (\ref{2.9}), (\ref{PQR1}) and (\ref{U1}) to (\ref{r4}) with $i=1$ we obtain the nonzero polynomial coefficients of $\partial_{x_{r,s}}$ in the $\pi_{\lambda}(\eta_1)$ are as follows
		\begin{center}
			\renewcommand{\arraystretch}{1.15}
			\begin{tabular}{@{}ll@{}}
				\toprule
				$(r,s)$ & the coefficient of $\partial_{x_{r,s}}$\\
				\midrule
				$(1,1)$ & $-y_1x_{1,1}-y_2$\\
				$(2,1)$ & $-y_1x_{2,1}-y_2x_{2,2}-y_3$\\
				$(2,-2)$ & $y_2x_{2,-1}-y_1x_{2,-2}-y_5$ \\
				$(3,1)$ & $-y_1x_{3,1}-y_2x_{3,2}-y_3x_{3,3}+y_4$\\
				$(3,-2)$ & $y_2x_{3,-1}-y_1x_{3,-2}-y_5x_{3,3}+y_6$\\
				$(3,-3)$ & $y_3x_{3,-1}+y_5x_{3,2}-y_1x_{3,-3}+y_9$\\
				$(4,1)$ &  $-y_1x_{4,1}-y_2x_{4,2}-y_3x_{4,3}+y_4x_{4,4}-y_7$\\
				$(4,-2)$ & $y_2x_{4,-1}-y_1x_{4,-2}-y_5x_{4,3}+y_6x_{4,4}-y_8$\\
				$(4,-3)$ & $y_3x_{4,-1}+y_5x_{4,2}-y_1x_{4,-3}+y_9x_{4,4}-y_{11}$\\
				$(4,-4)$ & $-y_4x_{4,-1}-y_6x_{4,2}-y_9x_{4,3}-y_1x_{4,-4}+y_{12}$\\
				\bottomrule
			\end{tabular}
		\end{center}
		\vspace{3mm}
		
		The above tabular together with (\ref{e61}) and (\ref{e62}) give the expression (\ref{3.90}) of $\pi_{\lambda}(E_{-\alpha_6})$. Finally, the expression (\ref{3.102}) of $\pi_{\lambda}(\alpha_6)$ is derived from (\ref{3.96}), (\ref{3.90}) and the equality
		\begin{eqnarray}
			\pi_{\lambda}(\alpha_6)=-[\pi_{\lambda}(E_{\alpha_6}),\pi_{\lambda}(E_{-\alpha_6})].
		\end{eqnarray}
		This completes the proof.
	\end{proof}

From now on, we write $\mathscr{C}_{\lambda}$ for $\mathscr{C}$ with the action of $\mathscr{G}^{E_6}$ given by $\pi_\lambda$, and we simply write
\begin{eqnarray}
	u(f)=\pi_{\lambda}(u)(f)
\end{eqnarray}
for $u\in\mathscr{G}^{E_6}$ and $f\in\mathscr{C}_{\lambda}$.
The following lemma is an useful observation of the representation formulas.
\begin{lemma}\label{lem3.3}
	For arbitrary $\lambda\in \mathscr{H}^*$ and every $u\in\mathscr{G}^{E_6}$ there exists
	\begin{eqnarray}
		u^{D}\in\mathscr{W}_{\mathscr{C}}\quad\text{and}\quad u^P_i\in\mathscr{C}_{\lambda}\;\;(i\in\overline{1,6})
	\end{eqnarray}
	independent of $\lambda$ such that
	\begin{eqnarray}
		\pi_{\lambda}(u)=u^D+\sum_{i=1}^6\lambda_iu^P_i.\label{3.31}
	\end{eqnarray}
	Moreover, 
	\begin{enumerate}[(i)]
		\item If $u$ is a positive root vector, then each $u_{i}^{P}=0$, thus $\pi_\lambda(u)\in\mathscr{W}_{\mathscr{C}}$.
		\item If $u$ is a negative root vector with root $-\sum_{i=1}^6n_i\alpha_i$, then $u_{i}^{P}=0$ whenever $n_i=0$.
	\end{enumerate}
\end{lemma}
\begin{proof}
	For $d_1,d_2\in\mathscr{W}_\mathscr{C}$ and $p_i,q_i\in\mathscr{C}$ for $i\in\overline{1,6}$, we have
	\begin{eqnarray}
		[d_1+\sum_{i=1}^6\lambda_ip_i,d_2+\sum_{i=1}^6\lambda_iq_i]=[d_1,d_2]+\sum_{i=1}^6 \lambda_i(d_1(q_i)-d_2(p_i)). \label{3.32}
	\end{eqnarray}
	Theorem \ref{thm3.2} gives the formulas of all simple root vectors and coroots, which coincide with (\ref{3.31}) and satisfy assertions (i)-(ii). Since these elements generate $\mathscr{G}^{E_6}$, Eq. (\ref{3.32}) guarantees that $\pi_{\lambda}(u)$ must be of the form (\ref{3.31}). 
	
	Since the positive nilpotent subalgebra of $\mathscr{G}^{E_6}$ is generated by these simple positive root vectors $E_{\alpha_i}$'s, which have no multiplication terms by (\ref{3.91})-(\ref{3.96}), Eq. (\ref{3.32}) implies the property (i) for all positive root vectors.
	
	Suppose that $u$ is a negative root vector with root $\alpha=-\sum_{i=1}^6n_i\alpha_i$ with $n_i\in\mathbb{N}$. We prove assertion (ii) by induction on $n=\sum_{i=1}^6n_i$. If $n=1$, then $u$ is a simple negative root vector, the explicit formulas (\ref{3.85})-(\ref{3.90}) in Theorem \ref{thm3.2} give (ii). Assume $n>1$. Then we can choose $i_0\in\overline{1,6}$ such that $\alpha_{i_0}-\sum_{i=1}^6n_i\alpha_i$ is a negative root, and there exists a corresponding negative root vector $u_0$ such that 
	\begin{align}
		[E_{-\alpha_{i_0}},u_0]=u.
	\end{align}
	By induction assumption on $u_0$ and (\ref{3.32}), we derive that $u$ also satisfies property (ii).
\end{proof}

The representation formulas above are the main tools for our subsequent analysis. Lemma \ref{lem3.3} says that the actions of the positive root vectors are first-order differential operators with polynomial coefficients independent of $(\lambda_1,\cdots,\lambda_6)$, and the negative root vectors may have a nontrivial multiplication term.

\begin{proposition}{\label{prop3.4}}
	For any $\lambda\in \mathscr{H}^*$, every monomial in $\mathscr{C}_{\lambda}$ is a weight vector. Every singular vector in $\mathscr{C}_{\lambda}$ is a nonzero scalar multiple of $1$. Moreover, the cyclic submodule $U(\mathscr{G}^{E_6})(1)$ is a highest-weight $\mathscr{G}^{E_6}$-module with highest weight $\lambda$. 
\end{proposition}
\begin{proof}
	The first assertion is derived by applying $(\ref{3.97})$-$(\ref{3.102})$ to each monomial.
	Suppose $f\in\mathscr{C}_{\lambda}$ is a singular vector. By (\ref{r3}) we have
	\begin{eqnarray}
		\xi_k(f)=\partial_{y_k}(f)=0\quad(k\in\overline{1,16}),
	\end{eqnarray}
	hence $f$ is independent of all the variables $\{y_i\mid i\in\overline{1,16}\}$. Next applying $A_{5,i}$ with $i=\pm1,\pm2,\pm3,\pm4$ to $f$ yields
	\begin{eqnarray}
		A_{5,i}(f)=\partial_{x_{4,i}}(f)=0,\label{3.114}
	\end{eqnarray}
	hence $f$ is independent of $x_{4,\pm1},x_{4,\pm2},x_{4,\pm3},x_{4,\pm4}$. Similarly, using $A_{4,i}$ ($i=\pm1,\pm2,\pm3$), $A_{3,i}$ ($i=\pm1,\pm2$), $A_{2,\pm1}$ successively, we eliminate dependence on $x_{3,i}$ ($i=\pm1,\pm2,\pm3$), $x_{2,i}$ ($i=\pm1,\pm2$) and finally $x_{1,\pm1}$. Hence $f$ must be a constant.
	
	A direct evaluation of the Cartan actions (\ref{3.97})-(\ref{3.102}) on the constant $1$ gives the weight $\lambda=\sum_{i=1}^6\lambda_i\omega_i$,
	and the second claim of this proposition follows.
\end{proof}

\subsection{Restricted dual modules}

We now are going to show that the constructed module $\mathscr{C}_{\lambda}$ is isomorphic to the restricted dual Verma module. For simplicity, for $u\in\mathscr{G}^{E_6}$ we write $u$ for $\pi_{\lambda}(u)$ when it acts on $\mathscr{C}_{\lambda}$ in this subsection.
Recall the definition of $restricted$ $dual$ $modules$ (see e.g., \cite{Hum2}) here. Let $N$ be a weight $\mathscr{G}^{E_6}$-module with finite-dimensional weight spaces, suppose $N$ has weight space decomposition
\begin{eqnarray}
	N=\bigoplus_{\mu\in \mathscr{H}^*}N_{\mu}.
\end{eqnarray}
Denote by $N_{\mu}^*$ the dual space of $N_{\mu}$ for each $\mu$. This space can be naturally identified with the subspace of the full dual space $N^*$ consisting of those linear functions that vanish on all weight spaces $N_{\nu}$ with $\nu\neq\mu$. Let $\tau$ be the standard anti-involution on $\mathscr{G}^{E_6}$ determined by the generators
\begin{eqnarray}
	\tau(E_{\alpha_i})=E_{-\alpha_i},\;\tau(E_{-\alpha_i})=E_{\alpha_i},\;\tau(h)=h,\quad (i\in\overline{1,6},\quad h\in \mathscr{H}), 
\end{eqnarray}
and the relation $\tau([u_1,u_2])=-[\tau(u_1),\tau(u_2)]$ for $u_1,u_2\in\mathscr{G}^{E_6}$. Hence for any $\gamma\in\Phi_{E_6}^+$,
\begin{eqnarray}
	\tau(E_\gamma)=(-1)^{\operatorname{ht}(\gamma)-1}E_{-\gamma},\;\;\tau(E_{-\gamma})=(-1)^{\operatorname{ht}(\gamma)-1}E_{\gamma}.
\end{eqnarray}
This map $\tau$ extends canonically to an anti‑automorphism of $U(\mathscr{G}^{E_6})$.
Set
\begin{eqnarray}
	N^{\vee}=\bigoplus_{\mu\in \mathscr{H}^*}N_{\mu}^{*},\label{dual1}
\end{eqnarray}
and define the action of $\mathscr{G}^{E_6}$ on $N^{\vee}$ by
\begin{eqnarray}
	(u\cdot f)(m)=f(\tau(u)\cdot m),\qquad u\in\mathscr{G}^{E_6},\quad f\in N^{\vee},\quad m\in N. \label{dual3}
\end{eqnarray}
Then $N^{\vee}$ becomes a $\mathscr{G}^{E_6}$-module.

Let $M(\lambda)$ be the Verma module over $\mathscr{G}^{E_6}$ generated by a highest weight vector $w_{\lambda}$ of weight $\lambda$. Choose 
\begin{eqnarray}
	\eta_1, \eta_2, \cdots, \eta_{16}; A_{-4,5}, A_{-3,5}, \cdots, A_{4,5}; A_{-3,4}, \cdots, A_{3,4}; A_{-2,3}, \cdots, A_{2,3}; A_{-1,2}, A_{1,2} \label{order}
\end{eqnarray}
as an ordered basis of the nilpotent subalgebra consisting of negative root vectors of $\mathscr{G}^{E_6}$.
For $a=(a_1,\cdots,a_{16})\in\mathbb{N}^{16}$ and $b=(b_{i,j}\mid 1\leqslant |j|\leqslant i\leqslant 4)\in\mathbb{N}^{20}$, define
\begin{eqnarray}
	\eta^{(a,b)}=\eta_1^{a_1}\eta_2^{a_2}\cdots\eta_{16}^{a_{16}}\big(\prod_{|j|=1}^4A_{j,5}^{b_{4,j}}\big)\big(\prod_{|j|=1}^3A_{j,4}^{b_{3,j}}\big)\big(\prod_{|j|=1}^2A_{j,3}^{b_{2,j}}\big)A_{-1,2}^{b_{1,-1}}A_{1,2}^{b_{1,1}},
\end{eqnarray}
where the product respects the ordering in (\ref{order}).
By the PBW theorem, the set $\{\eta^{(a,b)}w_{\lambda}\mid a\in\mathbb{N}^{16}, b\in\mathbb{N}^{20}\}$ forms a basis of $M(\lambda)$. We write $\operatorname{wt}(a,b)$ for the weight of $\eta^{(a,b)}w_{\lambda}$, then $\lambda-\operatorname{wt}(a,b)\in Q_{E_6}^+$. We also introduce the monomials
\begin{eqnarray}
	Y^aX^b=\prod_{i=1}^{16}y_i^{a_i}\prod_{1\leqslant |j|\leqslant i\leqslant 4}x_{i,j}^{b_{i,j}},
\end{eqnarray}
and denote by $(Y^aX^b)^*$ the corresponding dual element in $\mathscr{C}_{\lambda}^{\vee}$ as $(Y^aX^b)^*$. For convenience, we denote by $\epsilon_i$ ($i\in\overline{1,16}$) the standard tuples in $\mathbb{N}^{16}$ with
$1$ in the $i$-coordinate and $0$ elsewhere. Similarly, we write $\epsilon_{i,j}\in\mathbb{N}^{20}$ for the tuple that has $1$ in the $(i,j)$-entry and $0$ otherwise.

We now define a linear map $\varphi:M(\lambda)\longrightarrow\mathscr{C}_{\lambda}^{\vee}$ by
\begin{align}
	\varphi(\eta^{(a,b)}w_{\lambda})=(-1)^{\operatorname{ht}(\lambda-\operatorname{wt}(a,b))-||a||}a!b!(Y^{a}X^b)^*,
\end{align}
here $||a||=a_1+\cdots+a_{16}$.
One readily verifies that $\varphi$ is a linear isomorphism.

\begin{lemma}\label{ch}
	The map $\varphi$ sends weight vectors of $M(\lambda)$ to the weight vectors of $\mathscr{C}_{\lambda}^{\vee}$ and preserve the weight. 
\end{lemma}
\begin{proof}
	A straightforward calculation shows that $\eta^{(a,b)}w_{\lambda}$ is a weight vector with weight $\lambda-\sum_{i=1}^6\nu_i\alpha_i$, where the non‑negative integers $\nu_i$'s are given by the following expressions:
	\begin{align}
		\nu_1=&\:a_7+a_8+\sum_{i=11}^{16}a_i+\sum_{|s|=1}^4b_{4,s}, \label{ch1}        \\
		\nu_2=&\:a_5+a_6+\sum_{i=8}^{16}a_i+a_{16}+\sum_{1\leqslant -s\leqslant r\leqslant 4}b_{r,s} ,          \\
		\nu_3=&\:a_4+\sum_{i=6}^{16}a_{i}+a_{12}+\sum_{i=14}^{16}a_{i}+\sum_{|s|=1}^3(b_{4,s}+b_{3,s})+2b_{4,-4},     \\
		\nu_4=&\:\sum_{i=3}^{16}a_i+\sum_{i=9}^{16}a_i+a_{15}+a_{16}+\sum_{r=2,3,4\atop |s|=1,2}b_{r,s}+2(b_{4,-3}+b_{4,-4}+b_{3,-3}),           \\
		\nu_5=&\:\sum_{i=2}^{16}a_i+a_{10}+\sum_{i=13}^{16}a_{i}+\sum_{2\leqslant -s\leqslant r\leqslant 4}b_{r,s}+\sum_{r=1}^4b_{r,1},                    
	\end{align}
	and
	\begin{align}
		\nu_6=\sum_{i=1}^{16}a_i.  
	\end{align}
	
	On the other hand, the weight space decomposition of $\mathscr{C}_{\lambda}$ is obtained from the action of the coroots given in (\ref{3.97})-(\ref{3.102}). A direct verification shows that $Y^aX^b$ is a weight vector of weight $\mu=\sum_{i=1}^6\mu_i\omega_i$, where
	\begin{align}
		\mu_1=&\:\lambda_1+a_4+a_6-a_7-a_8+a_9+a_{10}-a_{11}-a_{13}\nonumber\\
		&\:-\sum\limits_{|r|=1}^3(b_{4,r}-b_{3,r})-2b_{4,4}, \label{ch3}\\
		\mu_2=&\:\lambda_2+a_3+a_4-a_5-a_6+a_7-a_8+a_{15}-a_{16}\nonumber\\
		&\:+\sum_{i=2}^4(b_{i,2}-b_{i,-2}+b_{i,1}-b_{i,-1})-2b_{1,-1},\\
		\mu_3=&\:\lambda_3+a_3-a_4+a_{5}-a_6+a_{11}-a_{12}+a_{13}-a_{14}\nonumber\\
		&\:+b_{4,4}-b_{4,-4}-b_{4,3}+b_{4,-3}-\sum_{|r|=1}^2(b_{3,r}-b_{2,r})-2b_{3,3}, \\
		\mu_4=&\:\lambda_4+a_2-a_{3}+a_6+a_{8}-a_9-a_{11}+a_{14}-a_{15}\nonumber\\
		&\:+\sum_{i=3}^4(b_{i,3}-b_{i,-3}-b_{i,2}+b_{i,-2})-\sum_{r\in\{1,-1\}}(b_{2,r}-b_{1,r})-2b_{2,2},
	\end{align}\begin{align}
		\mu_5=&\:\lambda_5+a_1-a_{2}+a_9-a_{10}+a_{11}+a_{12}-a_{13}-a_{14}\nonumber\\
		&\:+\sum_{i=2}^4(b_{i,2}-b_{i,-2}-b_{i,1}+b_{i,-1})-2b_{1,1},\\
		\mu_6=&\:\lambda_6-2a_1-\sum\limits_{i=2}^{9}a_i-a_{11}-a_{12}+b_{3,1}+b_{3,-2}+b_{3,-3}\nonumber\\
		&\:+b_{2,1}+b_{2,-2}+b_{1,1}+b_{4,1}+b_{4,-2}+b_{4,-3}+b_{4,-4}. \label{ch4}
	\end{align}
	Note that the twisted action by $\tau$ does not change $\mathscr{H}$, so $(Y^aX^b)^*$ is also a weight vector of weight $\mu$. From (\ref{ch1})-(\ref{ch4}), we have
	\begin{align}
		A^{E_6}\begin{pmatrix}
			\nu_1\\\vdots\\\nu_6
		\end{pmatrix}=\begin{pmatrix}
			\lambda_1-\mu_1\\\vdots\\\lambda_6-\mu_6
		\end{pmatrix}
	\end{align}
	(cf. (\ref{CM})), which implies $\lambda-\sum_{i=1}^6\nu_i\alpha_i=\mu$. This yields the required statement.
\end{proof}

Since $\varphi$ is bijective and sends a weight basis of $M(\lambda)$ to a weight basis of $\mathscr{C}_{\lambda}^{\vee}$, it follows that the formal characters of the $\mathscr{G}^{E_6}$-module $\mathscr{C}_{\lambda}$ and the Verma module $M(\lambda)$ coincide, that is,
\begin{eqnarray}
	\text{ch}(\mathscr{C}_{\lambda})=\text{ch}(M(\lambda)). \label{char}
\end{eqnarray}
Furthermore, we have
\begin{theorem}\label{thiso}
	As a module over $\mathscr{G}^{E_6}$, the restricted dual module $\mathscr{C}_{\lambda}^{\vee}$ is isomorphic to the Verma module $M(\lambda)$.
\end{theorem}
\begin{proof}
	It remains to show that the linear isomorphism $\varphi$ intertwines the $\mathscr{G}^{E_6}$-action. Because the Chevalley generators $\{E_{\alpha_i},E_{-\alpha_i},\alpha_i\mid i\in\overline{1,6}\}$ generate $\mathscr{G}^{E_6}$ as a Lie algebra, it suffices to verify $\varphi$ commutes with the actions of each of these generators.
	
	By the definition (\ref{dual1})-(\ref{dual3}) and the proof of Lemma \ref{ch}, $\eta^{(a,b)}w_{\lambda}$ and $(Y^{a}X^b)^*$ are weight vectors of the same weight $\mu$. For any $h\in \mathscr{H}$, we have
	\begin{align}
		&\:\varphi(h\eta^{(a,b)}w_{\lambda})=\varphi(\mu(h)\eta^{(a,b)}w_{\lambda})=\mu(h)\varphi(\eta^{(a,b)}w_{\lambda})\nonumber\\
		=&\:(-1)^{\operatorname{ht}(\lambda-\operatorname{wt}(a,b))-||a||}\mu(h)a!b!(Y^{a}X^b)^*\nonumber\\
		=&(-1)^{\operatorname{ht}(\lambda-\operatorname{wt}(a,b))-||a||}a!b!h(Y^{a}X^b)^*=h\varphi(\eta^{(a,b)}w_{\lambda}).
	\end{align}
	So $\varphi$ commutes with the action of $\mathscr{H}$.
	
	Next we consider the negative simple root vectors. According to (\ref{2.37})-(\ref{2.28}), we have
	\begin{eqnarray}
		[E_{-\alpha_1},\eta_4]=-\eta_7,\;[E_{-\alpha_1},\eta_6]=-\eta_8,\;[E_{-\alpha_1},\eta_{9}]=-\eta_{11},\;[E_{-\alpha_1},\eta_{10}]=-\eta_{13},
	\end{eqnarray}
	and $E_{-\alpha_1}$ commutes with other negative root vectors before $A_{-4,5}$ in ordering given in (\ref{order}). Hence we obtain
	\begin{align}
		E_{-\alpha_1}\eta^{(a,b)}w_{\lambda}=&\:(-a_4\eta^{(a+\epsilon_7-\epsilon_4,b)}-a_6\eta^{(a+\epsilon_8-\epsilon_6,b)}-a_{9}\eta^{(a+\epsilon_{11}-\epsilon_{9},b)}\nonumber\\
		&\:-a_{10}\eta^{(a+\epsilon_{13}-\epsilon_{10},b)}-\eta^{(a,b+\epsilon_{4,4})})w_{\lambda}.\label{-al1}
	\end{align}
	Here and in the following formulas, any term in which an exponent becomes negative is understood to be zero. Analogously, we have
	\begin{align}
		E_{-\alpha_2}\eta^{(a,b)}w_{\lambda}=&\:(a_3\eta^{(a+\epsilon_5-\epsilon_3,b)}+a_4\eta^{(a+\epsilon_6-\epsilon_4,b)}+a_7\eta^{(a+\epsilon_8-\epsilon_7,b)}\nonumber  \\
		&\:-a_{15}\eta^{(a+\epsilon_{16}-\epsilon_{15},b)}-b_{4,2}\eta^{(a,b+\epsilon_{4,-1}-\epsilon_{4,2})}+b_{4,1}\eta^{(a,b+\epsilon_{4,-2}-\epsilon_{4,1})}\nonumber\\
		&\:-b_{3,2}\eta^{(a,b+\epsilon_{3,-1}-\epsilon_{3,2})}+b_{3,1}\eta^{(a,b+\epsilon_{3,-2}-\epsilon_{3,1})}-b_{2,2}\eta^{(a,b+\epsilon_{2,-1}-\epsilon_{2,2})}\nonumber\\
		&\:+b_{2,1}\eta^{(a,b+\epsilon_{2,-2}-\epsilon_{2,1})}-\eta^{(a,b+\epsilon_{1,-1})})w_{\lambda}, \label{-al2}
	\end{align}
	\begin{align}
		E_{-\alpha_3}\eta^{(a,b)}w_{\lambda}=&\: (-a_3\eta^{(a+\epsilon_4-\epsilon_3,b)}-a_{5}\eta^{(a+\epsilon_6-\epsilon_5,b)} -a_{11}\eta^{(a+\epsilon_{12}-\epsilon_{11},b)}-a_{13}\eta^{(a+\epsilon_{14}-\epsilon_{13},b)}\nonumber\\
		&\:-b_{4,4}\eta^{(a,b+\epsilon_{4,3}-\epsilon_{4,4})}+b_{4,-3}\eta^{(a,b+\epsilon_{4,-4}-\epsilon_{4,-3})}-\eta^{(a,b+\epsilon_{3,3})})w_{\lambda}, 
	\end{align}
	\begin{align}
		E_{-\alpha_4}\eta^{(a,b)}w_{\lambda}=&\: (a_2\eta^{(a+\epsilon_3-\epsilon_2,b)}+a_{6}\eta^{(a+\epsilon_9-\epsilon_6,b)} +a_{8}\eta^{(a+\epsilon_{11}-\epsilon_{8},b)}-a_{14}\eta^{(a+\epsilon_{15}-\epsilon_{14},b)}\nonumber\\ 
		&\:-b_{4,3}\eta^{(a,b+\epsilon_{4,2}-\epsilon_{4,3})}+b_{4,-2}\eta^{(a,b+\epsilon_{4,-3}-\epsilon_{4,-2})}-b_{3,3}\eta^{(a,b+\epsilon_{3,2}-\epsilon_{3,3})} \nonumber\\
		&\:+b_{3,-2}\eta^{(a,b+\epsilon_{3,-3}-\epsilon_{3,-2})}-\eta^{(a,b+\epsilon_{2,2})})w_{\lambda}, 
	\end{align}
	\begin{align}
		E_{-\alpha_5}\eta^{(a,b)}w_{\lambda}=&\:(a_1\eta^{(a+\epsilon_2-\epsilon_1,b)}-a_{9}\eta^{(a+\epsilon_{10}-\epsilon_9,b)}-a_{11}\eta^{(a+\epsilon_{13}-\epsilon_{11},b)}\nonumber  \\
		&\:-a_{12}\eta^{(a+\epsilon_{14}-\epsilon_{12},b)}-b_{4,2}\eta^{(a,b+\epsilon_{4,1}-\epsilon_{4,2})}+b_{4,-1}\eta^{(a,b+\epsilon_{4,-2}-\epsilon_{4,-1})}\nonumber\\
		&\:-b_{3,2}\eta^{(a,b+\epsilon_{3,1}-\epsilon_{3,2})}+b_{3,-1}\eta^{(a,b+\epsilon_{3,-2}-\epsilon_{3,-1})}-b_{2,2}\eta^{(a,b+\epsilon_{2,1}-\epsilon_{2,2})}\nonumber\\
		&\:+b_{2,-1}\eta^{(a,b+\epsilon_{2,-2}-\epsilon_{2,-1})}-\eta^{(a,b+\epsilon_{1,1})})w_{\lambda}, 
	\end{align}
	and
	\begin{eqnarray}
		E_{-\alpha_6}\cdot \eta^{(a,b)}w_{\lambda}=\eta^{(a+\epsilon_1,b)}w_{\lambda}. \label{-al6}
	\end{eqnarray}
	On the other side, A direct computation by (\ref{3.91}) yields
	\begin{align}
		&\:E_{-\alpha_1}((Y^aX^b)^*)(Y^{a'}X^{b'})=(Y^aX^b)^*(E_{\alpha_1}(Y^{a'}X^{b'}))\nonumber\\
		=&\:(Y^aX^b)^*(a'_7Y^{a'-\epsilon_7+\epsilon_4}X^{b'}+a'_8Y^{a'-\epsilon_8+\epsilon_6}X^{b'}+a'_{11}Y^{a'-\epsilon_{11}+\epsilon_{9}}X^{b'}\nonumber\\
		&\:+a'_{13}Y^{a'-\epsilon_{13}+\epsilon_{10}}X^{b'}+b'_{4,4}Y^{a'}X^{b'-\epsilon_{4,4}})\nonumber\\
		=&\:a_7'\delta_{a,a'-\epsilon_7+\epsilon_4}\delta_{b,b'}+a_8'\delta_{a,a'-\epsilon_8+\epsilon_6}\delta_{b,b'}+a_{11}'\delta_{a,a'-\epsilon_{11}+\epsilon_9}\delta_{b,b'}\nonumber\\
		&\:+a_{13}'\delta_{a,a'-\epsilon_{13}+\epsilon_{10}}\delta_{b,b'}+b_{4,4}'\delta_{a,a'}\delta_{b,b'-\epsilon_{4,4}}\nonumber\\
		=&\:(a_7+1)\delta_{a+\epsilon_7-\epsilon_4,a'}\delta_{b,b'}+(a_8+1)\delta_{a+\epsilon_8-\epsilon_6,a'}\delta_{b,b'}+(a_{11}+1)\delta_{a+\epsilon_{11}-\epsilon_9,a'}\delta_{b,b'}\nonumber\\
		&\:+(a_{13}+1)\delta_{a+\epsilon_{13}-\epsilon_{10},a'}\delta_{b,b'}+(b_{4,4}+1)\delta_{a,a'}\delta_{b+\epsilon_{4,4},b'}.
	\end{align}
	Thus
	\begin{align}
		E_{-\alpha_1}((Y^aX^b)^*)=&\:(a_7+1)(Y^{a+\epsilon_7-\epsilon_4}X^{b})^*+(a_8+1)(Y^{a+\epsilon_8-\epsilon_6}X^b)^*\nonumber\\
		&\:+(a_{11}+1)(Y^{a+\epsilon_{11}-\epsilon_9}X^{b})^*+(a_{13}+1)(Y^{a+\epsilon_{13}-\epsilon_{10}}X^{b})^*\nonumber\\
		&\:+(b_{4,4}+1)(Y^{a}X^{b+\epsilon_{4,4}})^*.\label{al1}
	\end{align}
	Similarly, we have
	\begin{align}
		E_{-\alpha_2}((Y^aX^b)^*)=&\:-(a_5+1)(Y^{a+\epsilon_5-\epsilon_3}X^{b})^*-(a_6+1)(Y^{a+\epsilon_6-\epsilon_4}X^{b})^*\nonumber\\
		&\:-(a_8+1)(Y^{a+\epsilon_8-\epsilon_7}X^{b})^*+(a_{16}+1)(Y^{a+\epsilon_{16}-\epsilon_{15}}X^{b})^*\nonumber\\
		&\:+(b_{4,-1}+1)(Y^aX^{b+\epsilon_{4,-1}-\epsilon_{4,2}})^*-(b_{4,-2}+1)(Y^aX^{b+\epsilon_{4,-2}-\epsilon_{4,1}})^*\nonumber\\
		&\:+(b_{3,-1}+1)(Y^aX^{b+\epsilon_{3,-1}-\epsilon_{3,2}})^*-(b_{3,-2}+1)(Y^aX^{b+\epsilon_{3,-2}-\epsilon_{3,1}})^*\nonumber\\
		&\:+(b_{2,-1}+1)(Y^aX^{b+\epsilon_{2,-1}-\epsilon_{2,2}})^*-(b_{2,-2}+1)(Y^aX^{b+\epsilon_{2,-2}-\epsilon_{2,1}})^*\nonumber\\
		&\:+(b_{1,-1}+1)(Y^aX^{b+\epsilon_{1,-1}})^*,\label{al2}
	\end{align}
	\begin{align}
		E_{-\alpha_3}((Y^aX^b)^*)\nonumber=&\:(a_4+1)(Y^{a+\epsilon_4-\epsilon_3}X^{b})^*+(a_6+1)(Y^{a+\epsilon_6-\epsilon_5}X^{b})^*\nonumber\\
		&\:+(a_{12}+1)(Y^{a+\epsilon_{12}-\epsilon_{11}}X^{b})^*+(a_{14}+1)(Y^{a+\epsilon_{14}-\epsilon_{13}}X^{b})^*\nonumber\\
		&\:+(b_{4,3}+1)(Y^aX^{b+\epsilon_{4,3}-\epsilon_{4,4}})^*-(b_{4,-4}+1)(Y^aX^{b+\epsilon_{4,-4}-\epsilon_{4,-3}})^*\nonumber\\
		&\:+(b_{3,3}+1)(Y^aX^{b+\epsilon_{3,3}})^*,
	\end{align}
	\begin{align}
		E_{-\alpha_4}((Y^aX^b)^*)\nonumber=&\:-(a_3+1)(Y^{a+\epsilon_3-\epsilon_2}X^{b})^*-(a_9+1)(Y^{a+\epsilon_9-\epsilon_6}X^{b})^*\nonumber\\
		&\:-(a_{11}+1)(Y^{a+\epsilon_{11}-\epsilon_8}X^{b})^*+(a_{15}+1)(Y^{a+\epsilon_{15}-\epsilon_{14}}X^{b})^*\nonumber\\
		&\:+(b_{4,2}+1)(Y^aX^{b+\epsilon_{4,2}-\epsilon_{4,3}})^*-(b_{4,-3}+1)(Y^aX^{b+\epsilon_{4,-3}-\epsilon_{4,-2}})^*\nonumber\\
		&\:+(b_{3,2}+1)(Y^aX^{b+\epsilon_{3,2}-\epsilon_{3,3}})^*-(b_{3,-3}+1)(Y^aX^{b+\epsilon_{3,-3}-\epsilon_{3,-2}})^*\nonumber\\
		&\:+(b_{2,2}+1)(Y^aX^{b+\epsilon_{2,2}})^*,
	\end{align}
	\begin{align}
		E_{-\alpha_5}((Y^aX^b)^*)\nonumber=&\:-(a_2+1)(Y^{a+\epsilon_2-\epsilon_1}X^{b})^*+(a_{10}+1)(Y^{a+\epsilon_{10}-\epsilon_9}X^{b})^*\nonumber\\
		&\:+(a_{13}+1)(Y^{a+\epsilon_{13}-\epsilon_{11}}X^{b})^*+(a_{14}+1)(Y^{a+\epsilon_{14}-\epsilon_{12}}X^{b})^*\nonumber\\
		&\:+(b_{4,1}+1)(Y^aX^{b+\epsilon_{4,1}-\epsilon_{4,2}})^*-(b_{4,-2}+1)(Y^aX^{b+\epsilon_{4,-2}-\epsilon_{4,-1}})^*\nonumber\\
		&\:+(b_{3,1}+1)(Y^aX^{b+\epsilon_{3,1}-\epsilon_{3,2}})^*-(b_{3,-2}+1)(Y^aX^{b+\epsilon_{3,-2}-\epsilon_{3,-1}})^*\nonumber\\
		&\:+(b_{2,1}+1)(Y^aX^{b+\epsilon_{2,1}-\epsilon_{2,2}})^*-(b_{2,-2}+1)(Y^aX^{b+\epsilon_{2,-2}-\epsilon_{2,-1}})^*\nonumber\\
		&\:+(b_{1,1}+1)(Y^aX^{b+\epsilon_{1,1}})^*,
	\end{align}
	and
	\begin{eqnarray}
		E_{-\alpha_6}((Y^aX^b)^*)=(a_1+1)(Y^{a+\epsilon_1}X^b)^*. \label{al6}
	\end{eqnarray}
	According to (\ref{-al1}) and (\ref{al1}), we have
	\begin{align}
		&\:(-1)^{\operatorname{ht}(\lambda-\operatorname{wt}(a,b))+1-||a||}\varphi(E_{-\alpha_1}\eta^{(a,b)}w_{\lambda})\nonumber\\
		=&\:-(a+\epsilon_7)!b!(Y^{a+\epsilon_7-\epsilon_4}X^{b})^*-(a+\epsilon_8)!b!(Y^{a+\epsilon_8-\epsilon_6}X^b)^*-(a+\epsilon_{11})!b!(Y^{a+\epsilon_{11}-\epsilon_9}X^{b})^*\nonumber\\
		&\:-(a+\epsilon_{13})!b!(Y^{a+\epsilon_{13}-\epsilon_{10}}X^{b})^*-a!(b+\epsilon_{4,4})!(Y^{a}X^{b+\epsilon_{4,4}})^*\nonumber\\
		=&\:-a!b!E_{-\alpha_1}((Y^aX^b)^*)=-(-1)^{\operatorname{ht}(\lambda-\operatorname{wt}(a,b))-||a||}E_{-\alpha_1}(\varphi(\eta^{(a,b)}w_{\lambda})).
	\end{align}
	This yields
	\begin{align}
		\varphi(E_{-\alpha_1}\eta^{(a,b)}w_{\lambda})=E_{-\alpha_1}(\varphi(\eta^{(a,b)}w_{\lambda})).
	\end{align}
	The same verification is carried out for each $E_{-\alpha_i}$ ($i\in\overline{2,6}$) using the corresponding formulas (\ref{-al2})-(\ref{-al6}) and (\ref{al2})-(\ref{al6}). Therefore, $\varphi$ commutes with all negative simple root vectors.
	
	Finally, we examine the action of positive simple root vectors. Note that $1^*$ is a highest weight vector of $\mathscr{C}_{\lambda}^{\vee}$, we have
	\begin{align}
		\varphi(E_{\alpha_i}w_{\lambda})=0=E_{\alpha_i}(1^*)=E_{\alpha_i}(\varphi(w_{\lambda})).
	\end{align}
	Assume by induction on the sum of powers $||a||+||b||$ that $\varphi(E_{\alpha_i}\eta^{(a',b')}w_\lambda)=E_{\alpha_i}\varphi(\eta^{(a',b')}w_\lambda)$ for all $(a',b')\in\mathbb{N}^{36}$ with $||a'||+||b'||=r\geqslant0$. Take $(a,b)\in\mathbb{N}^{36}$ with $||a||+||b||=r+1$. By the PBW basis property we can write
	\begin{eqnarray}
		\eta^{(a,b)}w_{\lambda}=X\eta^{(a',b')}w_{\lambda}
	\end{eqnarray}
	for some negative root vector $X$ in (\ref{order}) and some $(a',b')\in\mathbb{N}^{36}$ with $||a'||+||b'||=r$. There are two cases. 
	
	If $E_{\alpha_i}$ commutes with $X$, then
	\begin{align}
		\varphi(E_{\alpha_i}\eta^{(a,b)}w_{\lambda})=&\:\varphi(E_{\alpha_i}X\eta^{(a',b')}w_{\lambda})=\varphi(XE_{\alpha_i}\eta^{(a',b')}w_{\lambda})\nonumber\\
		=&\:X\varphi(E_{\alpha_i}\eta^{(a',b')}w_{\lambda})=XE_{\alpha_i}\varphi(\eta^{(a',b')}w_{\lambda})\nonumber\\
		=&\:E_{\alpha_i}X(\varphi(\eta^{(a',b')}w_{\lambda}))=E_{\alpha_i}(\varphi(\eta^{(a,b)}w_{\lambda})).
	\end{align}
	If $E_{\alpha_i}$ does not commute with $X$, set
	\begin{align}
		u=[E_{\alpha_i},X]
	\end{align}
	with $u$ being a negative root vector or $u\in \mathscr{H}$. Therefore, by induction assumption we have
	\begin{align}
		\varphi(E_{\alpha_i}\eta^{(a,b)}w_{\lambda})=&\:\varphi((XE_{\alpha_i}+u)\eta^{(a',b')}w_{\lambda})\nonumber\\
		=&\:XE_{\alpha_i}\varphi(\eta^{(a',b')}w_{\lambda})+u\varphi(\eta^{(a',b')}w_{\lambda})\nonumber\\
		=&\:E_{\alpha_i}X\varphi(\eta^{(a',b')}w_{\lambda})=E_{\alpha_i}(\varphi(\eta^{(a,b)}w_{\lambda})).
	\end{align}
	In both cases we obtain $\varphi(E_{\alpha_i}\eta^{(a,b)}w_\lambda)=E_{\alpha_i}\varphi(\eta^{(a,b)}w_\lambda)$. By induction, $\varphi$ commutes with the action of every positive simple root vectors as well.
	Consequently, $\varphi$ is a $\mathscr{G}^{E_6}$-module isomorphism.	
\end{proof}
\begin{remark}
	The sign $(-1)^{\operatorname{ht}(\lambda-\operatorname{wt}(a,b))-||a||}$ in the formula for $\varphi$ is due to the basis we chose for the negative root vectors used in our PBW basis. For example, the chosen simple negative root vectors are
	\begin{eqnarray}
		&A_{4,5}=-E_{-\alpha_1},\;\;A_{-1,2}=-E_{-\alpha_2},\;\;A_{3,4}=-E_{-\alpha_3},\nonumber\\
		&A_{2,3}=-E_{-\alpha_4},\;\;A_{1,2}=-E_{-\alpha_5},\;\;\eta_1=E_{-\alpha_6}.
	\end{eqnarray}
\end{remark}

Since the duality functor $(~\bullet~)^{\vee}$ on category $\mathcal{O}$ is self‑inverse, Proposition \ref{prop3.4} and Theorem \ref{thiso} immediately imply the following.
\begin{corollary}\label{Coro3.6}
	As a module over $\mathscr{G}^{E_6}$, $\mathscr{C}_{\lambda}$ is isomorphic to the restricted dual Verma module $M(\lambda)^{\vee}$, hence the cyclic $\mathscr{G}^{E_6}$-submodule $U(\mathscr{G}^{E_6})(1)$ is the unique irreducible submodule of $\mathscr{C}_{\lambda}$.
\end{corollary}

\begin{lemma}\label{lem3.8}
	If a linear operator $\phi:\mathscr{C}_{\lambda}\longrightarrow\mathscr{C}_{\lambda}$ commutes with $\pi_{\lambda}(E_{\alpha_i})$ for each $i\in\overline{1,6}$, and the constant term of $\phi(f)$ is zero for any $f\in\mathscr{C}$, then $\phi=0$.
\end{lemma}
\begin{proof}
	The action of the Cartan subalgebra $\mathscr{H}$ determines a positive grading on $\mathscr{C}_{\lambda}$ by root height. Explicitly, let the weight space $(\mathscr{C}_{\lambda})_{\mu}$ with weight $\mu\in\mathscr{H}^{*}$ is nonzero if and only if $\lambda-\mu\in\sum_{r=1}^{6}\mathbb{N}\alpha_r$. For $n\in\mathbb{N}$ let
	\begin{align}
		(\mathscr{C}_{\lambda})_n=\bigoplus_{\operatorname{ht}(\lambda-\mu)=n}(\mathscr{C}_{\lambda})_\mu,
	\end{align}
	then
	\begin{align}
		\mathscr{C}_{\lambda}=\bigoplus_{n\in\mathbb{N}}(\mathscr{C}_{\lambda})_n\text{ and }(\mathscr{C}_{\lambda})_0=\mathbb{C}.
	\end{align}
	The action of every simple positive root vector $\pi_{\lambda}(E_{\alpha_i})$ decreases the grading degree by one. 
	If $\phi\neq0$, then there exist $n_0\in\mathbb{N}$ such that
	\begin{eqnarray}
		\operatorname{Im}\phi\cap (\bigoplus_{n=0}^{n_0}(\mathscr{C}_{\lambda})_n)\neq0\text{ and }\operatorname{Im}\phi\cap (\bigoplus_{n=0}^{n_0-1}(\mathscr{C}_{\lambda})_n)=0.
	\end{eqnarray}
	Take nonzero $f\in\operatorname{Im}\phi\cap (\bigoplus_{n=0}^{n_0}(\mathscr{C}_{\lambda})_n)$.
	All simple positive root vectors preserve $\operatorname{Im}\phi$ since they commute with $\phi$. If $\pi_{\lambda}(E_{\alpha_i})(f)\neq 0$ for some $i$, then 
	\begin{align}
		\pi_{\lambda}(E_{\alpha_i})(f)\in\operatorname{Im}\phi\cap (\bigoplus_{n=0}^{n_0-1}(\mathscr{C}_{\lambda})_n),
	\end{align}
	a contradiction. Thus $\pi_{\lambda}(E_{\alpha_i})f=0$ for all $i$, and $f$ is annihilated by all the positive root vectors since all positive root vectors can be generated by the bracket of simple positive root vectors. By Proposition \ref{prop3.4}, $f$ must be a constant. Therefore, $f$ must be zero by the second assumption of this lemma, a contradiction. Hence $\phi=0$.
\end{proof}

\subsection{${D_5}$-singular vectors}

In this subsection, we restrict the $\mathscr{G}^{E_6}$-modules $\mathscr{C}_{\lambda}$ to the subalgebra $\mathscr{G}^{D_5}\cong\mathfrak{o}(10)$, and determine all $\mathscr{G}^{D_5}$-singular vectors in $\mathscr{C}_{\lambda}$. 
Recall the polynomial algebra $\mathscr{C}_{\lambda}=\mathbb{C}[y_k,x_{i,j}\mid 1\leqslant k\leqslant 16,\; 1\leqslant |j|\leqslant i\leqslant 4]$ equipped with the action of $\mathscr{G}^{E_6}$ as described in Theorem \ref{thm3.1} (as well as in Theorem \ref{thm3.2}). The $\mathscr{G}^{D_5}$-singular vectors in $\mathscr{C}_{\lambda}$ are those weight vector in $\mathscr{C}_{\lambda}$ simultaneously annihilated by all positive root vectors $E_{\alpha}$ with $\alpha\in\Phi_{D_5}^+$.

For $u\in\mathscr{G}^{E_6}_0$, 
denote by $c(u)$ the $16\times 16$ matrix with $(i,j)$-entry $c_{i,j}(u)$. It is the matrix of $\sigma(u)$ on the $16$-dimensional space $\sum_{i=1}^{16}\mathbb{C}y_i$ with respect to the basis $\{y_1,\cdots,y_{16}\}$ by definition (\ref{248}), so $c(\cdot)$ preserves the Lie bracket on $\mathscr{G}^{E_6}_0$.
By Appendix \ref{app}, for each $r\in\overline{1,5}$ the matrix $c(A_{r,r})$ is diagonal with eigenvalues $\{\frac{1}{2},-\frac{1}{2}\}$. Thus we have 
\begin{eqnarray}
	[c(A_{r,r}),c(A_{r,s})]=c([A_{r,r},A_{r,s}])=c(A_{r,s})
\end{eqnarray}
for $1\leqslant |s|<r$, from which we see that $c(A_{r,s})$ increases the eigenvalue of $c(A_{r,r})$ by $1$. Consequently,
\begin{eqnarray}
	c(A_{r,s})c(A_{r,t})=0 \label{Gr}
\end{eqnarray}
for $1\leqslant |s|,|t|<r$. Set
\begin{eqnarray}
	&G_r=\sum\limits_{|s|=1}^{r-1}x_{r-1,s}c(A_{r,s})\quad(r\in\overline{2,5}),\label{Gr1}\\
	&G=(I_{16}+G_2)(I_{16}+G_3)(I_{16}+G_4)(I_{16}+G_5).\label{Gr2}
\end{eqnarray}
By definition the matrices $G_2,\cdots,G_5,G\in M_{16\times 16}(\mathscr{X})$. Eq. (\ref{Gr}) implies $G_r^2=0$ for each $r\in\overline{2,5}$, hence 
\begin{eqnarray}
	G^{-1}=(I_{16}-G_5)(I_{16}-G_4)(I_{16}-G_3)(I_{16}-G_2)\in M_{16\times 16}(\mathscr{X}). \label{G-1}
\end{eqnarray}
Denote
\begin{eqnarray}
	(\varkappa_1,\cdots,\varkappa_{16})=(y_1,\dots,y_{16})G^{-1}. 
\end{eqnarray}
The sets
\begin{align}
	\{\varkappa_1,\cdots,\varkappa_{16}\}\cup\{x_{i,j}\mid 1\leqslant|j|\leqslant i\leqslant 4\}\quad\text{and}\quad\{y_1,\cdots,y_{16}\}\cup\{x_{i,j}\mid 1\leqslant|j|\leqslant i\leqslant 4\}
\end{align} 
are polynomial equivalent, That is, any polynomial in the first set can be written as a polynomial in the second set, and vice versa. This also shows that $\varkappa_1,\cdots,\varkappa_{16}$ are algebra independent.
\begin{theorem}
	As a $\mathscr{G}^{D_5}$-module, any singular vector in $\mathscr{C}_{\lambda}$ is a polynomial in $\{\varkappa_1,\cdots,\varkappa_{16}\}$. Moreover, the monomial $\prod_{i=1}^{16}\varkappa_i^{r_i}$ is a weight vector of weight \begin{align}
		&\:(\mu_5-\mu_4+r_4+r_6-r_7-r_8+r_9+r_{10}-r_{11}-r_{13})\omega_1\nonumber\\+&\:(\mu_2+\mu_1+r_3+r_4-r_5-r_6+r_7-r_8+r_{15}-r_{16})\omega_2\nonumber\\
		+&\:(\mu_4-\mu_3+r_3-r_4+r_5-r_6+r_{11}-r_{12}+r_{13}-r_{14})\omega_3\nonumber\\+&\:(\mu_3-\mu_2+r_2-r_3+r_6+r_8-r_9-r_{11}+r_{14}-r_{15})\omega_4\nonumber\\+&\:(\mu_2-\mu_1+r_1-r_2+r_9-r_{10}+r_{11}+r_{12}-r_{13}-r_{14})\omega_5. \label{3.105}
	\end{align}
\end{theorem}
\begin{proof}
	We first show that any singular vector in $\mathscr{C}_{\lambda}$ must be a polynomial in $\{\varkappa_1,\cdots,\varkappa_{16}\}$.
	Take any positive root vectors $A_{i,j}\in\mathfrak{o}(10)$ with $1\leqslant |j| <i\leqslant 5$. By Theorem \ref{thm3.1} we have $\pi_{\lambda}(A_{i,j})=\varrho_{\mu}(A_{i,j})+\sigma(A_{i,j})\in\mathscr{W}_{\mathscr{C}}$ with $\varrho_{\mu}(A_{i,j})\in\mathscr{W}_{\mathscr{X}}$ given in (\ref{2.24}) and $\sigma(A_{i,j})\in\mathscr{W}_{\mathscr{Y}}$ given in Appendix \ref{app}. In the following, the derivations $\pi_{\lambda}(A_{i,j}),\varrho_{\mu}(A_{i,j})$ and $\sigma(A_{i,j})$ acting on matrices is by acting on each entry simultaneously.
	By (\ref{2.24}), we have
	\begin{align}
		\varrho_{\mu}(A_{i,j})(G_r)=\begin{cases}
			0, &\text{if }r<i,\\
			c(A_{i,j}), &\text{if }r=i,\\
			x_{r-1,i}c(A_{r,j})-x_{r-1,-j}c(A_{r,-i}), &\text{if }r>i.
		\end{cases} \label{4.55}
	\end{align}
	It follows from the Lie bracket of positive root vectors $[A_{r,s},A_{i,j}]=\delta_{s,i}A_{r,j}-\delta_{s,-j}A_{r,-i}$ for $r>i$ that the last row in (\ref{4.55}) precisely equals to $G_rc(A_{i,j})-c(A_{i,j})G_r$. This yields
	\begin{align}
		\varrho_{\mu}(A_{i,j})(I_{16}-G_r)=\begin{cases}
			0, &\text{if }r<i,\\
			-c(A_{i,j}), &\text{if }r=i,\\
			(I_{16}-G_r)c(A_{i,j})-c(A_{i,j})(I_{16}-G_r), &\text{if }r>i.
		\end{cases}\label{4.56}
	\end{align}
	Now we can take part the expression of $G^{-1}$ in (\ref{G-1})  into product of  factors with respect to the index $r>i$, $r=i$ and $r<i$ as
	\begin{eqnarray}
		G^{-1}=(I_{16}-G_5)\cdots(I_{16}-G_{i+1})\times (I_{16}-G_{i})\times (I_{16}-G_{i-1})\cdots (I_{16}-G_{2}). \label{3.108}
	\end{eqnarray}	
	Leibniz rule and (\ref{4.56}) therefore give
	\begin{eqnarray}
		\pi_{\lambda}(A_{i,j})(G^{-1})=\varrho_{\mu}(A_{i,j})(G^{-1})=-c(A_{i,j})G^{-1}. \label{4.57}
	\end{eqnarray}
	On the other hand, Eq. (\ref{248}) gives
	\begin{eqnarray}
		\pi_{\lambda}(A_{i,j})((y_1,\cdots,y_{16}))=\sigma(A_{i,j})((y_1,\cdots,y_{16}))=(y_1,\cdots,y_{16})c(A_{i,j}).\label{3.110}
	\end{eqnarray}
	Now from (\ref{4.57}) and (\ref{3.110}) we have
	\begin{align}
		&\pi_{\lambda}(A_{i,j})((\varkappa_1,\cdots,\varkappa_{16}))\nonumber\\
		=&\pi_{\lambda}(A_{i,j})((y_1,\cdots,y_{16})G^{-1})\nonumber\\
		=&\pi_{\lambda}(A_{i,j})((y_1,\cdots,y_{16}))G^{-1}+(y_1,\cdots,y_{16})\pi_{\lambda}(A_{i,j})(G^{-1})\nonumber\\
		=&(y_1,\cdots,y_{16})c(A_{i,j})G^{-1}-(y_1,\cdots,y_{16})c(A_{i,j})G^{-1}\nonumber\\
		=&(0,\cdots,0).
	\end{align}
	Equivalently, the polynomials $\varkappa_1,\cdots,\varkappa_{16}$ are killed by $A_{i,j}$.
	
	Suppose that $f\in\mathscr{C}_\lambda$ is a $D_5$-singular vector. Write $f$ as a polynomial in $\{\varkappa_1,\cdots,\varkappa_{16}\}\cup\{x_{r,s}\mid 1\leqslant|s|\leqslant r\leqslant 4\}$. Acting $A_{5,j}$ with $1\leqslant |j|\leqslant 4$ on $f$ yields
	\begin{eqnarray}
		\pi_{\lambda}(A_{5,j})(f)=\sum_{k=1}^{16}\pi_{\lambda}(A_{5,j})(\varkappa_k)\frac{\partial f}{\partial\varkappa_k}+\sum_{1\leqslant|s|\leqslant r\leqslant 4}\pi_{\lambda}(A_{5,j})(x_{r,s})\frac{\partial f}{\partial x_{r,s}}=\frac{\partial f}{\partial x_{4,j}}=0
	\end{eqnarray}
	by (\ref{2.24}), hence $f$ is independent of variables $x_{4,j}$ ($1\leqslant |j|\leqslant 4$). Similarly, acting $A_{4,j}$ ($1\leqslant |j|\leqslant 3$), $A_{3,j}$ ($1\leqslant |j|\leqslant 2$), $A_{2,\pm1}$ successively on $f$, we can derive that $f$ is independent of $x_{i,j}$ ($1\leqslant |j|\leqslant i\leqslant 4$). Therefore, $f$ must be a polynomial in $\{\varkappa_1,\cdots,\varkappa_{16}\}$.
	
	Next we prove the second assertion. Take $A_{i,i}\in\mathfrak{o}(10)$ with $i\in\overline{1,5}$. by (\ref{r2}) we have $\pi_{\lambda}(A_{i,i})=\varrho_{\mu}(A_{i,i})+\sigma(A_{i,i})$ with $\varrho_{\mu}(A_{i,i})-\mu_i\in\mathscr{W}_{\mathscr{X}}$ and $\sigma(A_{i,i})\in\mathscr{W}_{\mathscr{Y}}$. By (\ref{2.25}), we have
	\begin{align}
		\varrho_{\mu}(A_{i,i})(G_r)-\mu_iG_r=\begin{cases}
			0, &\text{if }r<i,\\
			-G_r, &\text{if }r=i,\\
			x_{r-1,i}c(A_{r,i})-x_{r-1,-i}c(A_{r,-i}), &\text{if }r>i.
		\end{cases} \label{3.112}
	\end{align}
	Since $[A_{r,s},A_{i,i}]=(\delta_{s,i}-\delta_{s,-i})A_{r,s}$ for $r>i$ and $1\leqslant |s| <r$, the last row in (\ref{3.112}) precisely equals to $G_rc(A_{i,i})-c(A_{i,i})G_r$. This yields
	\begin{align}
		\varrho_{\mu}(A_{i,i})(I_{16}-G_r)-\mu_i(I_{16}-G_r)=\begin{cases}
			0, &\text{if }r<i,\\
			G_r, &\text{if }r=i,\\
			(I_{16}-G_r)c(A_{i,i})-c(A_{i,i})(I_{16}-G_r), &\text{if }r>i.
		\end{cases}\label{3.113}
	\end{align}
	Taking $G^{-1}$ part into product of factors with respect to the index $r>i$, $r=i$ and $r<i$ as in (\ref{3.108}) again, then applying Leibniz rule and (\ref{3.113}) , we obtain
	\begin{equation}
		(\pi_{\lambda}(A_{i,i})-\mu_i)(G^{-1})=(\varrho_{\mu}(A_{i,i})-\mu_i)(G^{-1})=-c(A_{i,i})G^{-1}+G^{-1}c(A_{i,i}). \label{3,114}
	\end{equation}
	Eq. (\ref{248}) gives
	\begin{eqnarray}
		(\pi_{\lambda}(A_{i,i})-\mu_i)((y_1,\cdots,y_{16}))=\sigma(A_{i,i})(y_1,\cdots,y_{16})=(y_1,\cdots,y_{16})c(A_{i,i}). \label{3,115}
	\end{eqnarray}
	By (\ref{3,114}) and (\ref{3,115}) we have
	\begin{align}
		&(\pi_{\lambda}(A_{i,i})-\mu_i)((\varkappa_1,\cdots,\varkappa_{16}))\nonumber\\
		=&(\pi_{\lambda}(A_{i,i})-\mu_i)((y_1,\cdots,y_{16}))G^{-1}+(y_1,\cdots,y_{16})(\pi_{\lambda}(A_{i,i})-\mu_i)(G^{-1})\nonumber\\
		=&(y_1,\cdots,y_{16})c(A_{i,i})G^{-1}-(y_1,\cdots,y_{16})c(A_{i,i})G^{-1}+(y_1,\cdots,y_{16})G^{-1}c(A_{i,i})\nonumber\\
		=&(\varkappa_1,\cdots,\varkappa_{16})c(A_{i,i}). \label{3,116}
	\end{align}
	Comparing (\ref{3,115}) and (\ref{3,116}) we see that $\varkappa_k$ is a weight vector, whose weight is the same as that of $y_k$ ($(k\in\overline{1,16}$).
	Hence the monomial $\prod_{i=1}^{16}\varkappa_i^{r_i}$ is a weight vector, whose weight is the same as that of $\prod_{i=1}^{16}y_i^{r_i}$. Consequently, by (\ref{3.97})-(\ref{3.101}), the weight of $\prod_{i=1}^{16}\varkappa_i^{r_i}$ is precisely (\ref{3.105}). This completes the proof.	
\end{proof}
\begin{remark}
	Since positive root vectors are generated by the simple ones,
	one can also obtain all the $\mathscr{G}^{D_5}$-singular vectors by using the method of characteristics to solve the system of homogeneous linear first-order partial differential equations
	\begin{eqnarray}
		\pi_{\lambda}(E_{\alpha_i})(f)=0\quad (i\in\overline{1,5})
	\end{eqnarray}
	determined by (\ref{3.91})-(\ref{3.95}).
	In fact, such a method can be adapted to obtain the singular vectors restricting to any Levi subalgebra of $\mathscr{G}^{E_{6}}$ with respect to a subset of simple positive roots $\Sigma$, for which one need to solve the system $\pi_{\lambda}(E_{\alpha})(f)=0$ for $\alpha\in\Sigma$.
\end{remark}

\section{Finite-dimensional irreducible modules}\label{section4}

In this section, we turn to investigate the finite-dimensional irreducible $\mathscr{G}^{E_6}$-submodule. Recall that the isomorphism classes of finite-dimensional irreducible $\mathscr{G}^{E_6}$-modules are in bijection with the dominant integral weights, via their highest weights. 

According to Proposition \ref{prop3.4} and Corollary \ref{Coro3.6}, the $\mathscr{G}^{E_6}$-submodule $V(\lambda)=U(\mathscr{G}^{E_6})(1)$ of $\mathscr{C}_{\lambda}$ is an irreducible highest weight module with highest weight $\lambda$. Throughout this section, we take $\lambda$ to be dominant integral, equivalently, $\lambda_1,\cdots,\lambda_6\in\mathbb{N}$. We construct explicit polynomial bases for the six fundamental modules $V(\omega_i)$, $i\in\overline{1,6}$, then a finite spanning set for $V(\lambda)$ with arbitrary $\lambda$ is obtained by the compatibility of polynomial multiplication with addition of highest weights.

\subsection{Spanning set}

\begin{lemma}\label{lem4.1}
	If a finite-dimensional subspace $V\subset\mathscr{C}_{\lambda}$ contains $1$ and it is stable under the actions of all simple positive and negative root vectors, then $V=V(\lambda)$.
\end{lemma}
\begin{proof}
	The simple root vectors generate $\mathscr{G}^{E_6}$ as a Lie algebra, so $V$ is a finite-dimensional $\mathscr{G}^{E_6}$-module. By Weyl's theorem of complete reducibility, $V$ can be decomposed into a direct sum of irreducible submodules, as each irreducible component has a nonzero singular vector. All such vectors are multiple of $1$ according to Proposition \ref{prop3.4}, hence there is only one summand, then $V$ is irreducible and is isomorphic to $V(\lambda)$.
\end{proof}

Given $n$ subsets $S_1,S_2,\cdots,S_n\subset\mathscr{C}$, we define their $product$ $\prod_{i=1}^nS_i\subset\mathscr{C}$ to be the set of all product $\prod_{i=1}^nf_i$ of polynomials with $f_i\in S_i$. For a single set $S$, we write $S^n$ for the $n$-fold product.
\begin{proposition}\label{prop3.6}
	For dominant integral weight $\lambda$ and $\lambda'$, suppose $S,T\subset\mathscr{C}$ satisfy\begin{eqnarray}
		&V(\lambda)=\mathrm{Span}_{\mathbb{C}}\:S\quad\text{and} \quad V(\lambda')=\mathrm{Span}_{\mathbb{C}}\:T.
	\end{eqnarray} 
	Then
	\begin{eqnarray}
		&V(\lambda+\lambda')=\mathrm{Span}_{\mathbb{C}}\:ST.
	\end{eqnarray}
\end{proposition}
\begin{proof}
	Write $\lambda=\sum_{i=1}^6\lambda_i\omega_i$ and $\lambda'=\sum_{i=1}^6\lambda_i'\omega_i$. Take any $u\in\mathscr{G}^{E_6},p\in S, q\in T$. According to Lemma \ref{lem3.3}, we have
	\begin{eqnarray}
		u(p)=u^D(p)+\mathop{\textstyle\sum}\limits_{i=1}^6\lambda_iu^P_ip\in V(\mathop{\textstyle\sum}\limits_{i=1}^6\lambda_i\omega_i)=\mathrm{Span}_{\mathbb{C}}S.\label{3.126}
	\end{eqnarray}
	Analogously, we have
	\begin{eqnarray}
		u(q)=u^D(q)+\mathop{\textstyle\sum}\limits_{i=1}^6\lambda'_iu^P_iq\in V(\mathop{\textstyle\sum}\limits_{i=1}^6\lambda'_i\omega_i)=\mathrm{Span}_{\mathbb{C}}T.\label{3.127}
	\end{eqnarray}
	Now we compute
	\begin{align}
		u(pq)=&\:u^{D}(pq)+\mathop{\textstyle\sum}\limits_{i=1}^6(\lambda_i+\lambda'_i)u^P_ipq\nonumber\\
		=&\:(u^D(p)+\mathop{\textstyle\sum}\limits_{i=1}^6\lambda_iu^P_ip)q+p(u^D(q)+\mathop{\textstyle\sum}\limits_{i=1}^6\lambda'_iu^P_iq).
	\end{align}
	By (\ref{3.126}) and (\ref{3.127}), the first term is linear combination of elements of $S$ multiplied by $q$, and the second term is $p$ times a linear combination of element of $T$. Thus $u(pq)\in\mathrm{Span}_{\mathbb{C}}ST$. Since $u$ is arbitrary, we conclude that $\mathrm{Span}_{\mathbb{C}}ST$ is a $\mathscr{G}^{E_6}$-submodule of $\mathscr{C}_{\lambda+\lambda'}$.
	
	Both $\mathrm{Span}_{\mathbb{C}}S$ and $\mathrm{Span}_{\mathbb{C}}T$ are finite-dimensional by assumption, hence $\mathrm{Span}_{\mathbb{C}}ST$ is also finite-dimensional. By Corollary \ref{Coro3.6} we get
	\begin{eqnarray}
		\mathrm{Span}_{\mathbb{C}}ST=V(\lambda+\lambda').
	\end{eqnarray}
	This completes the proof. 
\end{proof}
\begin{theorem}\label{thmspan}
	Let $\mathcal{B}_i$ be a basis of $V(\omega_i)$ for $i\in\overline{1,6}$. Then for any dominant integral weight $\lambda=\sum_{i=1}^6\lambda_i\omega_i$, the finite-dimensional irreducible module $V(\lambda)$ is spanned by the set $\prod_{i=1}^6\mathcal{B}_i^{\lambda_i}$.
\end{theorem}
\begin{proof}
	This follows immediately from Proposition \ref{prop3.6} by induction on the sum $\sum_{i=1}^6{\lambda_i}$.
\end{proof}

\begin{remark}
	The dimensions of the fundamental modules $V(\omega_i)$ ($i\in\overline{1,6}$) are relatively small:
	\begin{eqnarray}
		&\dim V(\omega_1)=\dim V(\omega_6)=27,\;\;\dim V(\omega_2)=78,\\
		&\dim V(\omega_3)=\dim V(\omega_5)=351,\;\;\dim V(\omega_4)=2925.
	\end{eqnarray}
	Using the explicit representation formulas from Theorem \ref{thm3.2} together with a computer algebra program, one can compute a basis $\mathcal{B}_i$ for each $V(\omega_i)$. Theorem \ref{thmspan} then provides a finite spanning set for $V(\lambda)$ for any dominant integral weight $\lambda$. The extraction of a genuine basis from this spanning set is a problem that the author hopes to address in the future work.
\end{remark}
In the remainder of this section, we exhibit an explicit basis for each $V(\omega_i)$ that has particularly simple forms.

\subsection{Polynomial bases for two $27$-dimensional modules}
Set
\begin{eqnarray}
	F_0=\sum_{r=1}^4x_{4,r}x_{4,-r}\quad\text{and}\quad
	F_1=\zeta_6-\zeta_1F_0+\sum_{r=1}^{4}(\zeta_{6-r}x_{4,-r}+\zeta_{11-r}x_{4,r}). \label{4.10}
\end{eqnarray}
\begin{proposition}\label{prop4.5}
	\begin{enumerate}[(i)]
		\item The set
		\begin{align}
			\mathcal{B}_1=\{1,x_{4,r},x_{4,-r},F_0,F_1,\partial_{y_k}(F_1)\mid r\in\overline{1,4},\;k\in\overline{1,16}\} \label{4,10}
		\end{align}
		forms a basis of $V(\omega_1)$.
		\item The set
		\begin{eqnarray}
			\mathcal{B}_6=\{1,y_k,\zeta_{\ell}\mid k\in\overline{1,16},\;\ell\in\overline{1,10}\} \label{b60}
		\end{eqnarray}
		forms a basis of $V(\omega_6)$.
	\end{enumerate}	
\end{proposition}
\begin{proof}
	Take $\lambda=\omega_1$ firstly. By (\ref{3.122})-(\ref{3.124}), the corresponding parameters are
	\begin{eqnarray}
		\mu_1=\mu_2=\mu_3=\mu_4=0,\;\mu_5=\kappa=1.
	\end{eqnarray}
	Applying the representation formula (\ref{r4}) of $\eta_i$, we have
	\begin{eqnarray}
		\eta_7(\eta_{16}(1))=F_1\in V(\omega_1).
	\end{eqnarray}
	Applying the representation formula (\ref{r3}) of $\xi_i$, we see that the polynomials
	\begin{eqnarray}
		\partial_{y_k}(F_1)\in V(\omega_1), \quad k\in\overline{1,16}, \label{b12}
	\end{eqnarray}
	and have $y$-degree $1$. Moreover,
	\begin{eqnarray}
		&\partial_{y_1}\partial_{y_{10}}(F_1)=-F_0,\;\;\partial_{y_1}\partial_{y_{16}}(F_1)=x_{4,1},\;\;\partial_{y_2}\partial_{y_{16}}(F_1)=x_{4,2},\label{b13}\\
		&\partial_{y_3}\partial_{y_{16}}(F_1)=x_{4,3},\;\;\partial_{y_4}\partial_{y_{16}}(F_1)=-x_{4,4},\;\;\partial_{y_2}\partial_{y_{15}}(F_1)=-x_{4,-1},\label{b14}\\
		&\partial_{y_1}\partial_{y_{15}}(F_1)=x_{4,-2},\;\;\partial_{y_1}\partial_{y_{14}}(F_1)=x_{4,-3},\;\;\partial_{y_1}\partial_{y_{13}}(F_1)=x_{4,-4} \in V(\omega_1 ),\label{b15}
	\end{eqnarray}
	 these $9$ polynomials all have $y$-degree $0$ and are apparently linear independent. To see the linear independence of the polynomials in (\ref{b12}), suppose 
	 \begin{eqnarray}
	 	\sum_{k=1}^{16}c_k\partial_{y_k}(F_1)=0\quad\text{with }c_k\in\mathbb{C}.
	 \end{eqnarray}
	 The term with degree $2$ in $x_{4,r},x_{4,-r}$ ($r\in\overline{1,4}$) in $F_1$ is $-\zeta_1F_0$, hence by (\ref{zeta1})
	 \begin{eqnarray}
	 	\sum_{k=1}^{16}c_k\partial_{y_k}(\zeta_1)=c_1y_{10}+c_{10}y_1+c_2y_9+c_9y_{2}-c_3y_6-c_6y_3+c_4y_5+c_5y_4=0. \label{4.19}
	 \end{eqnarray}
	 The term with degree $0$ in $x_{4,r},x_{4,-r}$ ($r\in\overline{1,4}$) in $F_1$ is $\zeta_6$, hence by (\ref{zeta6})
	 \begin{align}
	 	\sum_{k=1}^{16}c_k\partial_{y_k}(\zeta_6)=c_7y_{16}+c_{16}y_7+c_8y_{15}+c_{15}y_8+c_{11}y_{14}+c_{14}y_{11}-c_{12}y_{13}-c_{13}y_{12}=0. \label{4.20}
	 \end{align}
	 Eqs. (\ref{4.19})-(\ref{4.20}) imply that every $c_k$ is zero. Eqs. (\ref{b12})-(\ref{b15}) then give $27$ linear independent elements in $V(\omega_1)$, proving (i).
	 
	 Now take $\lambda=\omega_6$. By (\ref{3.122})-(\ref{3.124}), the corresponding parameters are
	 \begin{eqnarray}
	 	\mu_1=\mu_2=\mu_3=\mu_4=\mu_5=0,\;\kappa=2.
	 \end{eqnarray}
	 By (\ref{r4}), we have
	 \begin{eqnarray}
	 	\eta_k(1)=-y_k,\;\;\eta_i\eta_j(1)=\sum_{\ell=1}^{10}\zeta_\ell\partial_{y_i}\partial_{y_j}(\zeta_{\ell}).
	 \end{eqnarray}
	 Choosing $(i,j)$ successively as $(1,10)$, $(1,13)$, $(1,14)$, $(1,15)$, $(2,15)$, $(7,16)$, $(4,16)$, $(3,16)$, $(2,16)$, $(1,16)$
	 produces $\zeta_1,\zeta_2,\zeta_3,\zeta_4,-\zeta_5,\zeta_6,-\zeta_7,\zeta_8,\zeta_9,\zeta_{10}$. Thus all $27$ polynomials in $\mathcal{B}_6$ lie in $V(\omega_6)$, and it is easy to see that they are linear independent. Since $\dim V(\omega_6)$=27, $\mathcal{B}_6$ forms a basis.
\end{proof}

\subsection{Polynomial bases for the $78$-dimensional module}

Next we consider $\lambda=\omega_2$. In this case, we have $\mu_1=\cdots=\mu_5=\frac{1}{2}$ and $\kappa=\frac{3}{2}$ by (\ref{3.122})-(\ref{3.124}). Retain the matrices $L,K,C\in M_{5\times 5}(\mathscr{X})$ defined in (\ref{2.8})-(\ref{2.9}). Let $T=L^{-1}C=L^{-1}K(L^t)^{-1}$ and let $t_{i,j}$ be the $(i,j)$-entry of $T$, then $T$ is skew-symmetric and each $t_{i,j}$ belongs to $\mathscr{X}$.

By (\ref{rho}), for any $u=\begin{pmatrix}
	P&Q\\R&-P^t
\end{pmatrix}\in\mathfrak{o}(10)$, the operator $d_u:=\varrho_{\mu}(u)-\frac{1}{2}\operatorname{tr}(P+TR)$ is a derivation on $\mathscr{X}$. Moreover, one can derive from (\ref{L}) and (\ref{C}) that 
\begin{eqnarray}
	d_u(L)=LP+CR-UL\quad\text{and}\quad d_u(C)=LQ-CP^t-UC, \label{4.23}
\end{eqnarray}
here $d_u$ acts on each entry of the matrices simultaneously, and $U$ is the upper triangular part of $(LP+CR)L^{-1}$ including the diagonal.
Eq. (\ref{4.23}) yields
\begin{eqnarray}
	d_uT=-L^{-1}d_u(L)L^{-1}C+L^{-1}d_u(C)=Q-PT-TP^t-TRT. \label{4.27}
\end{eqnarray}
For $i,j,k,l\in\overline{1,5}$, put
\begin{eqnarray}
	t_{ijkl}=t_{ij}t_{kl}-t_{ik}t_{jl}+t_{il}t_{jk}.
\end{eqnarray}
When $i,j,k,l$ are pairwise distinct, $t_{ijkl}$ is the Pfaffian of the $4\times 4$ principal submatrix of $T$ whose row and columns are both indexed by $i,j,k,l$ in that order.
From (\ref{4.27}) we see that for $r,s\in\overline{1,5}$,
\begin{eqnarray}
	&\varrho_{\mu}(A_{r,s})(1)=\frac{\delta_{r,s}}{2},\;\;\varrho_{\mu}(A_{r,s})(t_{ij})=-\delta_{i,r}t_{sj}-\delta_{j,r}t_{is}+\frac{1}{2}\delta_{r,s}t_{ij},\label{4.28}\\
	&\varrho_{\mu}(A_{r,s})(t_{ijkl})=-\delta_{i,r}t_{sjkl}-\delta_{j,r}t_{iskl}-\delta_{k,r}t_{ijsl}-\delta_{l,r}t_{ijks}+\frac{1}{2}\delta_{r,s}t_{ijkl}, \\
	&\varrho_{\mu}(A_{r,-s})(1)=0,\;\varrho_{\mu}(A_{r,-s})(t_{ij})=\delta_{i,r}\delta_{j,s}-\delta_{i,s}\delta_{j,r},\;\varrho_{\mu}(A_{r,-s})(t_{ijkl})=\frac{\partial t_{ijkl}}{\partial t_{rs}},\\
	&\varrho_{\mu}(A_{-r,s})(1)=t_{sr},\;\;\varrho_{\mu}(A_{-r,s})(t_{ij})=t_{ijsr},\;\;\varrho_{\mu}(A_{-r,s})(t_{ijkl})=0. \label{4.31}
\end{eqnarray}
Therefore, the subspace $\mathcal{T}$ of $\mathscr{X}$ spanned by
\begin{eqnarray}
	\{1,t_{ab},t_{ijkl}\mid 1\leqslant a<b\leqslant 5,\;1\leqslant i<j<k<l\leqslant 5 \} \label{4.32}
\end{eqnarray}is stable under $\varrho_{\mu}(\mathfrak{o}(10))$, and $1$ is annihilated by all the positive root vectors of $\mathfrak{o}(10)$. The first equation of $(\ref{4.28})$ shows that the weight of $1$ is $\omega_2$. So $1$ is a $\mathfrak{o}(10)$-singular vector. By Weyl's theorem of complete reducibility and that fact that $\dim L_{\mathfrak{o}(10)}(\omega_2)=16$, $\mathcal{T}\cong L_{\mathfrak{o}(10)}(\omega_2)$ and (\ref{4.32}) gives a basis. After reordering the elements of (\ref{4.32}) and adjusting the sign, we may write
\begin{align}
	(p_1,\cdots,p_{16})=(&\;t_{2345},-t_{1345},t_{1245},t_{1235},-t_{45},-t_{35},t_{1234},\nonumber\\
	&\;-t_{34},t_{25},t_{15},t_{24},t_{23},t_{14},t_{13},t_{12},1).
\end{align}
Using (\ref{4.28})-(\ref{4.31}), we can directly verify the following identity.
\begin{eqnarray}
	\varrho_{\mu}(u)p_i=-\sum_{j=1}^{16}c_{i,j}(u)p_j \label{4.34}
\end{eqnarray}
for any $u\in\mathscr{G}^{E_6}_0$, here $c_{i,j}(u)$'s are those coefficients in (\ref{247})-(\ref{248}), whose explicit values for $u=A_{r,s}$ can be found in Appendix \ref{app}. 
\begin{remark}
	In fact, the subspace $\mathcal{Y}=\sum_{k=1}^{16}\mathbb{C}y_k$ of $\mathscr{Y}$ is a $16$-dimensional $\mathfrak{o}(10)$-module with the action $\sigma$ (see (\ref{248})), and $y_1$ is the unique (up to a scalar) singular vector, of which the weight is $\omega_5$, hence $\mathcal{Y}\cong L_{\mathfrak{o}(10)}(\omega_5)$ as $\mathfrak{o}(10)$-modules. As a consequence, $\mathcal{T}$ and $\mathcal{Y}$ are the two mutually dual half-spin modules of $\mathfrak{o}(10)$.
\end{remark}

Define 
\begin{eqnarray}
	\Theta=\sum_{i=1}^{16}y_ip_i,\quad\vartheta_{u}=\sigma(u)\Theta,\quad\vartheta_{\eta_k}=\sigma(\eta_k)\Theta \label{4.35}
\end{eqnarray}
for $u\in\mathscr{G}^{E_6}_0$ and $k\in\overline{1,16}$. In particular, $\vartheta_{\widehat{\alpha}}=-3\Theta$.
Retain the basis of $\mathfrak{o}(10)$
\begin{eqnarray}
	\mathcal{B}_o=\{A_{i,j}\mid i,j\in\overline{1,5}\}\cup\{A_{i,-j},A_{-i,j}\mid 1\leqslant i<j\leqslant 5\}.
\end{eqnarray}
\begin{proposition}
	The set
	\begin{align}
		\mathcal{B}_2=\{p_\ell\mid \ell\in\overline{1,16}\}\cup\{\Theta\}\cup\{\vartheta_{u}\mid u\in\mathcal{B}_o\}\cup\{\vartheta_{\eta_k}\mid k\in\overline{1,16}\}
	\end{align}
	forms a basis of $V(\omega_2)$.
\end{proposition}
\begin{proof}
	 Firstly, each $\xi_i$ preserves $\operatorname{Span}_{\mathbb{C}}\mathcal{B}_2$ for all $i\in\overline{1,16}$. This follows from
	\begin{eqnarray}
		\pi_{\omega_2}(\xi_i)p_\ell=0,\;\pi_{\omega_2}(\xi_i)\Theta=p_i,\;\pi_{\omega_2}(\xi_i)\vartheta_{u}=\sum_{j=1}^{16}c_{i,j}(u)p_j,\;\pi_{\omega_2}(\xi_i)\vartheta_{\eta_k}=\vartheta_{[\xi_i,\eta_k]}, \label{4,38}
	\end{eqnarray}
	the last equality dues to the fact $\sigma(\eta_k)p_i=0$.
	
	Secondly, we consider the action of $\mathscr{G}^{E_6}_0$ on $\operatorname{Span}_{\mathbb{C}}\mathcal{B}_2$. By (\ref{248}), (\ref{r2}) and (\ref{4.34}), for $v\in\mathscr{G}^{E_6}_0$ we have
	\begin{eqnarray}
		\pi_{\omega_2}(v)(\Theta)=\sum_{i=1}^{16}(\sigma(v)(y_i)p_i+y_i\varrho_{\mu}(v)(p_i))=\sum_{i,j=1}^{16}c_{j,i}(v)y_jp_i-\sum_{i,j}^{16}c_{i,j}(v)y_ip_j=0.\label{4.38}
	\end{eqnarray}
	This equation now implies
	\begin{eqnarray}
		\pi_{\omega_2}(v)\vartheta_{\eta_k}=[\pi_{\omega_2}(v),\sigma(\eta_k)]\Theta=[\sigma(v),\sigma(\eta_k)]\Theta=\sum_{i=1}^{16}c_{i,k}(v)\vartheta_{\eta_i}, \label{4.40}
	\end{eqnarray}
	and
	\begin{eqnarray}
		\pi_{\omega_2}(v)\vartheta_u=[\pi_{\omega_2}(v),\sigma(u)]\Theta=[\sigma(v),\sigma(u)]\Theta=\vartheta_{[v,u]}\quad\text{for}\;u\in\mathcal{B}_o. \label{4.41}
	\end{eqnarray}
	Since $p_\ell$ has $y$-degree $0$, we also have
	\begin{eqnarray}
		\pi_{\omega_2}(u)p_\ell=\varrho_{\mu}(u)p_\ell=-\sum_{j=1}^{16}c_{\ell,j}(u)p_{j}.
	\end{eqnarray}
	Thus, the subspace $\operatorname{Span}_{\mathbb{C}}\mathcal{B}_2$ is stable under the action of $\pi_{\omega_2}(\mathscr{G}^{E_6}_0)$.
	
	Finally, we investigate $\pi_{\omega_2}(\eta_i)(\operatorname{Span}_{\mathbb{C}}\mathcal{B}_2)$. We have
	\begin{align}
		\pi_{\omega_2}(\eta_i)p_{\ell}=&\;\sigma(\eta_i)p_{\ell}+\sum_{r=1}^{16}y_r\varrho_{\mu}([\xi_r,\eta_i])p_\ell=-\sum_{r,j=1}^{16}c_{\ell,j}([\xi_r,\eta_i])y_rp_j\nonumber\\
		=&\;-\sum_{r,j=1}^{16}c_{r,j}([\xi_\ell,\eta_i])y_rp_{j}\nonumber\\
		=&\;-\sum_{j=1}^{16}\sigma([\xi_\ell,\eta_i])y_jp_j=-\vartheta_{[\xi_\ell,\eta_i]}, \label{4.43}
	\end{align}
	here the first line is by (\ref{3.14}) and (\ref{4.34}), the second line is derived from the Jacobi identity $[\xi_\ell,[\xi_r,\eta_i]]=[\xi_r,[\xi_\ell,\eta_i]]$, the third line is by (\ref{248}). The above identities then gives
    \begin{align}
    	\pi_{\omega_2}(\eta_i)\Theta=&\vartheta_{\eta_i}+\sum_{r=1}^{16}y_r\varrho_{\mu}([\xi_r,\eta_i])\Theta\nonumber\\
    	=&\;\vartheta_{\eta_i}-\sum_{r=1}^{16}y_r(\sigma(\xi_r)\sigma(\eta_i)-\sigma(\eta_i)\sigma(\xi_{r}))\Theta\nonumber\\
    	=&\;\vartheta_{\eta_i}-\sum_{r=1}^{16}y_{r}\partial_{y_{r}}\vartheta_{\eta_i}=-\vartheta_{\eta_i}, \label{4.44}
    \end{align}
    here the first line is by (\ref{3.14}) and then by (\ref{4.35}), the second line is by (\ref{4.38}), the last equality follows from $\sigma(\xi_\ell)=\partial_{y_{\ell}}$, $\sigma(\eta_i)p_\ell=0$ and the fact that $\vartheta_{\eta_i}$ has $y$-degree $2$. Furthermore, we have
    \begin{align}
    	\partial_{y_k}(\pi_{\omega_2}(\eta_i)\vartheta_u)=&\;\pi_{\omega_2}(\xi_k)\pi_{\omega_2}(\eta_i)\vartheta_u=\pi_{\omega_2}([\xi_k,\eta_i])\vartheta_u+\pi_{\omega_2}(\eta_i)\pi_{\omega_2}(\xi_k)\vartheta_u\nonumber\\
    	=&\;\vartheta_{[[\xi_k,\eta_i],u]}-\sum_{j=1}^{16}c_{k,j}(u)\vartheta_{[\xi_j,\eta_i]}\nonumber\\
    	=&\;-\sum_{r=1}^{16}c_{r,i}(u)\vartheta_{[\xi_k,\eta_r]}\nonumber\\
    	=&\;-\partial_{y_k}(\sum_{r=1}^{16}c_{r,i}(u)\vartheta_{\eta_r}),
    \end{align}
    here the second line is by the third equation in (\ref{4,38}), by (\ref{4.41}) and by (\ref{4.43}), the third line is the Jacobi identity, and the last line is by the last equation in (\ref{4,38}).
    Since $\pi_{\omega_2}(\eta_i)\vartheta_u+\sum_{r=1}^{16}c_{r,i}(u)\vartheta_{\eta_r}$ has $y$-degree $2$, it must be zero. Hence
    \begin{eqnarray}
    	\pi_{\omega_2}(\eta_i)\vartheta_u=-\sum_{r=1}^{16}c_{r,i}(u)\vartheta_{\eta_r}. \label{4.46}
    \end{eqnarray}
    In addition, we have
    \begin{align}
    	\partial_{y_k}(\pi_{\omega_2}(\eta_i)\vartheta_{\eta_j})=&\;\pi_{\omega_2}(\xi_k)\pi_{\omega_2}(\eta_i)\vartheta_{\eta_j}=\pi_{\omega_2}([\xi_k,\eta_i])\vartheta_{\eta_j}+\pi_{\omega_2}(\eta_i)\pi_{\omega_2}(\xi_k)\vartheta_{\eta_j}\nonumber\\
    	=&\;\sum_{r=1}^{16}c_{r,j}([\xi_k,\eta_i])\vartheta_{\eta_r}+\pi_{\omega_2}(\eta_i)\vartheta_{[\xi_k,\eta_j]}\nonumber\\
    	=&\;\sum_{r=1}^{16}c_{r,j}([\xi_k,\eta_i])\vartheta_{\eta_r}-\sum_{r=1}^{16}c_{r,i}([\xi_k,\eta_j])\vartheta_{\eta_r}\nonumber\\
    	=&\;0, 
    \end{align}
    here the second line is by the last equation in (\ref{4,38}) and by (\ref{4.40}), the third line is by (\ref{4.46}), and the last line is by the Jacobi identity. Since $\pi_{\omega_2}(\eta_i)\vartheta_{\eta_j}$ has $y$-degree $3$, it must be zero, that is
    \begin{eqnarray}
    	\pi_{\omega_2}(\eta_i)\vartheta_{\eta_j}=0. \label{4.48}
    \end{eqnarray}
    Combining (\ref{4.43}), (\ref{4.44}), (\ref{4.46}) and (\ref{4.48}) together, we see that the elements $\eta_i$'s also preserve $\operatorname{Span}_{\mathbb{C}}\mathcal{B}_2$. 
    
    Now Lemma \ref{lem4.1} implies that $\operatorname{Span}_{\mathbb{C}}\mathcal{B}_2=V(\omega_2)$. Since $|\mathcal{B}_2|\leqslant 78=\dim V(\omega_2)$, we conclude that $\mathcal{B}_2$ is a basis of $V(\omega_2)$.
\end{proof}

\subsection{Polynomial bases for remaining three fundamental modules}
The bases for $V(\omega_3)$, $V(\omega_4)$ and $V(\omega_5)$ arise from determinant expression in the differentiation of the bases of $V(\omega_1)$ and $V(\omega_6)$. We construct the required differential operators, which commute with the action of (simple) positive root vectors of $\mathscr{G}^{E_6}$ at first.

Retain the $16\times 16$ matrices $G_2,G_3,G_4,G_5$ and $G$ in (\ref{Gr1})-(\ref{Gr2}). Let $(g_1,\cdots,g_{16})$ be the first row of $G$. Note $g_1,\cdots,g_{16}\in\mathscr{X}$.
For example, one can compute $g_1=1$, $g_2=-x_{1,1}$, $g_3=-x_{2,1}+x_{1,1}x_{2,2}$, etc. Define   
\begin{align}
	\mho_1=&\;\partial_{x_{4,4}}+\sum_{r=1}^{3}(x_{3,r}\partial_{x_{4,r}}+x_{3,-r}\partial_{x_{4,-r}})-\big(\sum_{s=1}^3x_{3,s}x_{3,-s}\big)\partial_{x_{4,-4}},\label{4.51}\\
	\mho_3=&\;\partial_{x_{3,3}}+\sum_{r=1}^{2}(x_{2,r}\partial_{x_{3,r}}+x_{2,-r}\partial_{x_{3,-r}})-\big(\sum_{s=1}^2x_{2,s}x_{2,-s}\big)\partial_{x_{3,-3}},\\
	\mho_6=&\;\sum\limits_{k=1}^{16}g_k\partial_{y_k}.
\end{align}
 \begin{lemma}\label{lem4.8}
 	For any $\lambda\in \mathscr{H}^*$, $i\in\{1,3,6\}$ and $j\in\overline{1,6}$ we have
 	\begin{eqnarray}
 		[\mho_i,\pi_{\lambda}(E_{\alpha_j})]=0. \label{4.54}
 	\end{eqnarray}
 	Moreover, the constant terms of $\mho_if$ and $\pi_{\lambda}(E_{\alpha_i})f$ coincide for any $f\in\mathscr{C}$.
 \end{lemma}  
\begin{proof}
	Take any positive root vectors $A_{i,j}\in\mathfrak{o}(10)$ with $1\leqslant |j| <i\leqslant 5$, note that simple ones $E_{\alpha_1},\cdots,E_{\alpha_5}$ are of this form. Write $\pi_{\lambda}(A_{i,j})=\varrho_{\mu}(A_{i,j})+\sigma(A_{i,j})$ by Theorem \ref{thm3.1}. For positive root vectors $E_{\alpha_6}$, Theorem \ref{thm3.2} gives $\pi_{\lambda}(E_{\alpha_6})=\partial_{y_1}$.
	
	We firstly show that $\mho_6$ commutes with all the action of positive root vectors. Because every polynomial coefficient $g_k$ of $\mho_6$ has $y$-degree $0$, we immediately have
	\begin{eqnarray}
		[\partial_{y_1},\mho_6]=0.
	\end{eqnarray}
	Eq. (\ref{4.57}) gives
	\begin{eqnarray}
		\varrho_{\mu}(A_{i,j})(G)=-G\varrho_{\mu}(A_{i,j})(G^{-1})G=Gc(A_{i,j}).
	\end{eqnarray}
	Taking the first row of the matrices on both sides we get
	\begin{eqnarray}
		(\varrho_{\mu}(A_{i,j})g_1,\cdots,\varrho_{\mu}(A_{i,j})g_{16})=(g_1,\cdots,g_{16})c(A_{i,j}).\label{4.58}
	\end{eqnarray}
	On the other hand, (\ref{248}) gives
	\begin{align}
		[\sigma({A_{i,j}}),\partial_{y_k}]=-\sum_{l=1}^{16}c_{k,l}(A_{i,j})\partial_{y_l}. \label{4.59}
	\end{align}
	Now by (\ref{4.58})-(\ref{4.59}), we have
	\begin{eqnarray}
		[\pi_{\lambda}(A_{i,j}),\mho_6]=\sum_{k=1}^{16}\varrho_{\mu}(A_{i,j})g_k\partial_{y_k}-\sum_{k,l=1}^{16}g_kc_{k,l}(A_{i,j})\partial_{y_l}=0.
	\end{eqnarray}
	
	Next we proves (\ref{4.54}) for $i=1$. Because $\mho_1\in\mathscr{W}_{\mathscr{X}}$, it commutes with $\mathscr{W}_{\mathscr{Y}}$. In particular, $\mho_1$  commutes with $\pi_{\lambda}(E_{\alpha_6})$ and $\sigma(A_{i,j})$. It left to show that $\mho_1$ commutes with $\pi_{\lambda}(A_{i,j})$.
	Denote $\varepsilon_{r,s}$ by the $5\times 5$ matrix with $(r,s)$-entry being $1$ and others being $0$. By (\ref{2.8})-(\ref{cc}) and (\ref{4.51}) we have $\mho_1L=\varepsilon_{5,4}L$, and that $\mho_1C-\varepsilon_{5,4}C$ is upper triangular. Applying $\mho_1$ to the equality $CL^t+LC^t=0$ we obtain
	\begin{align}
		0=&\;(\mho_1C)L^t+CL^t\varepsilon_{5,4}^{t}+\varepsilon_{5,4}LC^t+L(\mho_1C^t)\nonumber\\=&\;L(\mho_1C^t-C^t\varepsilon_{5,4}^t)+(\mho_1C-\varepsilon_{5,4}C)L^t.
	\end{align}
	Hence $(\mho_1C-\varepsilon_{5,4}C)L^t$ is simultaneously skew-symmetric and upper triangular, and must be zero. Since $L$ is invertible, we have $\mho_1C=\varepsilon_{5,4}C$. Any positive root vector $A_{i,j}$ can be written in the form $\begin{pmatrix}
		P&Q\\
		0&-P^t
	\end{pmatrix}$ with $P,Q\in M_{5\times 5}(\mathbb{C})$, where $P$ is strictly lower triangular and $Q^{t}=-Q$. By (\ref{rho}) and the same arguments as above we can derive
	\begin{eqnarray}
		\varrho_{\mu}(A_{i,j})L=LP, \quad \varrho_{\mu}(A_{i,j})C=LQ-CP^t.
	\end{eqnarray}
	Consequently,
	\begin{eqnarray}
		&\mho_1\varrho_{\mu}(A_{i,j})L=\varepsilon_{5,4}LP=\varrho_{\mu}(A_{i,j})\mho_1L, \label{4,64}\\
		&\mho_1\varrho_{\mu}(A_{i,j})C=\varepsilon_{5,4}LQ-\varepsilon_{5,4}CP^t=\varrho_{\mu}(A_{i,j})\mho_1C. \label{4,65}
	\end{eqnarray}
	Since the lower triangular parts of $L$ and $C$ exhaust variables $x_{r,s}$ by (\ref{lc}), Eqs. (\ref{4,64})-(\ref{4,65}) implies that $\mho_1$ and $\varrho_{\mu}(A_{i,j})$ commute on $\mathscr{X}$. These two operators are independent of $y_k$, so they also commute on $\mathscr{C}$. Consequently,
	\begin{eqnarray}
		[\mho_1,\pi_{\lambda}(A_{i,j})]=[\mho_1,\varrho_{\mu}(A_{i,j})]+[\mho_1,\sigma(A_{i,j})]=0.
	\end{eqnarray}
	The proof for $i=3$ is done by substituting $\varepsilon_{5,4}$ to $\varepsilon_{4,3}$ in above arguments, we omit here.

	Taking the constant term is equivalent to impose all variables to be zero. Under such specialization, the matrix $G$ reduces to $I_{16}$, whenever $\mho_6$ reduces to $\partial_{y_1}$, which is exactly $\pi_{\lambda}(E_{\alpha_6})$. Meanwhile, the operator $\mho_1$ and $\mho_3$ reduce to $\partial_{x_{4,4}}$ and $\partial_{x_{3,3}}$, respectively. This coincides with $\pi_{\lambda}(E_{\alpha_1})$ and  $\pi_{\lambda}(E_{\alpha_3})$. Therefore, the second assertion of this lemma follows.
\end{proof}

\begin{lemma}
	For every $\nu\in \mathscr{H}^*$, $i\in\{1,3,6\}$ and $j\in\overline{1,6}$,
	\begin{align}
		\mho_i\pi_{\nu}(E_{\alpha_j})=&\;\pi_{\nu}(E_{\alpha_j})\mho_i,\label{mho1}\\
		\mho_i\pi_{\nu}(\alpha_j)=&\;\pi_{\nu-\alpha_i}(\alpha_j)\mho_i,\label{mho2}\\
		\mho_i\pi_{\nu}(E_{-\alpha_j})=&\;\pi_{\nu-\alpha_i}(E_{-\alpha_j})\mho_i-\delta_{i,j}\nu_i.\label{mho3}
	\end{align}
\end{lemma}
\begin{proof}
	The first formula is the first assertion in Lemma \ref{lem4.8}. For the second equation, we have
	\begin{align}
		&\;[\mho_i\pi_{\nu}(\alpha_j)-\pi_{\nu-\alpha_i}(\alpha_j)\mho_i,\pi_\nu(E_{\alpha_k})]\nonumber\\
		=&\;[\mho_i\pi_{\nu}(\alpha_j),\pi_\nu(E_{\alpha_k})]-[\pi_{\nu-\alpha_i}(\alpha_j)\mho_i,\pi_{\nu-\alpha_i}(E_{\alpha_k})]\nonumber\\
		=&\;(\alpha_k,\alpha_j)\mho_i\pi_{\nu}(E_{\alpha_k})-(\alpha_j,\alpha_k)\pi_{\nu-\alpha_i}(E_{\alpha_k})\mho_i=0,\quad (k\in\overline{1,6})
	\end{align}
	by (\ref{mho1}) and the fact that $\pi_{\nu}(E_{\alpha_k})$ is independent of $\nu$. On the other hand, for any $f\in\mathscr{C}$, taking the constant term, Lemma \ref{lem4.8} and (\ref{3.97})-(\ref{3.102}) give
	\begin{align}
		&\;\mho_i\pi_{\nu}(\alpha_j)f|_{x=0,y=0}=\pi_\nu(E_{\alpha_i})\pi_{\nu}(\alpha_j)f|_{x=0,y=0}\nonumber\\
		=&\;\pi_{\nu}(\alpha_j)\pi_\nu(E_{\alpha_i})f|_{x=0,y=0}-(\alpha_i,\alpha_j)\pi_\nu(E_{\alpha_i})f|_{x=0,y=0}\nonumber\\
		=&\;(\nu_j-(\alpha_i,\alpha_j))\mho_if|_{x=0,y=0}=\pi_{\nu-\alpha_i}(\alpha_j)\mho_if|_{x=0,y=0}.
	\end{align}
	Thus Lemma {\ref{lem3.8}} proves (\ref{mho2}). The Chevalley relation (\ref{Lie}) and the just proved identity (\ref{mho2}) give
	\begin{align}
		&\;[\mho_i\pi_{\nu}(E_{-\alpha_j})-\pi_{\nu-\alpha_i}(E_{-\alpha_j})\mho_i+\delta_{i,j}\nu_i,\pi_\nu(E_{\alpha_k})]\nonumber\\
		=&\;\delta_{j,k}\mho_i\pi_{\nu}(\alpha_j)-\delta_{j,k}\pi_{\nu-\alpha_i}(\alpha_j)\mho_i=0.
	\end{align}
	By the grading on the root height in the proof of Lemma \ref{lem3.8}, each simple negative root vector raises the degree by one, so its image always has zero constant terms.
	For any $f\in\mathscr{C}$, taking the constant of the polynomials we obtain
	\begin{align}
		&\;(\mho_i\pi_{\nu}(E_{-\alpha_j})-\pi_{\nu-\alpha_i}(E_{-\alpha_j})\mho_i)f|_{x=0,y=0}\nonumber\\
		=&\;\pi_\nu(E_{\alpha_i})\pi_{\nu}(E_{-\alpha_j})f|_{x=0,y=0}-\pi_{\nu-\alpha_i}(E_{-\alpha_j})\mho_if|_{x=0,y=0}\nonumber\\
		=&\;\pi_\nu(E_{-\alpha_j})\pi_\nu(E_{\alpha_i})f|_{x=0,y=0}-\delta_{i,j}\pi_{\nu}(\alpha_i)f|_{x=0,y=0}-\pi_{\nu-\alpha_i}(E_{-\alpha_j})\mho_if|_{x=0,y=0}\nonumber\\
		=&\;-\delta_{i,j}\nu_if|_{x=0,y=0}.
	\end{align}
	Now Lemma \ref{lem3.8} proves (\ref{mho3}).
\end{proof}

Recall the set $\mathcal{B}_1$ and $\mathcal{B}_6$ in Proposition \ref{prop4.5}, each set has $27$ elements. 
We may write $\mathcal{B}_1=\{v_1,\cdots,v_{27}\}$ and $\mathcal{B}_6=\{v_1',\cdots,v_{27}'\}$.
Define the following sets of polynomials given by determinants.
\begin{enumerate}[(i)]
	\item 
	\begin{eqnarray}
		\mathcal{B}_3=\left\lbrace \begin{vmatrix}
			v_i&v_j\\ \mho_1v_i&\mho_1v_j
		\end{vmatrix} \mid 1\leqslant i< j\leqslant 27\right\rbrace, 
	\end{eqnarray}
	\item 
	\begin{eqnarray}
		\mathcal{B}_4=\left\lbrace \begin{vmatrix}
			v_i&v_j&v_k\\ \mho_1v_i&\mho_1v_j&\mho_1v_k \\ \mho_3\mho_1v_i&\mho_3\mho_1v_j&\mho_3\mho_1v_k
		\end{vmatrix} \mid 1\leqslant i<j<k\leqslant 27\right\rbrace,
	\end{eqnarray}
	\item 
	\begin{eqnarray}
		\mathcal{B}_5=\left\lbrace \begin{vmatrix}
			v_i'&v_j'\\ \mho_6v_i'&\mho_6v_j'
		\end{vmatrix} \mid 1\leqslant i<j\leqslant 27\right\rbrace. 
	\end{eqnarray}
\end{enumerate}
\begin{proposition}
	The set $\mathcal{B}_i$ forms a basis of $V(\omega_i)$ for $i=3,4,5$.
\end{proposition}
\begin{proof}
	For (i) and (iii), we first prove the identity
	\begin{align}
		\pi_{2\nu-\alpha_i}(u)\left( \begin{vmatrix}
			q_1&q_2\\ \mho_iq_1&\mho_iq_2
		\end{vmatrix}\right) =\begin{vmatrix}
		\pi_{\nu}(u)q_1&q_2\\ \mho_i\pi_{\nu}(u)q_1&\mho_iq_2
		\end{vmatrix}+\begin{vmatrix}
		q_1&\pi_{\nu}(u)q_2\\ \mho_iq_1&\mho_i\pi_{\nu}(u)q_2
		\end{vmatrix}\label{det}
	\end{align}
	for all $u\in\mathscr{G}^{E_6}$, $q_1,q_2\in\mathscr{C}$ and $i\in\{1,3,6\}$. It suffices to show that (\ref{det}) holds for Chevalley generators $\{E_{\alpha_j},E_{-\alpha_j},\alpha_j\mid j\in\overline{1,6}\}$. By (\ref{mho3}), we get
	\begin{align}
		&\;\begin{vmatrix}
			\pi_{\nu}(E_{-\alpha_j})q_1&q_2\\ \mho_i\pi_{\nu}(E_{-\alpha_j})q_1&\mho_iq_2
		\end{vmatrix}+\begin{vmatrix}
			q_1&\pi_{\nu}(E_{-\alpha_j})q_2\\ \mho_iq_1&\mho_i\pi_{\nu}(E_{-\alpha_j})q_2
		\end{vmatrix}\nonumber\\
		=&\;\begin{vmatrix}
			\pi_{\nu}(E_{-\alpha_j})q_1&q_2\\ \pi_{\nu-\alpha_i}(E_{-\alpha_j})\mho_iq_1-\delta_{i,j}\nu_iq_1&\mho_iq_2
		\end{vmatrix}+\begin{vmatrix}
			q_1&\pi_{\nu}(E_{-\alpha_j})q_2\\ \mho_iq_1&\pi_{\nu-\alpha_i}(E_{-\alpha_j})\mho_iq_2-\delta_{i,j}\nu_iq_2
		\end{vmatrix}\nonumber\\
		=&\;\begin{vmatrix}
			\pi_{\nu}(E_{-\alpha_j})q_1&q_2\\ \pi_{\nu-\alpha_i}(E_{-\alpha_j})\mho_iq_1&\mho_iq_2
		\end{vmatrix}+\begin{vmatrix}
			q_1&\pi_{\nu}(E_{-\alpha_j})q_2\\ \mho_iq_1&\pi_{\nu-\alpha_i}(E_{-\alpha_j})\mho_iq_2
		\end{vmatrix}\nonumber\\
		=&\;\pi_{2\nu-\alpha_i}(E_{-\alpha_j})\big( \begin{vmatrix}
			q_1&q_2\\ \mho_iq_1&\mho_iq_2
		\end{vmatrix}\big), \label{4.64}
	\end{align}
	the last equality follows from Lemma \ref{lem3.3} and Leibniz rule.
 	The identity for each simple positive root vector and the coroot follows directly from (\ref{mho1}) and (\ref{mho2}). 
 	
 	Since $2\omega_1-\alpha_1=\omega_3$, taking $\nu=\omega_1$, $i=1$, $q_1,q_2\in\mathcal{B}_1$ in (\ref{det}) and applying Proposition \ref{prop4.5} (i), we derive that the subspace $\operatorname{Span}_{\mathbb{C}}\mathcal{B}_3$ of $\mathscr{C}$ spanned by $\mathcal{B}_3$ is stable under $\pi_{\omega_3}(\mathscr{G}^{E_6})$. Moreover, 
 	\begin{eqnarray}
 		\begin{vmatrix}
 			1&x_{4,4}\\\mho_11&\mho_1x_{4,4}
 		\end{vmatrix}=\begin{vmatrix}
 		1&x_{4,4}\\0&1
 		\end{vmatrix}=1\in\operatorname{Span}_{\mathbb{C}}\mathcal{B}_3.
 	\end{eqnarray}
 	Thus $\operatorname{Span}_{\mathbb{C}}\mathcal{B}_3=V(\omega_3)$ by Lemma \ref{lem4.1}. Then $|\mathcal B_3|\leqslant \binom{27}{2}=351=\dim V(\omega_3)$ implies that $\mathcal{B}_3$ is a basis of $V(\omega_3)$.
 	
 	Similarly, taking $\nu=\omega_6$, $i=6$ and $q_1,q_2\in\mathcal{B}_6$ in (\ref{det}), then applying $2\omega_6-\alpha_6=\omega_5$ and Proposition \ref{prop4.5} (ii), we see that $\pi_{\omega_5}(\mathscr{G}^{E_6})$ preserves $\operatorname{Span}_{\mathbb{C}}\mathcal{B}_5$. The element $1$ also belongs to $\operatorname{Span}_{\mathbb{C}}\mathcal{B}_5$ by
 	\begin{eqnarray}
 		\begin{vmatrix}
 			1&y_1\\\mho_61&\mho_6y_1
 		\end{vmatrix}=\begin{vmatrix}
 			1&y_1\\0&1
 		\end{vmatrix}=1.
 	\end{eqnarray}
 	Lemma \ref{lem4.1} gives $\operatorname{Span}_{\mathbb{C}}\mathcal{B}_5=V(\omega_5)$, and $|\mathcal B_5|\leqslant \binom{27}{2}=351=\dim V(\omega_5)$ implies that $\mathcal{B}_5$ is indeed a basis.
 	
 	For the assertion for $V(\omega_4)$, notice that $\omega_4=3\omega_1-2\alpha_1-\alpha_3$. Eqs. (\ref{4.10}) and (\ref{4,10}) show that every element in $\mathcal{B}_1$ is a polynomial in $y_1,\cdots,y_{16}$ and $x_{4,-4},\cdots,x_{4,4}$, so that $\mathcal{B}_1$ is annihilated by $\mho_3$. The formula (\ref{mho3}) consequently gives
 	\begin{eqnarray}
 		\mho_3\mho_1\pi_{\omega_1}(E_{-\alpha_j})=\pi_{\omega_1-\alpha_1-\alpha_3}(E_{-\alpha_j})\mho_3\mho_1-\delta_{3,j}\mho_1 \text{ on }\mathcal{B}_1.
 	\end{eqnarray}
 	Using the same method as in (\ref{det})-(\ref{4.64}) we can prove for any $u\in\mathscr{G}^{E_6}$ and $q_1,q_2,q_3\in\mathcal{B}_1$ that
 	\begin{align}
 		&\;\pi_{3\omega_1-2\alpha_1-\alpha_3}(u)\big(\begin{vmatrix}
 			q_1&q_2&q_3\\ \mho_1q_1&\mho_1q_2&\mho_1q_3 \\ \mho_3\mho_1q_1&\mho_3\mho_1q_2&\mho_3\mho_1q_3
 		\end{vmatrix}\big)\nonumber\\
 		=&\;\begin{vmatrix}
 			\pi_{\omega_1}(u)q_1&q_2&q_3\\ \mho_1\pi_{\omega_1}(u)q_1&\mho_1q_2&\mho_1q_3 \\ \mho_3\mho_1\pi_{\omega_1}(u)q_1&\mho_3\mho_1q_2&\mho_3\mho_1q_3
 		\end{vmatrix}+\begin{vmatrix}
 			q_1&\pi_{\omega_1}(u)q_2&q_3\\ \mho_1q_1&\mho_1\pi_{\omega_1}(u)q_2&\mho_1q_3 \\ \mho_3\mho_1q_1&\mho_3\mho_1\pi_{\omega_1}(u)q_2&\mho_3\mho_1q_3
 		\end{vmatrix}\nonumber\\
 		&\;+\begin{vmatrix}
 			q_1&q_2&\pi_{\omega_1}(u)q_3\\ \mho_1q_1&\mho_1q_2&\mho_1\pi_{\omega_1}(u)q_3 \\ \mho_3\mho_1q_1&\mho_3\mho_1q_2&\mho_3\mho_1\pi_{\omega_1}(u)q_3
 		\end{vmatrix}.
 	\end{align}
 	Therefore, $\operatorname{Span}_{\mathbb{C}}\mathcal{B}_4$ is stable under $\pi_{\omega_4}(\mathscr{G}^{E_6})$ by Proposition \ref{prop4.5} (i). Taking $1,x_{4,4},x_{4,3}\in\mathcal{B}_1$ we have
 	\begin{eqnarray}
 		\begin{vmatrix}
 			1&x_{4,4}&x_{4,3}\\
 			\mho_11&\mho_1x_{4,4}&\mho_1x_{4,3}\\
 			\mho_3\mho_11&\mho_3\mho_1x_{4,4}&\mho_3\mho_1x_{4,3}
 		\end{vmatrix}=\begin{vmatrix}
 		1&x_{4,4}&x_{4,3}\\
 		0&1&x_{3,3}\\
 		0&0&1
 		\end{vmatrix}=1\in\operatorname{Span}_{\mathbb{C}}\mathcal{B}_4.
 	\end{eqnarray}
 	Now Lemma \ref{lem4.1} yields $\operatorname{Span}_{\mathbb{C}}\mathcal{B}_4=V(\omega_4)$, and $|\mathcal B_4|\leqslant \binom{27}{3}=2925=\dim V(\omega_4)$ implies that $\mathcal{B}_4$ forms a basis of $V(\omega_4)$. This completes the proof.
 	
\end{proof}

\section{Homogeneous subspaces and graded Betti numbers} \label{section5}

In this section, we study the $y$-homogeneous subspace of $V(\lambda)$, where $\lambda$ is dominant integral, and determine their dimensions and the corresponding Hilbert series. The sixteen commuting negative root vectors $\eta_k$ raise the $y$-degree by one, making $V(\lambda)$ a finite length graded module over $U(\mathscr{G}^{E_6}_-)$. Using Kostant's homology theorem, we compute its graded Betti numbers and relate them to the homogeneous dimension formulas. We then characterize the highest weights for which the minimal graded free resolution is pure and establish a lower bound for the total Betti number when $\lambda\neq 0$, including all equality cases. For a reductive Lie algebra $\mathscr{G}$, we denote by $L_{\mathscr{G}}(\lambda)$ the abstract irreducible highest weight $\mathscr{G}$-module with highest weight $\lambda$.

\subsection{The dimension of the homogeneous subspaces}
For $n\in\mathbb{N}$, let $\mathscr{Y}[n]\subset\mathscr{Y}$ be the subspace of homogeneous polynomials with degree $n$. Then $\mathscr{C}[n]=\mathscr{Y}[n]\otimes \mathscr{X}$ is the subspace of $\mathscr{C}$ consisting of polynomials with degree $n$ in variable $y_1,\cdots,y_{16}$, which we will call {\it $y$-degree} for short. The {\it $y$-homogeneous subspace} of $V(\lambda)$ with $y$-degree $n$ is 
\begin{eqnarray}
	V(\lambda)[n]=V(\lambda)\cap\mathscr{C}[n].
\end{eqnarray}
Set
\begin{eqnarray}
	\chi_\lambda(q)=\sum_{n\in\mathbb{N}}\dim V(\lambda)[n]q^n, \label{5.2}
\end{eqnarray}
here $q$ is a formal variable. For $\nu\in \mathscr{H}^*$, denote by $V(\lambda)_{\nu}$ the weight space of $V(\lambda)$ with weight $\nu$. We use formal character $\operatorname{Ch}$ with formal variable $e$ as 
\begin{eqnarray}
	\operatorname{Ch}(V(\lambda))(e)=\sum_{\nu\in \mathscr{H}^*}\dim V(\lambda)_\nu e^{\nu}. \label{5.3}
\end{eqnarray}

\begin{lemma}\label{Lem5.1}
	We have
	\begin{eqnarray}
		\chi_\lambda(q)=\sum_{\nu\in \mathscr{H}^*}\dim V(\lambda)_\nu q^{(\lambda-\nu)(\widehat{\alpha}/3)}. \label{5.4}
	\end{eqnarray}
	Moreover, $\dim V(\lambda)[n]\neq 0$ if and only if
	\begin{eqnarray}
		0\leqslant n\leqslant N_\lambda, \label{5.5}
	\end{eqnarray}
	where $N_\lambda=2\lambda_1+2\lambda_2+3\lambda_3+4\lambda_4+3\lambda_5+2\lambda_6$.
\end{lemma}
\begin{proof}
	By (\ref{r1}), for every nonzero $f\in V(\lambda)_\nu$ we have
	\begin{eqnarray*}
		\mathcal{D}_y(f)=\frac{1}{3}(2\kappa-\pi_{\lambda}(\widehat{\alpha}))(f)=\frac{1}{3}(\lambda-\nu)(\widehat{\alpha})f,
	\end{eqnarray*}
	hence $V(\lambda)_\nu\subset V(\lambda)[(\lambda-\nu)(\widehat{\alpha}/3)]$, this implies (\ref{5.4}). Let $\lambda^*$ be the lowest weight of $V(\lambda)$, it is well known that 
	\begin{eqnarray}
		\lambda^*=-\lambda_6\omega_1-\lambda_2\omega_2-\lambda_5\omega_3-\lambda_4\omega_4-\lambda_3\omega_5-\lambda_1\omega_6.
	\end{eqnarray}
	The weight $\nu$ satisfies
	\begin{eqnarray}
		\lambda-\nu, \nu-\lambda^*\in Q^{+}_{E_6}. \label{5.6}
	\end{eqnarray}
	Hence (\ref{5.4}) shows that corresponding $y$-degree of $f$ lies between $0$ and $(\lambda-\lambda^*)(\widehat{\alpha}/3)$, and direct computation gives 
	\begin{eqnarray}
		(\lambda-\lambda^*)(\widehat{\alpha}/3)=2\lambda_1+2\lambda_2+3\lambda_3+4\lambda_4+3\lambda_5+2\lambda_6.
	\end{eqnarray}
	By the PBW theorem, we may write the lowest weight vector in $V(\lambda)$ as $E_{-\gamma_{i_r}}\cdots E_{-\gamma_{i_1}}1$ with each $\gamma_{i_j}\in\Phi_{E_6}^+$. Then every intermediate vector $E_{-\gamma_{i_s}}\cdots E_{-\gamma_{i_1}}1$ is nonzero with weight $\lambda-\sum_{j=1}^{s}\gamma_{i_j}$ for $s\in\overline{1,r}$. On the other hand, (\ref{CM}) and (\ref{2.104}) gives
	\begin{eqnarray}
		\alpha_i(\widehat{\alpha}/3)=\delta_{i,6} \quad\text{for}\;\; i\in\overline{1,6}, \label{5.9}
	\end{eqnarray}
	which yields
	\begin{align}
		\gamma(\widehat{\alpha}/3)=\begin{cases}
			0,&\text{if }\gamma\in\Phi_{D_5}^+,\\
			1,&\text{if }\gamma\in\Phi_{E_6}^{+}\setminus\Phi_{D_5}^{+}.
		\end{cases}\label{5.10}
	\end{align}
	 Therefore, along the sequence from $1$ to $E_{-\gamma_{i_r}}\cdots E_{-\gamma_{i_1}}1$, $y$-homogeneous polynomial with $y$-degree $n$ occurs for each $n\in\overline{0,N_\lambda}$. This completes the proof.
\end{proof}

Comparing (\ref{5.3}) and (\ref{5.4}) we find that $\chi_{\lambda}(q)$ can be recovered by replacing each $e^{\nu}$ by $q^{-\nu(\widehat{\alpha}/3)}$ in $\operatorname{Ch}(V(\lambda))(e)$ and then multiplying $q^{\lambda(\widehat{\alpha}/3)}$.
There are $20$ elements in $\Phi_{D_5}^+$ and $16$ elements in $\Phi_{E_6}^{+}\setminus\Phi_{D_5}^{+}$ (cf. (\ref{2.29})-(\ref{2.42})).

The Weyl's formula for the dimension of finite-dimensional irreducible modules gives the explicit expression
\begin{eqnarray}
	\dim L_{\mathscr{G}}(\mu)=\frac{\prod_{\alpha\in\Phi^+}\left\langle \mu+\rho,\alpha\right\rangle }{\prod_{\alpha\in\Phi^+}\left\langle \rho,\alpha\right\rangle }, \label{5.11}
\end{eqnarray}
where $\rho$ is the Weyl vector, see \cite[Corollary 24.3]{Hum} or \cite[Theorem 5.5.1]{Xu2}. Now Let $\rho_6$ be the Weyl vector of type $E_6$, precisely one has $\rho_6=\sum_{i=1}^6\omega_i$.
Denote by $W_{E_6}\supset W_{D_5}$ the Weyl groups of $\mathscr{G}^{E_6}\supset \mathscr{G}^{D_5}$ respectively. Let $s_i\in W_{E_6}$ be the reflection with respect to the simple root $\alpha_i$ for $i\in\overline{1,6}$, then $W_{E_6}$ is generated by $s_1,\cdots,s_6$, while $W_{D_5}$ is generated by $s_1,\cdots,s_5$. For each $w\in W_{E_6}$, define the $dot$ $action$ $\bullet$ of Weyl group on weights by $w\bullet\lambda=w(\lambda+\rho_6)-\rho_6$, and let $l(w)$ be the length of the element $w$. 

Choose the minimal set of representatives
\begin{eqnarray*}
	\widetilde{W}=\{w\in W_{E_6}\mid w^{-1}\alpha_i\in\Phi^{+}_{E_6}\text{ for all }i\in\overline{1,5}\}
\end{eqnarray*}
of the left coset $W_{D_5}\setminus W_{E_6}$.
Then $W_{E_6}=\sqcup_{w\in\widetilde{W}}W_{D_5}w$ and $|\widetilde{W}|=27$. One can produce all elements in $\widetilde{W}$ by computer, which can also be found in literatures (e.g. \cite{DZ1,DZ2}), we list the results with minimal length here.
\begin{eqnarray}
	&\operatorname{Id},\; s_6,\;s_6s_5,\;s_6s_5s_4,\;s_6s_5s_4s_2,\;s_6s_5s_4s_3,\;s_6s_5s_4s_2s_3,\;s_6s_5s_4s_3s_1,\label{w1}\\
	&s_6s_5s_4s_2s_3s_1,\;s_6s_5s_4s_2s_3s_4,\;s_6s_5s_4s_2s_3s_1s_4,\;s_6s_5s_4s_2s_3s_4s_5,\\
	&s_6s_5s_4s_2s_3s_1s_4s_3,s_6s_5s_4s_2s_3s_1s_4s_5,\;s_6s_5s_4s_2s_3s_4s_5s_6,\\
	&s_6s_5s_4s_2s_3s_1s_4s_3s_5,\;s_6s_5s_4s_2s_3s_1s_4s_5s_6,\;s_6s_5s_4s_2s_3s_1s_4s_3s_5s_4,\\
	&s_6s_5s_4s_2s_3s_1s_4s_3s_5s_6,\;s_6s_5s_4s_2s_3s_1s_4s_3s_5s_4s_2,\;s_6s_5s_4s_2s_3s_1s_4s_3s_5s_4s_6,\\
	&s_6s_5s_4s_2s_3s_1s_4s_3s_5s_4s_2s_6,\;s_6s_5s_4s_2s_3s_1s_4s_3s_5s_4s_6s_5,\\
	&s_6s_5s_4s_2s_3s_1s_4s_3s_5s_4s_2s_6s_5,\;s_6s_5s_4s_2s_3s_1s_4s_3s_5s_4s_2s_6s_5s_4,\\
	&s_6s_5s_4s_2s_3s_1s_4s_3s_5s_4s_2s_6s_5s_4s_3,\;s_6s_5s_4s_2s_3s_1s_4s_3s_5s_4s_2s_6s_5s_4s_3s_1.\label{w8}
\end{eqnarray}

For each $w\in\widetilde{W}$ and $i\in\overline{1,5}$, we have
\begin{eqnarray}
	(w\bullet\lambda)(\alpha_i)=(w(\lambda+\rho_6)-\rho_6)(\alpha_i)=(\lambda+\rho_6)(w^{-1}\alpha_i)-1\in\mathbb{N}, \label{5.12}
\end{eqnarray}
and $\lambda-w\bullet\lambda\in\sum_{i=1}^{6}\mathbb{N}\alpha_i$, whence
\begin{eqnarray}
	c_{\lambda}(w):=(\lambda-w\bullet\lambda)(\widehat{\alpha}/3)\in\mathbb{N}.\label{5.13}
\end{eqnarray}
Eq. (\ref{5.12}) suggests that $w\bullet\lambda|_{\mathscr{H}^{D_5}}$
is a dominant integral weight of $\mathscr{G}^{D_5}$.
\begin{theorem}\label{Thm5.2}
	For every dominant integral weight $\lambda\in \mathscr{H}^*$, we have
	\begin{eqnarray}
		(1-q)^{16}\chi_{\lambda}(q)=\sum_{w\in\widetilde{W}}(-1)^{l(w)}\dim L_{\mathscr{G}^{E_6}_0}(w\bullet\lambda)q^{c_{\lambda}(w)}. \label{5.14}
	\end{eqnarray}
	Consequently, for every $n\in\mathbb{N}$,
	\begin{eqnarray}
		\dim V(\lambda)[n]=\sum\limits_{w\in\widetilde{W},\atop c_{\lambda}(w)\leqslant n} (-1)^{l(w)}\dim L_{\mathscr{G}^{E_6}_0}(w\bullet\lambda)\binom{n-c_{\lambda}(w)+15}{15}. \label{5.15}
	\end{eqnarray}
\end{theorem}
\begin{proof}
	The parabolic form of the Weyl character formula gives
	\begin{eqnarray}
		\operatorname{Ch}(V(\lambda))(e)=\frac{\sum\limits_{w\in\widetilde{W}}(-1)^{l(w)}\operatorname{Ch}(L_{\mathscr{G}^{E_6}_0}(w\bullet\lambda))(e)}{\prod\limits_{\beta\in\Phi^+_{E_6}\setminus\Phi^+_{D_5}}(1-e^{-\beta})}, \label{5.16}
	\end{eqnarray}
	 (see e.g. \cite{GKRS}). Since every weight of $L_{\mathscr{G}^{E_6}_0}(w\bullet\lambda)$ is of the form $w\bullet\lambda-\nu$ with $\nu\in\sum_{i=1}^5\mathbb{N}\alpha_i$, Eq. (\ref{5.9}) implies that
	\begin{eqnarray}
		\operatorname{Ch}(L_{\mathscr{G}^{E_6}_0}(w\bullet\lambda))(e)\mapsto q^{-(w\bullet\lambda)(\widehat{\alpha}/3)}\dim L_{\mathscr{G}^{E_6}_0}(w\bullet\lambda)
	\end{eqnarray}
	under the specialization $e^{-\nu}\mapsto q^{\nu(\widehat{\alpha}/3)}=1$ for all $\nu\in\sum_{i=1}^5\mathbb{N}\alpha_i$.
	According to Lemma \ref{Lem5.1} and the observation thereafter, we obtain
	\begin{eqnarray}
		q^{-\lambda(\widehat{\alpha}/3)}\chi_\lambda(q)=\frac{\sum\limits_{w\in\widetilde{W}}(-1)^{l(w)}q^{-(w\bullet\lambda)(\widehat{\alpha}/3)}\dim L_{\mathscr{G}^{E_6}_0}(w\bullet\lambda)}{(1-q)^{16}},
	\end{eqnarray}
	which yields (\ref{5.14}). Then the expansion of the binomial
	\begin{eqnarray}
		(1-q)^{-16}=\sum_{r\in\mathbb{N}}\binom{r+15}{15}q^r
	\end{eqnarray}
	gives (\ref{5.15}).
\end{proof}

\subsection{Graded Betti numbers}
Since $\mathscr{G}^{E_6}_-=\sum_{k=1}^{16}\mathbb{C}\eta_k$ is abelian, the universal enveloping algebra $U(\mathscr{G}^{E_6}_-)$ is a free commutative associative $\mathbb{C}$-algebra, that is,
\begin{eqnarray}
	\mathscr{S}:=\mathbb{C}[z_1,\cdots,z_{16}]\cong U(\mathscr{G}^{E_6}_-).
\end{eqnarray}
We may regard $V(\lambda)$ as a $\mathscr{S}$-module via
\begin{eqnarray}
	z_if=\pi_{\lambda}(\eta_i)(f)\quad (f\in V(\lambda)).
\end{eqnarray}
Grading $\mathscr{S}$ as $\mathscr{S}=\bigoplus_{n\in\mathbb{N}}\mathscr{S}[n]$, where $\mathscr{S}[n]$ is the homogeneous subspace of $\mathscr{S}$ with degree $n$ as usual.
By (\ref{5.4}) and (\ref{5.10}), we have
\begin{eqnarray}
	z_iV(\lambda)[n]\subset V(\lambda)[n+1]\quad(n\in\mathbb{N}),
\end{eqnarray}
hence $V(\lambda)=\bigoplus_{{n}=0}^{N_{\lambda}}V(\lambda)[n]$ is a graded $\mathscr{S}$-module with respect to the degree of $\mathscr{S}$. Because $V(\lambda)$ is finite-dimensional, a power of the maximal ideal $\mathfrak{m}$ of $\mathscr{S}$ generated by $\{z_1,\cdots,z_{16}\}$ annihilates it. Thus $V(\lambda)$ has finite length and codimension $16$, and the Hilbert series of it is exactly $\chi_\lambda(q)$ defined in (\ref{5.2}). 

The {\it Betti numbers} are $b_{i,j}(\lambda)=\dim_{\mathbb{C}}\operatorname{Tor}_i^{\mathscr{S}}(\mathbb{C},V(\lambda))_j$ and $b_i(\lambda)=\sum_{j\in\mathbb{Z}}b_{i,j}(\lambda)$ for $i\in\mathbb{N}$ and $j\in\mathbb{Z}$. Equivalently, suppose that
\begin{eqnarray}
	\cdots\longrightarrow F_2\longrightarrow F_1\longrightarrow F_0\longrightarrow V(\lambda)\longrightarrow 0 \label{reso}
\end{eqnarray}
is the minimal graded free resolution of $V(\lambda)$, then one has it $i$-th free $\mathscr{S}$-module
\begin{eqnarray}
	F_i\cong \bigoplus_{j\in\mathbb{Z}}\mathscr{S}(-j)^{b_{i,j}(\lambda)}, \label{ith}
\end{eqnarray}
Here $\mathscr{S}(-j)$ is the graded free $\mathscr{S}$-module of rank $1$ with $\mathscr{S}(-j)_n=\mathscr{S}_{n-j}$, so that its canonical generator has degree $j$.
The resolution (\ref{reso}) is called {\it pure} if each nonzero $F_i$ is generated by elements of the same degree. Equivalently, for each $i\in\mathbb{N}$ the Betti numbers $b_{i,j}(\lambda)$ are nonzero for at most one $j\in\mathbb{Z}$.

The Chevalley-Eilenberg homology complex $\bigwedge^{\bullet}\mathscr{G}^{E_6}_-\otimes_{\mathbb{C}}V(\lambda)$
coincides with the Koszul complex computing $\operatorname{Tor}^{\mathscr{S}}(\mathbb{C},V(\lambda))$. Namely, we have
\begin{eqnarray}
	\operatorname{Tor}_i^{\mathscr{S}}(\mathbb{C},V(\lambda))\cong H_i(\mathscr{G}^{E_6}_-,V(\lambda)). \label{Ti}
\end{eqnarray}
Kostant's theorem (see e.g. \cite[Theorem 5.14]{K} or \cite[Theorem 5.1]{SVW}) gives
\begin{eqnarray}
	H_i(\mathscr{G}^{E_6}_-,V(\lambda))\cong\bigoplus_{w\in\widetilde{W},\atop l(w)=i} L_{\mathscr{G}^{E_6}_0}(w\bullet\lambda) \label{Hi} 
\end{eqnarray}
The component of degree $j$ in $\bigwedge^{i}\mathscr{G}^{E_6}_-\otimes_{\mathbb{C}}V(\lambda)$ is precisely $\bigwedge^{i}\mathscr{G}^{E_6}_-\otimes_{\mathbb{C}}V(\lambda)[j-i]$, on which the eigenvalue of $\widehat{\alpha}/3$ is
\begin{eqnarray}
	-i+\lambda(\widehat{\alpha}/3)-(j-i)=\lambda(\widehat{\alpha}/3)-j.
\end{eqnarray}
Since $\widehat{\alpha}$ commutes with $\mathscr{G}^{E_6}_0$, the eigenvalue of $\widehat{\alpha}/3$ on each $L_{\mathscr{G}^{E_6}_0}(w\bullet\lambda)$ is
\begin{eqnarray}
	(w\bullet\lambda)(\widehat{\alpha}/3)=\lambda(\widehat{\alpha}/3)-c_{\lambda}(w)
\end{eqnarray}
(cf. (\ref{5.13})). Therefore, (\ref{Ti})-(\ref{Hi}) gives 
\begin{eqnarray}
	b_{i,j}(\lambda)=\sum\limits_{w\in\widetilde{W},l(w)=i,\atop c_{\lambda}(w)=j}\dim L_{\mathscr{G}^{E_6}_0}(w\bullet\lambda)\quad\text{and}\quad b_{i}(\lambda)=\sum\limits_{w\in\widetilde{W},l(w)=i}\dim L_{\mathscr{G}^{E_6}_0}(w\bullet\lambda). \label{bij}
\end{eqnarray}
The list (\ref{w1})-(\ref{w8}) implies that $b_{i,j}(\lambda)=0$ for all $i>16$. Applying Theorem \ref{Thm5.2}, we get
\begin{corollary}
	For every dominant integral weight $\lambda\in \mathscr{H}^*$,
	\begin{eqnarray}
		(1-q)^{16}\chi_{\lambda}(q)=\sum_{i,j}(-1)^{i}b_{i,j}(\lambda)q^{j}.\label{5.38}
	\end{eqnarray}
\end{corollary}

Define 
\begin{eqnarray}
	\mathcal{B}_{\lambda}(p,q)=\sum_{i,j}b_{i,j}(\lambda)p^iq^j
\end{eqnarray}
with formal variable $p$ and $q$.
\begin{lemma}\label{Lem5.4}
	For positive integer $m$ and $\Lambda=m\lambda+(m-1)\rho_6$,
	\begin{eqnarray}
		\mathcal{B}_{\Lambda}(p,q)=m^{20}\mathcal{B}_{\lambda}(p,q^m).\label{bpq}
	\end{eqnarray}
	Equivalently, 
	\begin{align}
		b_{i,j}(\Lambda)=\begin{cases}
			m^{20}b_{i,j/m}(\lambda),&\text{if }m\mid j,\\
			0,&\text{if }m\nmid j.
		\end{cases}\label{bij2}
	\end{align}
\end{lemma}
\begin{proof}
	Note $\Lambda+\rho_6=m(\lambda+\rho_6)$.
	For $\alpha\in\Phi^{+}_{D_5}$, the numerator of the Weyl dimension formula (\ref{5.11}) at $w\bullet\Lambda$ is
	\begin{align}
		&\left\langle w\bullet\Lambda+\rho_6,\alpha\right\rangle =\left\langle \Lambda+\rho_6,w^{-1}\alpha \right\rangle\nonumber\\
		=&m\left\langle \lambda+\rho_6,w^{-1}\alpha \right\rangle=m\left\langle w\bullet\lambda+\rho_6,\alpha\right\rangle,
 	\end{align}
	this yields $c_{\Lambda}(w)=mc_\lambda(w)$ and $\dim L_{\mathscr{G}_0^{E_6}}(w\bullet\Lambda)=m^{20}\dim L_{\mathscr{G}_0^{E_6}}(w\bullet\lambda)$.
	Therefore, Eq. (\ref{bij}) gives (\ref{bij2}), and (\ref{bpq}) follows by definition.
\end{proof}

\begin{theorem}\label{Thm5.5}
	The minimal graded free resolution $F_{\bullet}$ of $V(\lambda)$ is pure if and only if there exists $k\in\mathbb{N}$ such that $\lambda=k\rho_6$. Under such a condition,
	\begin{eqnarray}
		F_i\cong\mathscr{S}(-(k+1)i)^{(k+1)^{20}\binom{16}{i}},\quad (i\in\mathbb{N}). \label{pure}
	\end{eqnarray}
\end{theorem}
\begin{proof}
	Since $V(0)\cong\mathbb{C}$, the minimal graded free resolution is precisely the Koszul resolution, hence $b_{i,j}(0)=\delta_{i,j}\binom{16}{i}$. Then Lemma \ref{Lem5.4} proves (\ref{pure}) for $\lambda=k\rho_6$, and the sufficiency also follows.
	
	For necessity, suppose that $F_{\bullet}$ is pure. Since $\dim L_{\mathscr{G}_0^{E_6}}(w\bullet\lambda)$ is positive, Eq. (\ref{bij}) forces $c_{\lambda}(w)=c_{\lambda}(w')$ whenever $l(w)=l(w')$. We choose the elements with length $4$, $5$, $6$, $7$, $8$ in the list (\ref{w1})-(\ref{w8}) and derive
	\begin{eqnarray}
	 &c_{\lambda}(s_6s_5s_4s_2)-c_{\lambda}(s_6s_5s_4s_3)=\lambda_2-\lambda_3,\\
	 &c_{\lambda}(s_6s_5s_4s_2s_3)-c_{\lambda}(s_6s_5s_4s_3s_1)=\lambda_2-\lambda_1,\\
	 &c_{\lambda}(s_6s_5s_4s_2s_3s_1)-c_{\lambda}(s_6s_5s_4s_2s_3s_4)=\lambda_1-\lambda_4,\\
	 &c_{\lambda}(s_6s_5s_4s_2s_3s_1s_4)-c_{\lambda}(s_6s_5s_4s_2s_3s_4s_5)=\lambda_1-\lambda_5,\\
	 &c_{\lambda}(s_6s_5s_4s_2s_3s_1s_4s_5)-c_{\lambda}(s_6s_5s_4s_2s_3s_4s_5s_6)=\lambda_1-\lambda_6.
	\end{eqnarray}
	These equations give $\lambda_1=\cdots=\lambda_6$, that is, $\lambda=k\rho$ for some $k\in\mathbb{N}$.
\end{proof}

\begin{theorem}
	For two dominant integral weight $\lambda,\nu\in \mathscr{H}^*$ with $\lambda_r\geqslant\nu_r$ for all $r\in\overline{1,6}$, we have $b_i(\lambda)\geqslant b_i(\nu)$ for all $i\in\mathbb{N}$. In particular,
	\begin{eqnarray}
		b_i(\lambda)\geqslant\binom{16}{i}. \label{5.51}
	\end{eqnarray}
	If $\lambda\neq0$, then we have
	\begin{eqnarray}
		\sum_{i=0}^{16}b_i(\lambda)\geqslant 435204, \label{5.52}
	\end{eqnarray}
	and the equality holds if and only if $\lambda=\omega_1$ or $\omega_6$.
\end{theorem}
\begin{proof}
	By definition, $w^{-1}\alpha\in\Phi^+_{E_6}$ for each $\alpha\in\Phi^{+}_{D_5}$. Hence $\left\langle w\bullet\lambda+\rho,\alpha\right\rangle =\left\langle \lambda+\rho, w^{-1}\alpha \right\rangle $ is a linear polynomial in $\lambda_1,\cdots,\lambda_6$  with nonnegative integer coefficients, and the constant term is positive. Therefore,
	\begin{eqnarray}
		\left\langle w\bullet\lambda+\rho,\alpha\right\rangle\geqslant\left\langle w\bullet\nu+\rho,\alpha\right\rangle.
	\end{eqnarray}
	By the Weyl's formula of dimension (\ref{5.11}), we get
	\begin{eqnarray}
		\dim L_{\mathscr{G}_0^{E_6}}(w\bullet\lambda)\geqslant \dim L_{\mathscr{G}_0^{E_6}}(w\bullet\nu).
	\end{eqnarray}
	Then (\ref{bij}) implies $b_{i}(\lambda)\geqslant b_{i}(\nu)$ for all $i$. Taking $\nu=0$ we obtain
	\begin{eqnarray}
		b_i(\lambda)\geqslant b_i(0)=\binom{16}{i}, 
	\end{eqnarray}
	here the last equality is explained in the proof of Theorem \ref{Thm5.5}. 
	
	If $\lambda$ is nonzero, then $\lambda_r>0$ for some $r\in\overline{1,6}$, so $b_i(\lambda)\geqslant b_i(\omega_r)$. Summing the Weyl dimensions in (\ref{bij}) over the $27$ representatives listed in (\ref{w1})-(\ref{w8}) gives
	\begin{eqnarray}
		&\sum\limits_{i=0}^{16}b_i(\omega_1)=\sum\limits_{i=0}^{16}b_i(\omega_6)=435204,\label{5.55}\\
		&\sum\limits_{i=0}^{16}b_i(\omega_3)=\sum\limits_{i=0}^{16}b_i(\omega_5)=1690500,\\
		&\sum\limits_{i=0}^{16}b_i(\omega_2)=970200,\quad\sum\limits_{i=0}^{16}b_i(\omega_4)=5850000.\label{5.57}
	\end{eqnarray}
	Hence (\ref{5.52}) follows. If one of $\lambda_2,\lambda_3,\lambda_4,\lambda_5$ is nonzero, then (\ref{5.55})-(\ref{5.57}) shows that (\ref{5.52}) holds without equality.
	Otherwise, we assume firstly that $\lambda_1$ is nonzero, the arguments above shows that the equality in (\ref{5.52}) holds if and only if
	\begin{eqnarray}
		\left\langle w\bullet\lambda+\rho,\alpha\right\rangle=\left\langle w\bullet\omega_1+\rho,\alpha\right\rangle
	\end{eqnarray}
	for all $w\in\widetilde{W}$ and $\alpha\in\Phi^+_{D_5}$. Taking $w=\operatorname{Id},\alpha=\alpha_1$ we deduce $\lambda_1=1$, taking $w=s_6,\alpha=\alpha_5$ we deduce $\lambda_6=0$, hence $\lambda=\omega_1$. For the case $\lambda_6\neq0$, we can show $\lambda=\omega_6$ by the same method. Consequently, the equality in (\ref{5.52}) holds if and only if $\lambda=\omega_1$ or $\omega_6$.
\end{proof}

\begin{remark}
	The results established in this subsection gives
	\begin{eqnarray}
		b_i(\lambda)\geqslant b_i(0)=\binom{16}{i}\qquad (i\in\overline{1,16}),
	\end{eqnarray}
	since $V(0)\cong \mathscr{S}/(z_1,\cdots,z_{16})$ has the Koszul resolution. Hence every $V(\lambda)$ satisfies the Buchsbaum–Eisenbud–Horrocks conjecture (see e.g. \cite[Conjecture 4.1]{BG}). The question proposed by Charalambous, Evans and Miller formulated in \cite[Question B]{BG} predicts that
	\begin{eqnarray}
		\sum_{i\in\mathbb{N}}b_i(M)\geqslant 2^{c}+2^{c-1}
	\end{eqnarray}
	for a nonzero finitely generated graded module $M$ of codimension $c>1$ over a polynomial algebra $\mathscr{P}$, provided that $M$ is not a quotient of $\mathscr{P}$ by a homogeneous regular sequence up to a shift of grading.  The lower bound in (\ref{5.52}) shows that
	\begin{eqnarray}
		\sum_{i=0}^{16}b_i(\lambda)\geqslant 435204>98304=2^{16}+2^{15}\quad\text{for any }\lambda\neq0.
	\end{eqnarray}
	The constant $435204$ is optimal for the family of $\mathscr{S}$-modules we considered here.
\end{remark}

\subsection{Combinatorial identities}
We now apply the preceding formulas to selected highest weights. Set
\begin{eqnarray}
	\Lambda=(k+1)\lambda+k\rho_6
\end{eqnarray}
with $k\in\mathbb{N}$ throughout this subsection. We use the convention that $\binom{m}{n}=0$ for integers $m<n$.
\begin{proposition}\label{prop5.6}
	For any dominant integral weight $\lambda\in \mathscr{H}^*$ and every $k\in\mathbb{N}$, 
	\begin{eqnarray}
		(k+1)^{20}\big(\sum_{r=0}^kq^r\big)^{16}\chi_\lambda(q^{k+1})=\chi_{\Lambda}(q).\label{5.63}
	\end{eqnarray}
	Consequently, for every $n\in\mathbb{N}$
	\begin{align}
		(k+1)^{20}\sum_{r=0}^{N_{\lambda}}\sum\limits_{j=0}^{16}(-1)^j\binom{16}{j}\binom{n-(k+1)(r+j)+15}{15}\dim V(\lambda)[r]=\dim V(\Lambda)[n].\label{5.33}
	\end{align}
\end{proposition}
\begin{proof}
	Taking $p=-1$ and $m=k+1$ in (\ref{5.38}) and (\ref{bpq}), we get
	\begin{align}
		(1-q)^{16}\chi_{\Lambda}(q)&=\mathcal{B}_{\Lambda}(-1,q)=(k+1)^{20}\mathcal{B}_{\lambda}(-1,q^{k+1})\nonumber\\
		&=(k+1)^{20}(1-q^{k+1})^{16}\chi_{\lambda}(q^{k+1}).\label{5.65}
	\end{align}
	Since 
	\begin{eqnarray}
		\frac{1-q^{k+1}}{1-q}=1+q+\cdots +q^k,
	\end{eqnarray}
	then (\ref{5.65}) gives (\ref{5.63}). Using the expansion
	\begin{align}
		\big(\sum_{r=0}^{k}q^r\big)^{16}=(1-q^{k+1})^{16}(1-q)^{-16}=\sum_{j=0}^{16}\sum_{r\in\mathbb{N}}(-1)^{j}\binom{16}{j}\binom{r+15}{15}q^{(k+1)j+r}
	\end{align}
	to extract the coefficients of $q^n$ on both sides of (\ref{5.63}) yields (\ref{5.33}).
\end{proof}
\begin{corollary}\label{Cor5.9}
	With the notation of Proposition \ref{prop5.6}, we have
	\begin{align}
		\sum_{r=0}^{N_\lambda}\sum_{j=0}^{16}(-1)^j\binom{16}{j}\binom{(k+1)(N_\lambda+16-r-j)}{16}\dim V(\lambda)[r]=(k+1)^{16}\dim V(\lambda).\label{5.36}
	\end{align}
\end{corollary}
\begin{proof}
	The Weyl dimension formula gives
	\begin{eqnarray}
		\dim V(\Lambda)=(k+1)^{36}\dim V(\lambda),
	\end{eqnarray}
	because every positive root factor is multiplied by $k+1$.
	Lemma \ref{Lem5.1} shows that the maximal $N_{\Lambda}$ such that $V(\Lambda)[N_{\Lambda}]$ is nonzero is
	\begin{eqnarray}
		2\Lambda_1+2\Lambda_2+3\Lambda_3+4\Lambda_4+3\Lambda_5+2\Lambda_6=(k+1)N_\lambda+16k.
	\end{eqnarray}
	Sum (\ref{5.33}) over $n\in\overline{0,N_\Lambda}$, then (\ref{5.36}) follows by
	\begin{align}
		\sum_{n=0}^{N_\Lambda}\binom{n-(k+1)(r+j)+15}{15}=\binom{N_\Lambda-(k+1)(r+j)+16}{16}.
	\end{align}
	This completes the proof.
\end{proof}

Now we reduce to the $\mathscr{G}^{E_6}$-modules $V(k\rho_6)$, which are the only case admitting pure graded minimal free resolution as $\mathscr{S}$-modules as showing in Theorem \ref{Thm5.5}. Taking $\lambda=0$, we have $V(0)=\mathbb{C}$ and $\chi_{0}(q)=1$.
In particular,
\begin{align}
	\dim V(0)[n]=\begin{cases}
		1,&\text{if }n=0,\\
		0, &\text{if }n>0.
	\end{cases}
\end{align}
Proposition \ref{prop5.6} and Corollary \ref{Cor5.9} reduces to the following results.
\begin{proposition}\label{Prop5.3}
	For every $k,n\in\mathbb{N}$, we have
	\begin{eqnarray}
		\chi_{k\rho_6}(q)=(k+1)^{20}\big(\sum_{r=0}^{k}q^r\big)^{16},
	\end{eqnarray}
	and
	\begin{eqnarray}
		\dim V(k\rho_6)[n]=(k+1)^{20}\sum_{j=0}^{16}(-1)^{j}\binom{16}{j}\binom{n-(k+1)j+15}{15}.
	\end{eqnarray}
	In particular, we have the identity
	\begin{eqnarray}
		\sum_{j=0}^{16}(-1)^j\binom{16}{j}\binom{(16-j)(k+1)}{16}=(k+1)^{16}. \label{5.76}
	\end{eqnarray}
\end{proposition}
\begin{remark}
	The identity (\ref{5.76}) is indeed a specific case of the following classical identity
	\begin{eqnarray}
		\sum_{j=0}^{r}(-1)^j\binom{r}{j}\binom{(r-j)(k+1)}{r}=(k+1)^{r} \label{5.30}
	\end{eqnarray}
	with $r=16$. One can derive (\ref{5.30}) by taking the coefficients of $x^r$ in both sides of the binomial expansion
	\begin{eqnarray}
		\sum_{j=0}^r(-1)^j\binom{r}{j}\sum_{i=0}^{(r-j)(k+1)}\binom{(k+1)(r-j)}{i}x^{i}=((x+1)^{k+1}-1)^r.
	\end{eqnarray}
\end{remark}

\appendix
\section{Spin representation of $\mathscr{G}^{E_6}_0$ and its extension to $\mathscr{G}^{E_6}$} \label{app}
The Lie algebra homomorphism $\sigma:\mathscr{G}^{E_6}_0\longrightarrow\mathscr{W}_{\mathscr{Y}}$ in Section \ref{section2.3} is explicitly given as follows.
\begin{align}
	&\sigma(\alpha_6)=-2y_1\partial_{y_{1}}-\sum\limits_{i=2}^{9}y_i\partial_{y_i}-y_{11}\partial_{y_{11}}-y_{12}\partial_{y_{12}},\\
	&\sigma(A_{5,4})=y_4\partial_{y_7}+y_6\partial_{y_8}+y_9\partial_{y_{11}}+y_{10}\partial_{y_{13}},\\
	&\sigma(A_{4,3})=y_3\partial_{y_4}+y_5\partial_{y_6}+y_{11}\partial_{y_{12}}+y_{13}\partial_{y_{14}},\\
	&\sigma(A_{3,2})=-y_2\partial_{y_3}-y_6\partial_{y_9}-y_8\partial_{y_{11}}+y_{14}\partial_{y_{15}},\\
	&\sigma(A_{2,1})=-y_1\partial_{y_2}+y_9\partial_{y_{10}}+y_{11}\partial_{y_{13}}+y_{12}\partial_{y_{14}},\\
	&\sigma(A_{2,-1})=-y_3\partial_{y_5}-y_4\partial_{y_6}-y_7\partial_{y_8}+y_{15}\partial_{y_{16}},\\
	&\sigma(A_{5,3})=-y_3\partial_{y_7}-y_5\partial_{y_8}+y_9\partial_{y_{12}}+y_{10}\partial_{y_{14}},\\
	&\sigma(A_{4,2})=y_{2}\partial_{y_4}-y_5\partial_{y_9}+y_8\partial_{y_{12}}+y_{13}\partial_{y_{15}},\\
	&\sigma(A_{3,1})=-y_1\partial_{y_3}-y_6\partial_{y_{10}}-y_8\partial_{y_{13}}-y_{12}\partial_{y_{15}},\\
	&\sigma(A_{3,-1})=y_2\partial_{y_5}-y_4\partial_{y_9}-y_7\partial_{y_{11}}+y_{14}\partial_{y_{16}},\\
	&\sigma(A_{5,2})=-y_2\partial_{y_7}+y_5\partial_{y_{11}}+y_6\partial_{y_{12}}+y_{10}\partial_{y_{15}},\\
	&\sigma(A_{4,1})=y_1\partial_{y_4}-y_5\partial_{y_{10}}+y_{8}\partial_{y_{14}}-y_{11}\partial_{y_{15}},\\
	&\sigma(A_{4,-1})=-y_2\partial_{y_6}-y_{3}\partial_{y_9}+y_7\partial_{y_{12}}+y_{13}\partial_{y_{16}},\\
	&\sigma(A_{3,-2})=-y_1\partial_{y_5}+y_4\partial_{y_{10}}+y_7\partial_{y_{13}}+y_{12}\partial_{y_{16}},\\
	&\sigma(A_{5,1})=-y_1\partial_{y_7}+y_5\partial_{y_{13}}+y_6\partial_{y_{14}}-y_9\partial_{y_{15}},\\
	&\sigma(A_{5,-1})=y_2\partial_{y_8}+y_3\partial_{y_{11}}+y_4\partial_{y_{12}}+y_{10}\partial_{y_{16}},\\
	&\sigma(A_{4,-2})=y_1\partial_{y_6}+y_3\partial_{y_{10}}-y_7\partial_{y_{14}}+y_{11}\partial_{y_{16}},\\
	&\sigma(A_{5,-2})=-y_1\partial_{y_8}-y_3\partial_{y_{13}}-y_4\partial_{y_{14}}+y_9\partial_{y_{16}},\\
	&\sigma(A_{4,-3})=y_1\partial_{y_9}-y_2\partial_{y_{10}}+y_7\partial_{y_{15}}-y_8\partial_{y_{16}},\\
	&\sigma(A_{5,-3})=-y_1\partial_{y_{11}}+y_2\partial_{y_{13}}+y_4\partial_{y_{15}}-y_6\partial_{y_{16}},\\
	&\sigma(A_{5,-4})=y_1\partial_{y_{12}}-y_2\partial_{y_{14}}+y_3\partial_{y_{15}}-y_5\partial_{y_{16}},\label{2.101}\\
	&\sigma(A_{r,r})=\mathop{\textstyle\sum}\limits_{i=1}^{16}(\frac{1}{2}+a_{r,i})y_i\partial_{y_{i}}\quad(r\in\overline{1,5}),
\end{align}
where $a_{r,i}$'s are given in the following table. Moreover, for $1\leqslant |j|<i\leqslant 5$ the formulas $\sigma(A_{j,i})$ are obtained by replacing $y_p\partial_{y_q}$ by $y_q\partial_{y_p}$ in $\sigma(A_{i,j})$.
\begin{center}
	\renewcommand{\arraystretch}{1.15}
		\begin{tabular}{|c|c|c|c|c|c|c|c|c|c|c|c|c|c|c|c|c|}
			\hline
			$i$&$1$&$2$&$3$&$4$&$5$&$6$&$7$&$8$&$9$&$10$&$11$&$12$&$13$&$14$&$15$&$16$\\
			\hline
			$a_{1,i}$&$-1$&$0$&$0$&$0$&$-1$&$-1$&$0$&$-1$&$-1$&$0$&$-1$&$-1$&$0$&$0$&$0$&$-1$\\
			\hline
			$a_{2,i}$&$0$&$-1$&$0$&$0$&$-1$&$-1$&$0$&$-1$&$0$&$-1$&$0$&$0$&$-1$&$-1$&$0$&$-1$\\
			\hline
			$a_{3,i}$&$0$&$0$&$-1$&$0$&$-1$&$0$&$0$&$0$&$-1$&$-1$&$-1$&$0$&$-1$&$0$&$-1$&$-1$\\
			\hline
			$a_{4,i}$&$0$&$0$&$0$&$-1$&$0$&$-1$&$0$&$0$&$-1$&$-1$&$0$&$-1$&$0$&$-1$&$-1$&$-1$\\
			\hline
			$a_{5,i}$&$0$&$0$&$0$&$0$&$0$&$0$&$-1$&$-1$&$0$&$0$&$-1$&$-1$&$-1$&$-1$&$-1$&$-1$\\
			\hline
		\end{tabular}
\end{center}
Furthermore, its extension to $\mathscr{G}^{E_6}$ are given by
\begin{eqnarray}
	\sigma(\xi_i)=\partial_{y_i} \qquad(i\in\overline{1,16}),
\end{eqnarray}
and
\begin{eqnarray}
	&\sigma(\eta_1)=y_1\mathcal{D}_y-\zeta_1\partial_{y_{10}}-\zeta_2\partial_{y_{13}}-\zeta_3\partial_{y_{14}}-\zeta_4\partial_{y_{15}}-\zeta_{10}\partial_{y_{16}},\\
	&\sigma(\eta_2)=y_2\mathcal{D}_y-\zeta_1\partial_{y_9}-\zeta_2\partial_{y_{11}}-\zeta_3\partial_{y_{12}}+\zeta_5\partial_{y_{15}}-\zeta_9\partial_{y_{16}}, \\
	&\sigma(\eta_3)=y_3\mathcal{D}_y+\zeta_1\partial_{y_6}+\zeta_2\partial_{y_8}+\zeta_4\partial_{y_{12}}+\zeta_5\partial_{y_{14}}-\zeta_8\partial_{y_{16}}, \\
	&\sigma(\eta_4)=y_4\mathcal{D}_y-\zeta_1\partial_{y_5}+\zeta_3\partial_{y_8}-\zeta_4\partial_{y_{11}}-\zeta_5\partial_{y_{13}}+\zeta_7\partial_{y_{16}}, \\
	&\sigma(\eta_5)=y_5\mathcal{D}_y-\zeta_1\partial_{y_4}-\zeta_2\partial_{y_7}-\zeta_{10}\partial_{y_{12}}+\zeta_9\partial_{y_{14}}-\zeta_8\partial_{y_{15}}, \\
	&\sigma(\eta_6)=y_6\mathcal{D}_y+\zeta_1\partial_{y_3}-\zeta_3\partial_{y_7}+\zeta_{10}\partial_{y_{11}}-\zeta_9\partial_{y_{13}}+\zeta_7\partial_{y_{15}}, \\
	&\sigma(\eta_7)=y_7\mathcal{D}_y-\zeta_2\partial_{y_5}-\zeta_3\partial_{y_6}+\zeta_4\partial_{y_9}+\zeta_5\partial_{y_{10}}-\zeta_6\partial_{y_{16}}, \\  &\sigma(\eta_8)=y_8\mathcal{D}_y+\zeta_2\partial_{y_3}+\zeta_3\partial_{y_4}-\zeta_{10}\partial_{y_9}+\zeta_9\partial_{y_{10}}-\zeta_6\partial_{y_{15}}, \\
	&\sigma(\eta_9)=y_9\mathcal{D}_y-\zeta_1\partial_{y_2}+\zeta_4\partial_{y_7}-\zeta_{10}\partial_{y_{8}}-\zeta_8\partial_{y_{13}}+\zeta_7\partial_{y_{14}}, \\
	&\sigma(\eta_{10})=y_{10}\mathcal{D}_y-\zeta_1\partial_{y_1}+\zeta_5\partial_{y_7}+\zeta_9\partial_{y_8}+\zeta_8\partial_{y_{11}}-\zeta_7\partial_{y_{12}}, \\
	&\sigma(\eta_{11})=y_{11}\mathcal{D}_y-\zeta_2\partial_{y_2}-\zeta_4\partial_{y_4}+\zeta_{10}\partial_{y_6}+\zeta_8\partial_{y_{10}}-\zeta_6\partial_{y_{14}}, \\
	&\sigma(\eta_{12})=y_{12}\mathcal{D}_y-\zeta_3\partial_{y_2}+\zeta_4\partial_{y_3}-\zeta_{10}\partial_{y_5}-\zeta_7\partial_{y_{10}}+\zeta_6\partial_{y_{13}}, \\
	&\sigma(\eta_{13})=y_{13}\mathcal{D}_y-\zeta_2\partial_{y_1}-\zeta_5\partial_{y_4}-\zeta_9\partial_{y_6}-\zeta_8\partial_{y_9}+\zeta_6\partial_{y_{12}}, \\
	&\sigma(\eta_{14})=y_{14}\mathcal{D}_y-\zeta_3\partial_{y_1}+\zeta_5\partial_{y_3}+\zeta_9\partial_{y_5}+\zeta_7\partial_{y_9}-\zeta_6\partial_{y_{11}}, \\
	&\sigma(\eta_{15})=y_{15}\mathcal{D}_y-\zeta_4\partial_{y_1}+\zeta_5\partial_{y_2}-\zeta_8\partial_{y_5}+\zeta_7\partial_{y_6}-\zeta_6\partial_{y_8}, \\
	&\sigma(\eta_{16})=y_{16}\mathcal{D}_y-\zeta_{10}\partial_{y_1}-\zeta_9\partial_{y_2}-\zeta_8\partial_{y_3}+\zeta_7\partial_{y_4}-\zeta_6\partial_{y_7}.
\end{eqnarray}

\section*{Declarations}
\subsection*{Data Availability Statement} 
No new data were generated or analyzed in this study. Data sharing is not applicable to this article.

\subsection*{Conflict of Interest} 
The author declares no competing interests.

\subsection*{Acknowledgments}
This work is supported by the China Postdoctoral Science Foundation (2026M793371) and the Postdoctoral Fellowship Program of CPSF (GZC20261681).

\subsection*{Declaration of generative AI and AI-assisted technologies}
The author used GPT-6 Astra to assist with organizing formulas and to obtain feedback through a simulated pre-submission review of the manuscript. The author reviewed and revised the AI-assisted content and takes full responsibility for the accuracy, integrity, and entire content of this article.

\end{document}